\documentclass{amsart}
\usepackage{graphicx}
\usepackage{stix}

\numberwithin{equation}{section}
\begin{document}

\title[Geometry on  Real Projective Cayley-Klein Spaces]{Geometry on Real Projective Cayley-Klein Spaces}
\author{Manfred Evers}
\curraddr[Manfred Evers]{Bendenkamp 21, 40880 Ratingen, Germany}
\email[Manfred Evers]{manfred\_evers@yahoo.com}
\date{\today}

\begin{abstract}
We investigate several topics of the geometry on real Cayley-Klein spaces.  An important concern for us is to define a distance function on the projective space in such a way that the distance between two anisotropic subspaces of the same dimension can be easily calculated and case distinctions are avoided as far as possible.
\end{abstract}

\maketitle \hspace*{\fill}\vspace*{-1 mm}\\

\section*{Introduction}
We introduce a real semi Cayley-Klein space as a pair $(\textrm{P}\boldsymbol{V},\beta)$ consisting of a projective space $\textrm{P}\boldsymbol{V}$ over a real vector space $\boldsymbol{V}$ and a symmetric bilinear form $\beta{:\,}\boldsymbol{V}\times\boldsymbol{V}\to\mathbb{R}$ different from the zero function. A point $P{\,\in\,}\textrm{P}\boldsymbol{V}$ is an isotropic point of $(\textrm{P}\boldsymbol{V},\beta)$ if for any nonzero vector  $\boldsymbol{p} \in \boldsymbol{V}$ representing $P$ (\,we write $P= \mathbb{R}\boldsymbol{p}\,$) the equation $\beta(\boldsymbol{p},\boldsymbol{p})=0$ applies, otherwise $P$ is anisotropic. Let $\mathcal{Q}$ be the set of isotropic points of  $(\textrm{P}\boldsymbol{V},\beta)$.

A reflection $\phi$ in an anisotropic point $P=\mathbb{R}\boldsymbol{p}$ of $(\textrm{P}\boldsymbol{V},\beta)$ is an automorphism of $\textrm{P}\boldsymbol{V}$ which is different from the identical mapping and maps any point $Q=\mathbb{R}\boldsymbol{q}$ of  $\textrm{P}\boldsymbol{V}$ to a point $Q'=\mathbb{R}\big(\beta(\boldsymbol{p},\boldsymbol{p})\boldsymbol{q}-2\beta(\boldsymbol{p},\boldsymbol{q})\boldsymbol{p}\big)$. These point reflections are idempotent mappings.
Points that stay fixed under all point reflections are called singular points in $(\textrm{P}\boldsymbol{V},\beta)$. All singular points are isotropic, but not all isotropic points have to be singular. The singular points form the radical of $\beta$, denoted by rad$(\beta)$. This is a projective subspace of $\textrm{P}\boldsymbol{V}$, which might be empty, but on the other hand it can also be a hyperplane in $\textrm{P}\boldsymbol{V}$. 

If $P=\mathbb{R}\boldsymbol{p}$ and $Q=\mathbb{R}\boldsymbol{q}$ are two anisotropic points in $(\textrm{P}\boldsymbol{V},\beta)$, then $\xi(P,Q):=$ $1{-}\beta(p,q)/(\beta(p,p)\beta(q,q ))$ is a real number which can be interpreted as a distance between the points $P$ and $Q$ and is called \textit{quadrance} of the points $P$ and $Q$.\footnote{$^)$ We adopt the definition of the function and the name \textit{quadrance} from N. Wildberger \cite{W1, W2}. }$^)$ If $\phi$ is an arbitrary reflection in an anisotropic point, then $\xi(P,Q)=\xi(\phi(P),\phi(Q))$ for all anisotropic points $P$ and $Q$.

The quadrance function $\xi{:\,}(\textrm{P}\boldsymbol{V}{\smallsetminus}{\hspace{-4.3pt}\smallsetminus}\,\mathcal{Q})\times(\textrm{P}\boldsymbol{V}{\smallsetminus}{\hspace{-4.3pt}\smallsetminus}\,\mathcal{Q}) \to \mathbb{R}$ can strongly degenerate. For example, if the set rad$(\beta)$ of singular points is a hyperplane in $\textrm{P}\boldsymbol{V}$, then $\xi$ is the zero function.  In this case, the distance function is obviously too coarse to allow meaningful distance measurement. But if we define a nonzero bilinear form $\beta'$ on rad$(\beta)$, we can significantly improve the situation: We can now introduce reflections in points which are anisotropic in the semi Cayley-Klein space $(\textrm{rad}(\beta),\beta')$. But more importantly, by the help of $\beta'$ the function $\xi$ can be modified to a continuous function $\xi'{:\,}(\textrm{P}\boldsymbol{V}{\smallsetminus}{\hspace{-4.3pt}\smallsetminus}\,\mathcal{Q})\times(\textrm{P}\boldsymbol{V}{\smallsetminus}{\hspace{-4.3pt}\smallsetminus}\,\mathcal{Q}) \to \mathbb{R}$ with\\
1. $\xi'(P,Q)\ne 0$ for almost all distinct anisotropic points $P,Q$,\\
2. $\xi'(P,Q) = \xi(P,Q)$ for all anisotropic points $P,Q \in \textrm{P}\boldsymbol{V}$ with $\xi(P,Q)\ne 0$, and \\
3. if $\phi$ is a reflection in an anisotropic point of $(\textrm{P}\boldsymbol{V},\beta)$, then $\xi'(P,Q)=\xi'(\phi(P),\phi(Q))$ \\
\hspace*{3.2 mm}for all $P, Q \in \textrm{P}\boldsymbol{V}{\smallsetminus}{\hspace{-4.3pt}\smallsetminus}\,\mathcal{Q}$. 

We now come to the notion of a real Cayley-Klein space. A {Cayley-Klein space} with underlying projective space $\textrm{P}\boldsymbol{V}$ is a nested sequence of finitely many, say $\rho+1$, semi Cayley-Klein spaces $(A_i,{{\beta_i}})$, such that $A_0 = \textrm{P}\boldsymbol{V}$, $A_{i+1}$ is a nonempty radical rad(${\beta_i}$), $0\le i\le \rho-1$, and $A_{\rho+1}{:=\textrm{rad}(\beta_\rho)}=\emptyset.$\vspace*{-1 mm}\\

The first section gives a short introduction to the subject. We assume that the reader is familiar with projective geometry, but in order to introduce the terminology and fix notations, we give some basic definitions, rules and theorems. As a main basis serves the book \textit{Vorlesungen \"uber h\"ohere Geometrie} \cite{Gi} by O. Giering, which was published in 1982 and deals with the topic of Cayley-Klein spaces in detail. There are a number of more recent works on this subject; we cite \cite{KH,Ko,OS,Ri,St1,St2}. An outline of the historical development is given in \cite{CP}.

In the second section, we study the geometry on semi Cayley-Klein spaces. This includes the determination of distances between projective subspaces. Another focus is on the investigation of geometric figures such as circles and simplices.

In the third section, the geometry on Cayley-Klein spaces is examined. In contrast to the calculations in the second section, distance measurement exhibits a fine structure.

In Section 4 it is shown how geometric algebra (GA) can be used to calculate reflections and their compositions. \vspace*{2 mm}  

\section{Fundamentals / terminology and notation\vspace*{0 mm}}
\subsection{} \textbf{Real and complex projective spaces.\,}
Let $\mathbb{F}$ be the field of real or complex numbers and $\boldsymbol{V}$ be a vector space of finite dimension $n{+}1$ over $\mathbb{F}$. We introduce the projective space of $\boldsymbol{V}$ by $\textrm{P}\boldsymbol{V} = \{\,\mathbb{F}\boldsymbol{v}\; |$ $\; \boldsymbol{v}\in \boldsymbol{V} {\smallsetminus}{\hspace{-4.3pt}\smallsetminus} \,\{\boldsymbol{0}\}\}$. The dimension $n$ of this projective space is defined by $n =\textrm{dim}(\boldsymbol{V}) - 1$. A subset ${U}$ of $\textrm{P}\boldsymbol{V}$ is called a \textit{subspace} if there exists a linear subspace $\boldsymbol{U}$ of $\boldsymbol{V}$ with ${U} =\textrm{P}\boldsymbol{U}$; we write ${U}\leq\textrm{P}\boldsymbol{V}$. In particular, the empty set is  a projective subspace of $\textrm{P}\boldsymbol{V}$. It is quite common to notate this subspace by $0$; the dimension of $0$ is $\textrm{-}1$. If ${U}_1 = \textrm{P}\boldsymbol{U}_1$ and ${U}_2 =\textrm{P}\boldsymbol{U}_2$ are subspaces of $\textrm{P}\boldsymbol{V}$, their \textit{join} is the projective subspace ${\textrm{P}(\boldsymbol{U}_1+\boldsymbol{U}_2)}$; we denote it by ${U}_1 \sqcup {U}_2$. The intersection of two subspaces ${U}_1, {U}_2$ of $\textrm{P}\boldsymbol{V}$ is always a subspace, which is called the \textit{meet} of ${U}_1$ and ${U}_2$ and which we denote by ${U}_1 \sqcap {U}_2$\footnote{$^)$ The symbols $\sqcap$ and $\sqcup$ are preferred to $\wedge$ and $\vee$,\\\hspace*{3.5 mm} since the latter are used for exterior products \big(in Grassmann algebras, see Section 4\big). }$^)$. The elements of $\textrm{P}\boldsymbol{V}$ are called \textit{points}. Let $\boldsymbol{B}:=(\boldsymbol{b}_1,\cdots,\boldsymbol{b}_{n+1})$ be an ordered basis of $\boldsymbol{V}$. Given a vector $\boldsymbol{v} \in \boldsymbol{V}$ with coordinates $({v}_1{,\dots\,}{v}_{n+1})$ with respect to this basis, we denote the vector $\boldsymbol{v}$ by $[v_1,\dots,v_{n+1}]_{\boldsymbol{B}}$  and the point $\mathbb{F}\boldsymbol{v}$ by ${[v_1:\dots:v_{n+1}]_{\boldsymbol{B}}}$.\\
If a subspace contains just one point $P$, then we  follow a convention and denote this subspace by $P$ instead of $\{P\}$. A $k$-dimensional subspace of $\textrm{P}\boldsymbol{V}$ is also called a $k$-plane (in $\textrm{P}V$). There are separate names for two special cases: A \textit{hyperplane} of $\textrm{P}\boldsymbol{V}$  is a subspace (a plane) of dimension $n{\,-\,}1$, and if $P$ and $Q$ are different points in $\textrm{P}\boldsymbol{V}$, the projective one-dimensional subspace $P \sqcup Q$ is called a \textit{line} passing through $P$ and $Q$.
If ${U}$ is a subspace of $\textrm{P}\boldsymbol{V}$, then a subspace ${U'}\le\textrm{P}\boldsymbol{V}$ with ${U}\sqcup {U}'= \textrm{P}\boldsymbol{V}$ and ${U}\sqcap {U}'=\emptyset$ is called a \textit{projective complement} of ${U}$.\\
Given a subset ${S}$ of $\textrm{P}\boldsymbol{V}\!$, the set $\textrm{span}({S}):= \bigcap \{{U}\,|\,{U}\;\textrm{is a subspace of}\;\textrm{P}\boldsymbol{V} \textrm{and}\; {{S} \subseteq {U}\,}\}$ is a projective subspace of $\textrm{P}\boldsymbol{V}$, called the \textit{span} of ${S}$. Put $\text{dim}(S) := \text{dim}(\text{span}(S))$. If $U\le \textrm{P}\boldsymbol{V}$ is a $k$-plane, a set of $k{+}1$ points spanning $U$ is called a \textit{minimal generating set} of $U$.
 We say that the points of a set ${S}\subseteq \textrm{P}\boldsymbol{V}\!$ are \textit{in a general position} if $\text{dim}(\tilde{\mathcal{S}}) = k-1$ for each subset $\tilde{{S}}$ of $\mathcal{S}$ with $\# \tilde{{S}} = k \le n{+}1$. An $(n{+}2)$-tuple $(P_0,\dots,P_{n+1})$ of points is called a \textit{projective frame} of $\textrm{P}\boldsymbol{V}$ if these points are in general position. We give an example: Suppose $\boldsymbol{V}=\mathbb{R}^{n+1}$ with canonical basis $\boldsymbol{B}=(\boldsymbol{e}_1,\dots,\boldsymbol{e}_{n+1})$. Put $E_i:= \mathbb{R}\boldsymbol{e}_i, 1\le i\le n{+}1$, and $E_{0}:=\mathbb{R}(\boldsymbol{e}_1+\cdots+\boldsymbol{e}_{n+1}).$ Then $(E_0,E_1,\dots,E_{n+1})$ is a projective frame. It is called \textit{the canonical frame} of  $\textrm{P}\mathbb{R}^{n+1}$. The point $E_{0} = [1:1:\dots:1]_{\boldsymbol{B}}$ is called \textit{unit-point} of the frame. Each $(n{+}1)$-elementary subset of $\{E_0,E_1,\dots,E_{n+1}\}$ is a minimal generating set  of $\textrm{P}\boldsymbol{V}$.\vspace*{1.5 mm}
\subsection{} \textbf{The cross ratio.\,}\;
Given four vectors $\boldsymbol{a},\boldsymbol{b},\boldsymbol{c},\boldsymbol{d}\in \mathbb{F}^2\,{\smallsetminus}{\hspace{-4.3pt}\smallsetminus} \,\{\boldsymbol{0}\}$, the cross ratio of the four points $\mathbb{F}\boldsymbol{a},\mathbb{F}\boldsymbol{b},\mathbb{F}\boldsymbol{c},\mathbb{F}\boldsymbol{d}\in \textrm{P}\mathbb{F}^2$ is given by the point\vspace*{0.7 mm}\\
$P=\hspace*{1 mm}(\mathbb{F}\boldsymbol{a},\mathbb{F}\boldsymbol{b};\mathbb{F}\boldsymbol{c},\mathbb{F}\boldsymbol{d})$\vspace*{0.7 mm}\\
$\hspace*{1.6mm}:=\Big(\det\!\left(\begin{array}{rr} 
a_1 & a_2  \\  
c_1 & c_2\end{array}\right)\det\left(\begin{array}{rr} 
d_1 & d_2  \\  
b_1 & b_2\end{array}\right){\,:\,}\det\left(\begin{array}{rr} 
a_1 & a_2  \\  
d_1 & d_2\end{array}\right)\det\left(\begin{array}{rr} 
c_1 & c_2  \\  
b_1 & b_2\end{array}\right)\,\Big){\;\in\;}\textrm{P}\mathbb{F}^2\,$,\vspace*{1.5 mm}\\
cf. \cite[Ch. 13]{Sto} for the case $\mathbb{F}=\mathbb{R}$. \vspace*{0.5 mm}\\
If $\varphi{:\,}\mathbb{F}^2\to \mathbb{F}^2$ is an automorphism, then\vspace*{0.5 mm}\\
\centerline{$(\mathbb{F}\boldsymbol{a},\mathbb{F}\boldsymbol{b};\mathbb{F}\boldsymbol{c},\mathbb{F}\boldsymbol{d})=(\mathbb{F}\varphi(\boldsymbol{a}),\mathbb{F}\varphi(\boldsymbol{b});\mathbb{F}\varphi(\boldsymbol{c}),\mathbb{F}\varphi(\boldsymbol{d}))$.}\vspace*{0.5 mm}\\
\textit{Remark}: {$\big(\det\left(\begin{array}{rr} 
a_1 & a_2  \\  
c_1 & c_2\end{array}\right)\det\left(\begin{array}{rr} 
d_1 & d_2  \\  
b_1 & b_2\end{array}\right)\big):\big({\det\left(\begin{array}{rr} 
a_1 & a_2  \\  
d_1 & d_2\end{array}\right)\det\left(\begin{array}{rr} 
c_1 & c_2  \\  
b_1 & b_2\end{array}\right)\big)}$}\vspace*{1 mm}\\
is usually interpreted as the ratio\vspace*{1 mm} of two numbers, but it should be noted that the divisor can be zero.\vspace*{1 mm}

\subsection{} \textbf{Projective collineations.\,}\;
Let $\varphi{:\,}\boldsymbol{V}\to \boldsymbol{V}'$ be an injective linear map between two $\mathbb{F}$-vector spaces $\boldsymbol{V}$ and $\boldsymbol{V}'$. Then $\varphi$ induces a map ${\phi}{:\,}\textrm{P}\boldsymbol{V}\to \textrm{P}\boldsymbol{V}'$ between the projective spaces $\textrm{P}\boldsymbol{V}$ and $\textrm{P}\boldsymbol{V}'$ by ${\phi}(\mathbb{F}\boldsymbol{v})=\mathbb{F}{\varphi(\boldsymbol{v})}$. Since ${\phi}{:\,}\textrm{P}V\to \textrm{P}V'$ maps  a set of collinear points of $\textrm{P}V$ to a set of collinear points of $\textrm{P}V'$,
${\phi}$ is called a \textit{projective collineation}.\\If $\varphi{:\,}V\to V'$ is bijective,  ${\phi}{:\,}\textrm{P}V\to \textrm{P}V'$ is an \textit{isomorphism}. \\
If $A,B,C,D$ are four points in a $1$-dimensional projective space $\textrm{P}\boldsymbol{V}$ and $\phi{:\,}\textrm{P}\boldsymbol{V}\to\textrm{P}\mathbb{F}^2$ is an isomorphism, then we put $(A,B;C,D):=(\phi(A),\phi(B);\phi(C),\phi(D))$. It can be easily checked that the point $(A,B;C,D)$ (the cross ratio of the four points) does not depend on the special choice of the isomorphism.\\
Four collinear points $A,B,C,D$ in a projective space \textit{form a harmonic range} if $(A,C;B,D)$ $=(-1:1).$

The automorphism group Aut$(\textrm{P}\boldsymbol{V})$ of $\textrm{P}V$ can be identified with the projective linear group PGL$(\boldsymbol{V})$ = GL$(\boldsymbol{V})/$Z, where the center\;Z\ of GL$(\boldsymbol{V})$ consists of all nonzero multiples of the identity. Given a basis $\boldsymbol{B}$ of $\textrm{P}\boldsymbol{V}$, an element $\phi$ of $\textrm{Aut}(\textrm{P}\boldsymbol{V})$ can be represented by an invertible $(n{+}1)\times(n{+}1)$-matrix $\mathfrak{A}_\phi{\,=\,}(\mathfrak{a}_{ij})_{1\le i,j\le n{+}1}$ which acts on a row $(x_1,\dots,x_{n+1})$ of $n{+}1$ entries by matrix multiplication from the right: \\ $\phi([x_1:\dots:x_{n+1}]_{\boldsymbol{B}})=[y_1:\dots:y_{n+1}]_{\boldsymbol{B}}\,,\, (x_1,\dots,x_{n+1})\mathfrak{A}_\phi{\,=\,}(y_1,\dots,y_{n+1})$. \\
All these automorphisms are collineations. Furthermore, these mappings preserve the cross ratio of quadruples of collinear points. Two subsets of $\textrm{P}\boldsymbol{V}$ are called \textit{projectively equivalent} if one can be mapped onto the other by an automorphism.\\
Given two projective frames $(P_0,\dots,P_{n+1})$, $(Q_0,\dots,Q_{n+1})$, there exists precisely one automorphism $\phi\in \text{Aut}(\textrm{P}\boldsymbol{V})$ with $\phi(P_i)=Q_i$, $0\le i\le n{+}1$. \\
The only automorphism that fixes a projective frame is the identity; any other automorphism has at most $n{+}1$ independent fixed points. An automorphism $\phi\ne \textrm{id}$ is called \textit{biaxial collineation with axes} $A$ \textit{and} $B$ if $A$ and $B$ are nonempty complementary subspaces of $\textrm{P}\boldsymbol{V}$ and $\phi$ fixes all points of  $A$ and of $B$. Such a collineation also fixes all lines incident with $A$ and $B$. If $A$ is a hyperplane and $B$ consists of one point, $\phi$ is called \textit{dilation}  with axis $A$ and center $B$. A dilation $\phi$ is a special case of a \textit{central collineation} (also called \textit{perspectivity}); all the hyperplanes through one point, its center, are invariant under $\phi$. There is another kind of central collineation which is called \textit{elation}. As a dilation, an elation also leaves a hyperplane, its axis, pointwise invariant and fixes all hyperplanes through its center; but the center of an elation is a point on its axis.\\
Given any central collineation
$\phi$ with center $Z=\mathbb{R}\boldsymbol{z}$ and axis $U=\mathbb{R}\boldsymbol{U}$, we can find a linear form $m\!:\boldsymbol{V}\to \mathbb{R}$ with $\boldsymbol{U}=\textrm{ker}(m)$ and $m(\boldsymbol{z})\ne -1$ such that\vspace*{-0.0 mm}
\\
\centerline{$\phi(\mathbb{R}\boldsymbol{p}) =  \mathbb{R}(\boldsymbol{p}+m(\boldsymbol{p}){_{\,}}\boldsymbol{z})\;$for all $\boldsymbol{p}\in \boldsymbol{V}$.}\vspace*{0.0 mm}\\
The points $P{\,=\,}\mathbb{R}\boldsymbol{p}, Z, \phi(P)$ are collinear. And $\phi(P) = P$ precisely when $P = Z$ or $m(\boldsymbol{p}) = 0$.\\
Aut$(\textrm{P}\boldsymbol{V})$ is generated by central collineations; each automorphism is the product of at most $n{+}1$ central collineations. 

We denote the set of all projective subspaces of $\textrm{P}\boldsymbol{V}\!$ by $\textrm{sub}(\textrm{P}\boldsymbol{V})$. An automorphism $\phi\in \text{Aut}(\text{P}\boldsymbol{V})$ induces a mapping\,  $\text{sub}(\text{P}\boldsymbol{V})\to \text{sub}(\text{P}\boldsymbol{V})$, which we also denote by $\phi$, with $\phi({U}_1\sqcup {U}_2) = \phi({U}_1)\sqcup \phi({U}_2)$ and 
$\phi({U}_1\sqcap {U}_2) = \phi({U}_1)\sqcap \phi({U}_2)$.\\
A mapping ${\kappa{:}\textrm{sub}(\textrm{P}\boldsymbol{V})\!\to}$ $ \textrm{sub}(\textrm{P}\boldsymbol{V})$ is called a \textit{correlation} if it maps subspaces of dimension $k$ to subspaces of dimension $n{-}k{-}1$ such that $\kappa({U}_1\sqcup {U}_2)=\kappa({U}_1)\sqcap \kappa({U}_2)$ and $\kappa({U}_1\sqcap {U}_2)=\kappa({U}_1)\sqcup \kappa({U}_2)$.
\subsection{} \textbf{Quadrics.\,}\;
We already introduced $\boldsymbol{V}$ as an $\mathbb{F}$-vector space of dimension $n{+}1$. Let $\beta{:\,}\boldsymbol{V}\times \boldsymbol{V}\to \mathbb{F}$ be a symmetric bilinear form.  Given a basis $\boldsymbol{B} = (\boldsymbol{b}_1,\dots,\boldsymbol{b}_{n+1})$, we define a matrix 
$\mathfrak{A} = (\mathfrak{a}_{ij})\in \mathbb{F}^{n+1}\times \mathbb{F}^{n+1} $ by $\mathfrak{a}_{ij}{\,:=\,}\beta(\boldsymbol{b}_i,\boldsymbol{b}_j)$. The rank of this matrix does not depend on the special choice of the basis $\boldsymbol{B}$,
so we can assign this number to the bilinear form $\beta$. It is always possible to find an orthogonal basis $(\boldsymbol{b}_1,\dots,\boldsymbol{b}_{n+1})$ for $\beta$, i.e. a basis of $\boldsymbol{V}$ with $\beta(\boldsymbol{b}_i,\boldsymbol{b}_j)=0\; \text{for}\,  i\ne j.$ If $\boldsymbol{B}$ is an orthogonal basis, the associated matrix is diagonal. In the following we always assume that rank$(\beta)=\textrm{rank}(\mathfrak{A})$ is not zero.\\
The set $Q_\beta:=\{ P{\,=\,\mathbb{F}\boldsymbol{b}}\in\text{P}\boldsymbol{V}\,| \,\beta(\boldsymbol{b},\boldsymbol{b})=0\}$ is called the \textit{quadric} associated with $\beta$. In the case of $\mathbb{F}=\mathbb{R}$ this quadric can be the empty set; but if $\mathbb{F}=\mathbb{C}$, this quadric is always nonempty.\\
Two points $P=\mathbb{F}\boldsymbol{p}, Q=\mathbb{F}\boldsymbol{q}\in \textrm{P}\boldsymbol{V}$ are called \textit{conjugate} with respect to $\beta$, if $\beta(\boldsymbol{p},\boldsymbol{q})=0$. Two sets $\mathcal{S}_1, \mathcal{S}_2 \subset \textrm{P}\boldsymbol{V}$ are \textit{conjugate} with respect to $\beta$, if $ P_1 \,\textrm{is conjugate to}\, P_2$ for all points $P_1{\,\in\,}\mathcal{S}_1$ and $P_2{\,\in\,}\mathcal{S}_2$. If  $\mathcal{S}$ is a subset of $\textrm{P}\boldsymbol{V}$, the set $\textrm{polar}_{\beta}(\mathcal{S}):= \{P\in \textrm{P}\boldsymbol{V}\,$ $|\,P \,\textrm{is conjugate to}\,\mathcal{S} \}$ is a projective subspace of $\textrm{P}\boldsymbol{V}$ and is called \textit{the} \textit{polar} of $\mathcal{S}$ with respect to $\beta$. Obviously, $\mathcal{S}\subseteq \textrm{polar}_{\beta}(\textrm{polar}_{\beta}(\mathcal{S}))$. The mapping $\textrm{polar}_\beta{: \textrm{sub}(\textrm{P}\boldsymbol{V})\!\to}$ $ \textrm{sub}(\textrm{P}\boldsymbol{V})$ is a correlation iff rank$(\beta)= n{+}1$.\vspace*{-3 mm}\\

\subsection{} \textbf{Real semi Cayley-Klein spaces.}\;We are especially interested in the case $\mathbb{F}=\mathbb{R}$. We choose a symmetric bilinear form $\beta{:\,}\boldsymbol{V}\times \boldsymbol{V}\to\mathbb{R}$ with rank$(\beta) > 0$ and call it the \textit{absolute bilinear form} of $\boldsymbol{V}$. If $\boldsymbol{B}=(\boldsymbol{b}_1,\dots,\boldsymbol{b}_{n+1})$ is an orthogonal basis for $\beta$, we put $\beta_+{:=\,}\#\{\boldsymbol{b}_i|\beta(\boldsymbol{b}_i,\boldsymbol{b}_i){>}0\}$, \,$\beta_-{:=\,}\#\{\boldsymbol{b}_i |\, \beta(\boldsymbol{b}_i,\boldsymbol{b}_i) < 0\}$,\; $\beta_0{:=\,}\#\{\boldsymbol{b}_i |\, \beta(\boldsymbol{b}_i,\boldsymbol{b}_i) = 0\}$. These three numbers do not change if we change the basis for another orthogonal basis (Sylvester's law of inertia). The triple $(\beta_+,\beta_-,\beta_0)$ is called the \textit{signature} of $\beta$. 
Now it is obvious that we can even find an orthogonal basis $\boldsymbol{B}=(\boldsymbol{b}_1,\dots,\boldsymbol{b}_{n+1})$ such that
$\beta(\boldsymbol{b}_i,\boldsymbol{b}_i)=1$ for all $i{\,\in\,}\{1,\dots,\beta_+\}$, $\beta(\boldsymbol{b}_i,\boldsymbol{b}_i)=-1$ for all $i{\,\in\,} \{\beta_{+}{+\,}1,\dots,\beta_+{+\,}\beta_-\}$ and $\beta(\boldsymbol{b}_i,\boldsymbol{b}_i)=0$ for all $i> \textrm{rank}(\beta)$.\\
$\mathrm{P}\boldsymbol{V}$ can now be decomposed into three disjoint subsets: $\textrm{P}\boldsymbol{V} = Q_\beta{\,\dot\cup\,}Q_{\beta}^{\,+}{\,\dot\cup\,}Q_{\beta}^{\,-}$ where 
$Q_\beta$ is the quadric associated with $\beta$, $Q_\beta^{\,+}{:=\,}\{\mathbb{R}\boldsymbol{v}\,|\,\beta(v,v)>0\}$ and $Q_\beta^{\,-}:=\{\mathbb{R}\boldsymbol{v}\,|\,\beta(v,v)<0\}$. Since $\beta{:\,}\boldsymbol{V}\times \boldsymbol{V}\to\mathbb{R}$ is a nonzero function, $\mathcal{Q}_{\beta}^{\;\times}:=\mathcal{Q}_{\beta}^{\;+}\cup \mathcal{Q}_{\beta}^{\;-}$ is a dense subset of $\textrm{P}\boldsymbol{V}$ (with respect to the topology on $\textrm{P}\boldsymbol{V}$, which is the quotient topology of the canonical topology on $\boldsymbol{V}$).  We call $(\textrm{P}\boldsymbol{V},\beta)$ a \textit{semi Cayley-Klein space}.\vspace*{0.7 mm}

Let $\boldsymbol{U}$ be a linear subspace of $\boldsymbol{V}$. We assume that the restriction ${\beta'}:= \beta|\boldsymbol{U}\times \boldsymbol{U}$ is a nonzero function. Then $(\textrm{P}{\boldsymbol{U}}, {\beta}')$ is a \textit{semi CK subspace} of $(\textrm{P}\boldsymbol{V}\!,\beta)$.\vspace*{+0.5 mm}

Let $U=\textrm{P}\boldsymbol{U}\le \textrm{P}\boldsymbol{V}$ be a plane of dimension $s\le n$ and let $\boldsymbol{u}_1,\dots,\boldsymbol{u}_{s+1}$ be a basis of $\boldsymbol{U}$. Define a matrix $\mathfrak{N}=(\mathfrak{n}_{ij})$ by $\mathfrak{n}_{ij}:= \beta(\boldsymbol{u}_i,\boldsymbol{u}_j)$.
The rank of this matrix $\mathfrak{N}$ does not dependent of the special choice of the basis. If $\det(\mathfrak{N})= 0$, we call $U$ \textit{isotropic}, otherwise  \textit{anisotropic}. 
As a special case, a point $P\in \mathbb{R}\boldsymbol{V}$ is isotropic iff $P\in\mathcal{Q}_{\beta}$.\\ If $P$ is an anisotropic point, then $P\cap \textrm{polar}_\beta(P)=\emptyset$.\\$P$ is called \textit{double point} \textit{of} \textit{a plane} $U$ if $U\subseteq \textrm{polar}_\beta(P)$. \vspace*{+0.9 mm}\\
A plane $U$ is {anisotropic} iff it does not contain any double point.\vspace*{+0.6 mm}\\ 
\textit{Proof}: We assign a matrix $\mathfrak{N}$ to $U$ as described above.\\ 
If there is a double point $P$ in $U$, we may assume that it is $\mathbb{R}\boldsymbol{u}_1$. In this case, the first row of the matrix $\mathfrak{N}$ is a zero-row and $\det(\mathfrak{N})=0$.\\ Now we assume that $U$ is isotropic and $\det(\mathfrak{N})=0$. We may also assume that the matrix $\mathfrak{N}$ is a diagonal matrix. Then there exists a vector $\boldsymbol{u}_i$ with $\beta(\boldsymbol{u}_i,\boldsymbol{u}_j)=0$ for all $j$, and $\mathbb{R}\boldsymbol{u}_i$ must be a double point of $U$.\;$\Box$\vspace*{+0.6 mm} \\
A point $P{\,=\,}\mathbb{R}(\boldsymbol{p})\in \textrm{P}\boldsymbol{V}$ is called \textit{singular} if it is a double point of $\textrm{P}\boldsymbol{V}$. The singular points form a plane in $\textrm{P}\boldsymbol{V}$ which is called the \textit{vertex} of $\mathcal{Q}_{\beta}$ or the \textit{radical} rad$(\beta)$ of $\beta$. \\
\textit{Remark}: A point $P$ is a double point of the $0$-dimensional space $P$ iff $P \in \mathcal{Q}_{\beta}$.\\
Properties of polars: (1) The polar of a singular point is $\textrm{P}\boldsymbol{V}$, the polar of a nonsingular point is a hyperplane in $\textrm{P}\boldsymbol{V}$. (2) $A\le \textrm{P}\boldsymbol{V}$ is an anisotropic plane iff $\textrm{polar}_\beta(A){:=\,}\bigcap_{P\in A}\textrm{polar}_\beta(P)$ is a projective complement of $A$ in $\textrm{P}\boldsymbol{V}$.\vspace*{0.6 mm}\\
A subspace ${U}\le\textrm{P}\boldsymbol{V}$ of positive dimension located entirely in the quadric $\mathcal{Q}_{\beta}$ is called \textit{totally isotropic}. It is obvious, that all subspaces of $\textrm{rad}(\beta)$ are totally isotropic if their dimension is positive. But even if rad$(\beta)$ is empty, $(\textrm{P}\boldsymbol{V},\beta)$ can have totally isotropic subspaces. An example is $(\,\textrm{P}\mathbb{R}^{4},\beta\,)$ with $\beta\big((v_1,v_2,v_3,v_4),(w_1,w_2,w_3,w_4)\big)=v_1w_1+v_2w_2-v_3w_3-v_4w_4.$\\
The maximal dimension of a totally isotropic subspace is $n -\text{max}(\beta_+,\beta_-)$.\vspace*{0 mm}
\subsection{} \textbf{Automorphisms on a semi CK space $(\,\mathrm{P}\boldsymbol{V}, \beta\,)$.}\;
Let $(\textrm{P}\boldsymbol{V}\!,\beta),(\textrm{P}\tilde{\boldsymbol{V}}\!,\tilde{\beta})$ semi CK spaces. Let $\phi{:\,}\textrm{P}\boldsymbol{V}\to\textrm{P}\tilde{\boldsymbol{V}}$ be a projective collineation and $\varphi{:\;}\boldsymbol{V}\to \tilde{\boldsymbol{V}}$ a representative of $\phi$. 
Then $\phi$ is called a \textit{projective collineation} from the semi CK space $(\textrm{P}\boldsymbol{V}\!,\beta)$ to the semi CK space $(\textrm{P}\boldsymbol{\tilde{\boldsymbol{V}}}\!,\tilde{\beta})$ if  $\tilde{\beta}({\varphi}(\boldsymbol{p}),{\varphi}(\boldsymbol{q})){\,=\,}0 \Leftrightarrow\beta(\boldsymbol{p},\boldsymbol{q}){\,=\,}0$ for all points $P=\mathbb{R}\boldsymbol{p}$ and $Q=\mathbb{R}\boldsymbol{p}$.
If this mapping $\phi$ is bijective, then $\phi$ is an isomorphism  between \textit{semi CK spaces} and, in the case of $(\textrm{P}\boldsymbol{V}\!,\beta) =$ $(\textrm{P}\boldsymbol{\tilde{\boldsymbol{V}}}\!,\tilde{\beta})$, $\phi$ is an automorphism on $(\textrm{P}\boldsymbol{V}\!,\beta)$.

Special automorphisms on a semi CK space are reflections: 
Let $Z= \mathbb{R}\boldsymbol{z}\in \textrm{P}\boldsymbol{V}$ be an anisotropic point and let $\phi\in \textrm{Aut}(\textrm{P}\boldsymbol{V})$ be a central collineation with center $Z$ and axis $\textrm{polar}_\beta(Z)$. Then $\textrm{polar}_\beta(Z)$ is a hyperplane not incident with $Z$, and $\phi$ is a \textit{dilation}.
Because rad$(\beta) \le \textrm{polar}_\beta(Z)$, all points of rad$(\beta)$ are fixed points. There exists some number $\lambda\in \mathbb{R}{\smallsetminus}{\hspace{-4.3pt}\smallsetminus}\{0,1\}$ such that $\phi$ maps a point $P=\mathbb{R}\boldsymbol{p}\in \textrm{P}\boldsymbol{V}$ to a point\vspace*{0.5 mm}\\
\centerline{$P' = \mathbb{R}\big(\lambda \beta(\boldsymbol{z},\boldsymbol{z})\boldsymbol{p}+(1-\lambda)\beta(\boldsymbol{z},\boldsymbol{p})\boldsymbol{z}\big) $.}\vspace*{0.8 mm}\\
If $\lambda = -1$, $\phi$ is an involution and $\phi$ is called a \textit{reflection in the point} $Z$.
It can be easily checked that a reflection maps an anisotropic point to an anisotropic point lying in the same connected component of $\mathcal{Q}_{\beta}^{\;\times}$. If $P$ is the only isotropic point on the (isotropic) line $P\sqcup Z$, then $P$ is a fixed point of the reflection $\phi$, and if the line $P\sqcup Z$ contains still another isotropic point $Q$, then $P$ and $Q$ are interchanged by $\phi$.\\
Given an anisotropic plane  $U$, then a biaxial collineation $\phi$ with axes $U$ and $\textrm{polar}_\beta(U)$ is called a \textit{reflection in} $U$. If $Z$ is an anisotropic point of $U$ and $P$ is any point of $\textrm{P}\boldsymbol{V}$, then $\phi(P)$ is the image of $P$ under a reflection in $Z$.
\subsection{} \textbf{Barycentric coordinates of points in a semi CK space.}\;
Let from now on $\boldsymbol{V}$ be the real vector space $\mathbb{R}^{n+1}$ with canonical basis $\boldsymbol{B}=(\boldsymbol{e}_1,\dots,\boldsymbol{e}_{n+1})$. Given a vector $\boldsymbol{v}=[v_1,\dots,v_{n+1}]_{\boldsymbol{B}} \in \boldsymbol{V}$, we write $\boldsymbol{v}=(v_1,\dots,v_{n+1})$, and for the point $\mathbb{R}\boldsymbol{v}\in \textrm{P}\boldsymbol{V}$ we use the notation $(v_1{\,:\,}\dots{\,:\,}v_{n+1})$. Let $(\beta_+,\beta_-,\beta_0)$ be a triple of nonnegative integers with $\beta_+{\,+\,}\beta_-{\,+\,}\beta_0=n{+}1$
and let $\mathfrak{A}^{(\beta_+,\beta_-,\beta_0)}=\mathfrak{A}= (\mathfrak{a}_{ij})$ denote the $(n{+}1)\times (n{+}1)$-diagonal matrix with $\mathfrak{a}_{ii}{\,=\,}1$ for $1\le i\le \beta_+$, $\mathfrak{a}_{ii}{\,=\,}{-}1$ for $\beta_+{\,<\,}i\le \beta_+{+\,}\beta_-$ and $\mathfrak{a}_{ii}{\,=\,}0$ for $i{\,>\,}\beta_+{+\,}\beta_-$. Associated with $\mathfrak{A}$ is the bilinear form $\beta{:\,}\mathbb{R}^{n+1}{\times\,}\mathbb{R}^{n+1}\to \mathbb{R}$,\, $\beta(\boldsymbol{e}_i,\boldsymbol{e}_j)=\mathfrak{a}_{ij}$.\\$(\textrm{P}\boldsymbol{V},\beta)$ is a semi CK space.\vspace*{+0.6 mm}\\
\textit{Remark}: Every real semi CK space can be subsumed here (up to isomorphism).\vspace*{+0.6 mm}\\
Instead of $\beta(\boldsymbol{v},\boldsymbol{w})$ we will usually write $\boldsymbol{v}{\scriptstyle{[\mathfrak{A}]}}\boldsymbol{w}$. (We regard ${\scriptstyle{[\mathfrak{A}]}}$ as a generalized inner product.)
For $\mathcal{Q}_{\beta}$ we write $\mathcal{Q}_{\mathfrak{A}}$, the polar of a set $\mathcal{S}\subset \mathrm{P}\boldsymbol{V}$ with respect $\beta$ to will be denoted by $\mathcal{S}^{\mathfrak{A}}$ instead of polar$_\beta(\mathcal{S})$, and the semi CK space $(\textrm{P}\boldsymbol{V},\beta)$ we denote by $(\textrm{P}\boldsymbol{V},\mathfrak{A})$.\vspace*{1.5 mm}

We assign a vector $P^\circ\in \boldsymbol{V}$  to each anisotropic point $P$: 
First, we define a function $\chi{:}\;\boldsymbol{V} \rightarrow \{-1,0,1\}$ by\\ \vspace*{2.5 mm}  
{\!$\chi(p_0{,\dots,}p_{n+1}) = 
\begin{cases}
 \;\,0,\;\text{if } (p_1,\dots,p_{n+1}) = (0,\dots,0)\;, \\
 \;\,1,\;\text{if } (p_1,\dots,p_{n+1}) > (0,\dots,0)\; \textrm{with respect to the lexicographic order,} \\
-1,\text{if } (p_1,\dots,p_{n+1}) < (0,\dots,0)\; \textrm{with respect to the lexicographic order}, \\
\end{cases}$}\vspace*{-1.5 mm}\\
then we put $\displaystyle P^\circ := \frac{\chi(\boldsymbol{p})}{\sqrt{|{\boldsymbol{p} {\scriptstyle{[\mathfrak{A}]}\,}\boldsymbol{p}|}}}\,\boldsymbol{p}$.\vspace*{+1.0 mm}\\
We introduce a function sgn${:\,}\textrm{P}\boldsymbol{V}\to \{-1,0,1\}$ by sgn$(P)=0$
if $P$ is isotropic, and \\
sgn$(P)=P^\circ{\scriptstyle{[\mathfrak{A}]}\,}P^\circ$ if $P$ is anisotropic. \vspace*{0.9 mm} \\ 
\hspace*{1.6 mm}Since the matrix $\mathfrak{A}$ is nonzero, there exist $\,n{+}1$ independent anisotropic points $P_1,\dots,P_{n+1}$ generating $\textrm{P}\boldsymbol{V}$. Given a point $Q=\mathbb{R}\boldsymbol{q}$, there exists an $(n{+}1)$-tuple $(\tilde{q}_1,\dots,\tilde{q}_{n+1})$ of real numbers such that $\boldsymbol{q} = \tilde{q}_1 P_1^\circ+\cdots+\tilde{q}_{n+1}P_{n{+}1}^\circ$. The tuple $(\tilde{q}_1,\dots,\tilde{q}_{n+1})$ is uniquely determined by the point $Q$ and the tuple $(P_1,\dots,P_{n+1})$ except for multiplication by a nonzero real number, and is called \textit{tuple of homogeneous coordinates} of $Q$ with respect to $(P_1,\dots,P_{n+1})$. We write $Q=\sum_{1\le i\le n+1}\tilde{q}_iP_i$ and, after having fixed the tuple $(P_1,\dots,P_{n+1})$, we also write $Q = [\tilde{q}_1:\dots:\tilde{q}_{n+1}]$.\\  
If $Q$ is anisotropic, there exists a uniquely determined $(n{+}1)$-tuple $(\tilde{q}_1,\dots,\tilde{q}_{n+1})$ such that $Q^\circ = \tilde{q}_1 P_1^\circ+\cdots+\tilde{q}_{n+1}P_{n{+}1}^\circ$. This tuple is the tuple of \textit{absolute coordinates} of $Q$ with respect to $(P_1,\dots,P_{n+1})$.
\subsection{}\textbf{The dual of a semi CK space.}\;Let $\boldsymbol{V}^\ast = \textrm{Hom}(\boldsymbol{V},\mathbb{R})$ be the vector space dual to $\boldsymbol{V}$ and let $\boldsymbol{e}^1,\dots,\boldsymbol{e}^{n+1}$ be the dual basis of $\boldsymbol{e}_1{,\dots,}\boldsymbol{e}_{n+1}$. Put $E^i{:=\,}\mathbb{R}\boldsymbol{e}^i, i{\,=\,}1,\dots,n{+}1$. The projective space $\textrm{P}\boldsymbol{V}^\ast$ is called the \textit{dual} of $\textrm{P}\boldsymbol{V}$.\vspace*{0.5 mm}

Let $(\textrm{P}\boldsymbol{V},\mathfrak{A})$ be a semi CK space, and let $\mathfrak{A}^\sharp = \textrm{adj}(\mathfrak{A})$ denote the adjugate of $\mathfrak{A}$. We call the semi CK space $(\textrm{P}\boldsymbol{V}^\ast,{{\mathfrak{A^{\sharp}}}})$ the \textit{dual }of $(\textrm{P}\boldsymbol{V},{{\mathfrak{A}}})$.\vspace*{0.5 mm} 

We now assume that rank$(\mathfrak{A})=n{+}1.$ If $P_1, P_2$ are anisotropic points in $\textrm{P}\boldsymbol{V}$, then $H_1{\,:=\,}P_1^\mathfrak{A}, H_2{\,:=\,}P_2^\mathfrak{A}$ are anisotropic hyperplanes in $\textrm{P}\boldsymbol{V}$ and there are uniquely determined anisotropic elements $x_1, x_2\in \textrm{P}\boldsymbol{V}^\ast$ with $x_1(Q_1)=x_2(Q_2)=0$ for all $Q_1\in H_1$ and all $Q_2\in H_2$. For $x_1, x_2$ and $P_1, P_2$ the equation $\displaystyle\frac{(x_1{\scriptstyle{[\mathfrak{A}^{-1}]}}x_2)^2}{(x_1{\scriptstyle{[\mathfrak{A}^{-1}]}}x_1)(x_2{\scriptstyle{[\mathfrak{A}^{-1}]}}x_2)}{\,=\,}\frac{(P_1{\scriptstyle{[\mathfrak{A}]}}P_2)^2}{(P_1{\scriptstyle{[\mathfrak{A}]}}P_1)(P_2{\scriptstyle{[\mathfrak{A}]}}P_2)}$  applies.
\subsection{} \textbf{Perspectivity and orthology.}\;
Let $U\le \textrm{P}\boldsymbol{V}$ be a plane of dimension $s>1$. 
$U$ is arguesian: Let $(P_1,\dots,P_{s{+}1})$, $(Q_1,\dots,Q_{s{+}1})$ be two systems of points, each generating $U$.  $(P_1,\dots,P_{s{+}1})$ and $(Q_1,\dots,Q_{s{+}1})$ are \textit{perspective} if there exists exactly one point $Z$ such that all the sets $P_i\sqcup Q_i\sqcup Z\,$, $1\le{i\le{s{+}1},}$ are lines. If such a point $Z$ exists, it is called the \textit{perspector} of $(P_1,\dots,P_{s{+}1})$ and $(Q_1,\dots,Q_{s{+}1})$.\\
Suppose this perspector $Z$ is different from all the points $P_i, Q_i,\; i= 1,\dots,s{+}1 $, then there exists a uniquely determined central collineation  $\phi{:\,}U{\,\to\,}U$ mapping $(P_1,\dots,P_{s{+}1},Z)$ onto $(Q_1,\dots,Q_{s{+}1},Z)$. The axis of $\phi$, a hyperplane of $U$, is called \textit{perspectrix} of $(P_1,\dots,P_{s{+}1})$ and $(Q_1,\dots,Q_{s{+}1})$.\\ A perspectrix can be assigned to these two tuples even if the perspector is one of the points $P_i, Q_i, 1\le i\le s{+}1.$ 
For example, if $Z=P_{s+1}\ne Q_{s+1}$, then $Q_1\sqcup\dots\sqcup Q_{s}$ is the perspectrix.\vspace*{-2 mm}\\

We now assume that $U$ is anisotropic. Given an $(s{+}1)$-tuple $R_1,\dots,R_{s+1}$ of points spanning $U$, put $R_i^{\;\delta}{:=\,}(\bigsqcup_{j\ne i}R_j)^\mathfrak{A}\,\sqcap\,U$, $1\le i\le s{+}1$.
If $(P_1,\dots,P_{s{+}1})$ and $(Q_1,\dots,Q_{s{+}1})$ are perspective with center $Z$ and perspectrix $S\le U$, then $(P_1^{\;\delta},\dots,P_{s+1}^{\;\,\delta}) , (Q_1^{\;\delta},\dots,Q_{s+1}^{\;\,\delta})$ are also perspective; their perspector is $S^\mathfrak{A}\sqcap U$ and their perspectrix is $Z^\mathfrak{A} \sqcap A$.\vspace*{-0.0 mm}\\
Two tuples $(P_1,\dots,P_{s+1})$ and $(Q_1,\dots,Q_{s+1})$ are called \textit{orthologic} with \textit{ortho\-logy center} $Z$ if $(P_1,\dots,P_{k+1})$ and $(Q_1^{\;\delta},\dots,Q_{k+1}^{\;\,\delta})$ are perspective at $Z$.\vspace*{0.0 mm}\\
If $(P_1{,\dots,}P_{s+1})$ and $(Q_1{,\dots,}Q_{s+1})$ are orthologic, then $(Q_1{,\dots,}Q_{s+1})$ and $(\!P_1{,\dots,}P_{s+1})$ are also orthologic. If $S$ is the perspectrix of $(\!P_1{,\dots,}P_{s+1})$ and $(Q_1^{\;\delta}{,\dots,}Q_{s+1}^{\;\,\delta})$, then $S^\mathfrak{A}\sqcap A$ is the associated orthology center.\vspace*{1.9 mm}
  
\section{Metric structures on $\mathrm{P}V$\vspace*{1.7 mm}}
\subsection{} \textbf{Metric structures on the projective line.}\; 
K. von Staudt \cite{Stau}, E. Laguerre \cite{La} and A. Cayley \cite{Ca} showed us how quadrics in $\textrm{P}\boldsymbol{V}$ can be used to introduce metric structures on $\textrm{P}\boldsymbol{V}$. 
An important role within the connection between quadrics and distances plays the cross ratio of points and of lines. Later, systematic studies by Felix Klein \cite{Kl} led to a classification of metric geometries. We give a rough illustration of the ideas of these great mathematicians for the simple case $n = 1$. 
We start with the inner product ${\scriptstyle{[\mathfrak{A}]}}\,$ on the projective line $\textrm{P}\mathbb{R}^2$, which is given by the matrix  $\mathfrak{A} = \left(\begin{array}{rr} 
1 & 0  \\ 
0 & r \\ 
\end{array}\right)$. (Every quadric on the projective line is projectively equivalent to a quadric $\mathcal{Q}_{\mathfrak{A}}$ with a suitable number $r \in \{-1,0,1\}$.) The set of isotropic points is empty if $r>0$, and consists of the two points $J_-=(-\sqrt{-r}:1)$ and $J_+=(\sqrt{-r}:1)$ if $r<0$. 
If $r=0$, there is exactly one isotropic point. This point is a double point of $\textrm{P}\mathbb{R}^2$; the points  $J_-$ and $J_+$ merge, so to speak,   for $r=0$ to one point $J = (0:1)$.
In the case of $r>0$, Laguerre \cite{La} encouraged us to think of $\textrm{P}\mathbb{R}^2$ embedded in the complex projective space $\textrm{P}\mathbb{C}^2$ and to think of the points $J_-$ and $J_+$ as complex valued because these two complex valued points (in fact, his focus was on the angle between lines rather than on the distance between points) are helpful for introducing a metric.
The following classification comes from F. Klein: depending on whether $r> 0$ or $r = 0$ or $r< 0$, we speak of the \textit{elliptic} or \textit{parabolic} or \textit{hyperbolic case}.\vspace*{-2 mm}\\

Let us first assume $r\ne 0$. In this case $(\textrm{P}\mathbb{R}^2,\mathfrak{A})$ is an anisotropic CK space.
We take two independent anisotropic points $P_1=\mathbb{R}\boldsymbol{p}_1=\mathbb{R}(p_{11}{\,,\,}p_{12})$ and $P_2=\mathbb{R}\boldsymbol{p}_2=\mathbb{R}(p_{21} , p_{22})$, $\chi(p_{11} , p_{12})=\chi(p_{21} , p_{22})$, and determine homogeneous barycentric coordinates for $J_-$ and $J_+$ with respect to $(P_1,P_2)$:\vspace*{1 mm}\\
\centerline{$\displaystyle J_\mp = \mathbb{R}\big({(p_{21}\pm\sqrt{-r}\, p_{22})} (p_{11},p_{12}) - {(p_{11}\pm\sqrt{-r}\, p_{12})}(p_{21},p_{22})\big)$}\vspace*{1 mm}
{$\hspace*{23,5 mm}\displaystyle = {(p_{21}\pm\sqrt{-r}\, p_{22})}{\sqrt{p_{11}^{\;2}+r\,p_{12}^{\;2}}}\,P_1-{(p_{11}\pm\sqrt{-r}\, p_{12})}{\sqrt{p_{21}^{\;2}+r\,p_{22}^{\;2}}}\,P_2.$}\vspace*{1 mm}\\
The cross ratio of the four points $P_1,P_2,J_-,J_+$ is \vspace*{1 mm}\\
\centerline{$\displaystyle (P_1,P_2;J_-,J_+) = {((p_{11}-\sqrt{-r}\,p_{12})(p_{21}+\sqrt{-r}\,p_{22}))}{\,:\,}{((p_{11}+\sqrt{-r}\,p_{12})(p_{21}-\sqrt{-r}\,p_{22}))}\,.$}\\

\noindent In the elliptic case, this cross ratio is a complex number of modulus $1$ and, following Laguerre, Cayley and Klein, we can define a distance $d_{12} \in \mathbb{R}$ between the points $P_1$ and $P_2$ by $d_{12} = c\,\ln(\sqrt{(P_1,P_2;J_-,J_+)}$ with a constant $c$ that can be fixed to $c = \frac{1}{2\mathbf{i}}$ by requiring that the distance $d$ of two points on the elliptic line is within the range $0 \le d\le \frac{1}{2} \pi$ and that the distance of these points is $\frac{1}{2} \pi$ precisely when one point is conjugate to the other. If we take this distance function, we get the equation (see \cite[ch. 5]{Ko}\vspace*{0.4 mm}\\ 
\centerline{$\cos(d_{12}) = P_1^{\,\circ} {\scriptstyle{[\mathfrak{A}]}\,} P_2^{\,\circ}\;\;\;\;\;\;(*).  $}\vspace*{0.4 mm}  \\ 
At the beginning of the $20^{\,th}$ century, the Hungarian mathematician C. V\"or\"os\,\footnote{$^)$ I take this information from \cite{H1}.}$^)$ {proposed to proceed in the hyperbolic case quite analogously to the elliptic case: All hyperbolic lines have the same length  $\pi\mathbf{i}$, and the distance function satisfies the equation \vspace*{1mm}\\
\centerline{$\displaystyle \cosh^2(d_{12}) =\, \frac{(\boldsymbol{p}_1 {\scriptstyle{[\mathfrak{A}]}\,} \boldsymbol{p}_2)^2}{({\boldsymbol{p}_1 {\scriptstyle{[\mathfrak{A}]}\,} \boldsymbol{p}_1)}({\boldsymbol{p}_2 {\scriptstyle{[\mathfrak{A}]}\,} \boldsymbol{p}_2})}\, = \textrm{sgn}({P_1})\,\textrm{sgn}({P_2})\,(P_1^{\,\circ} {\scriptstyle{[\mathfrak{A}]}\,} P_2^{\,\circ})^2\,\;\;\;(**).  $}  \vspace*{0mm}\\ 
If $\textrm{sgn}({P_1})\textrm{sgn}({P_2}) = 1$, then $\cosh(d_{12}) > 1$ and $d_{12}$ is a real number. If $\textrm{sgn}({P_1})\textrm{sgn}({P_2}) = -1$, then $\cosh(d_{12})$ is purely imaginary. In this case, equation $(**)$ can be satisfied by allowing $d_{12}$ to be a complex number of the form $t{\,+\,}\frac{1}{2}{\pi}\mathbf{i}\,,\,t{\,\in\,}\mathbb{R}$.\\
We adopt V\"or\"os' proposal and, in order to achieve as much harmony as possible in formulae that are valid in the elliptic plane and corresponding ones that are valid in the hyperbolic plane, we replace the elliptic distance $d_{12}$ by $d_{12}\mathbf{i}\,.$
Now formula $(**)$ is valid also on the elliptic line.\\
The points conjugate to $P_1^{\,}$ and $P_2^{\,}$ are $P_1^{\,\mathfrak{A}}=(p_{12}{:}r\,p_{11})$ and $P_2^{\,\mathfrak{A}}=(p_{22}{:}r\,p_{21})$, respectively. The cross ratio $(P_1^{\,\mathfrak{A}},P_2^{\,\mathfrak{A}};J_-,J_+)$ is the same as the cross ratio $(P_1^{\,},P_2^{\,};J_-,J_+)$. \\
We calculate the cross ratio of the points $P_1^{\,},P_2^{\,\mathfrak{A}}, P_2^{\,},P_1^{\,\mathfrak{A}}$ and get\vspace*{0.7 mm}\\
\centerline{$\displaystyle (P_1^{\,},P_2^{\,\mathfrak{A}};P_2^{\,},P_1^{\,\mathfrak{A}}) = \,\frac{r\,(p_{11}p_{22}-p_{12} p_{21})^2}{(p_{11}^{\;2}+r\,p_{12}^{\;2})(p_{21}^{\;2}+r\,p_{22}^{\;2})} = - \sinh^2(d_{12})\;\;\;\; (***)\,.$}\vspace*{1.2 mm}\\
We adopt from N.J. Wildberger \cite{W2}  the name \textit{quadrance} of $P_1$ and $P_2$ for the real number\vspace*{1 mm}\\ \centerline{$\displaystyle\xi_0(P_1,P_2):=- \sinh^2(d_{12})={1 - \displaystyle \textrm{sgn}({P_1})\,\textrm{sgn}({P_2})\,{(P_1^{\,\circ} {\scriptstyle{[\mathfrak{A}]}\,} P_2^{\,\circ})^2}}$}\\
{$\hspace*{57.6mm}\displaystyle =\, \frac{\det\!\left(\begin{array}{rr} 
\boldsymbol{p}_1{\,\scriptstyle{[\mathfrak{A}]}\,}\boldsymbol{p}_1 & \boldsymbol{p}_1{\,\scriptstyle{[\mathfrak{A}]}\,}\boldsymbol{p}_2  \\ 
\boldsymbol{p}_2{\,\scriptstyle{[\mathfrak{A}]}\,}\boldsymbol{p}_1 & \boldsymbol{p}_2{\,\scriptstyle{[\mathfrak{A}]}\,}\boldsymbol{p}_2 \\ 
\end{array}\right)}{(\boldsymbol{p}_1{\,\scriptstyle{[\mathfrak{A}]}\,}\boldsymbol{p}_1)(\boldsymbol{p}_2{\,\scriptstyle{[\mathfrak{A}]}\,}\boldsymbol{p}_2)}$\;.}  \vspace*{1.8 mm}\\ 
If $Q = x P_1 + y{\,}P_2, x,y\in\mathbb{R}$, is anisotropic, then (see \cite{Ev3} for a proof)\vspace*{1 mm}\\
\centerline{$x^2\,:\,y^2 = ({\textrm{sgn}({P_2})}\sinh^2(d(Q,P_1))){\;:\,}({\textrm{sgn}({P_1}})\sinh^2(d(Q,P_2)))\;$.}\vspace*{-0.4 mm}\\

Let us consider now the parabolic case. If $r=0$, there is precisely one isotropic point, the double point $J = (0{\,:\,}1)$. Given two distinct anisotropic points $P_1, P_2 \in \textrm{P}\mathbb{R}^2$, we can find real numbers $p_{12}, p_{22}$ such that $P_1 = (1{\,:\,}p_{12})$ and $P_2 = (1{\,:\,}p_{22})$. If we now apply the formula $\cosh(d(P_1,P_2)) = P_1^{\,\circ} {\scriptstyle{[\mathfrak{A}]}\,} P_2^{\,\circ}$, we find out that in the parabolic case the projective line is a null line: The distance $d(P_1,P_2)$
and the quadrance $\xi_0(P_1,P_2)$ of any two anisotropic points $P_1,P_2$ is 0. We can still postulate that the total length of a projective line is  $\pi \mathbf{i}$, we just have to demand that the length of a line segment is $0$ iff the double point $J$ lies outside this segment.  
We can accept this "gross structure", the distance of two anisotropic points on a parabolic line being always 0, but it is possible to find a "fine structure"  by using possibilities that Nonstandard Analysis offers us.\vspace*{1.5mm}\\
\textit{Excursus}: Let $\mathbb{R}^\ast$ and $\mathbb{C}^\ast$ denote the hyperreal and hypercomplex numbers. Any function $f{:\,}\mathbb{C}\to \mathbb{C}$ can be extended to a function $f^\ast{:\,}\mathbb{C}^\ast\to \mathbb{C}^\ast$ in a canonical way (s. \cite{Str}). In the following write $f$ instead of $f^\ast$; we skip $\ast$ .\\
If $f$ is meromorphic in a neighborhood of $0$ in $\mathbb{C}$, we can find a series $(a_n)_{n{\,=\,}0,\dots,\infty}\in \mathbb{C}^\mathbb{N}$ and an integer $m$, both uniquely determined, such that $a_0\ne 0$ and $f(\varepsilon)=\varepsilon^{m}(a_0+a_1\varepsilon+a_2\varepsilon^2+\cdots)$ for all infinitesimals $\varepsilon\ne 0$. Let $\star_{\varepsilon}(f(\varepsilon)):=\varepsilon^{m} a_0$ be the leading term of $f(\varepsilon)$.\vspace*{-1.8mm} \\

Now we define an $\mathbf{\varepsilon}$\textit{-quadrance} of anisotropic points $P_1=(p_{11}{:}p_{12})$ and $P_2=(p_{21}{:}p_{22})$\vspace*{1.5mm}\\
by\centerline{\ $\displaystyle \xi_{\mathbf{\varepsilon}}(P_1,P_2) := \star_{\mathbf{\varepsilon}}\big(1- \frac{(P_1^{\,\circ} {\scriptstyle{[\mathfrak{A}_\mathbf{\varepsilon}]}\,} P_2^{\,\circ})^2}{({P_1^{\,\circ} {\scriptstyle{[\mathfrak{A}_\mathbf{\varepsilon}]}\,} P_1^{\,\circ})}{(P_2^{\,\circ} {\scriptstyle{[\mathfrak{A}_\mathbf{\varepsilon}}]}\, P_2^{\,\circ})}}\big)$,}\vspace*{1mm}\\
where $\mathbf{\varepsilon}$ is a nonzero infinitesimal real number and $\displaystyle \mathfrak{A}_\mathbf{\varepsilon}:=\textrm{diag}(1,\mathbf{\varepsilon})$.\\
While the two functions $\xi_{\mathbf{\varepsilon}}$ and $\xi_0$ agree on the elliptic and on the hyperbolic line, they differ on the parabolic line.\vspace*{0.5mm}\\
On the parabolic line the following applies:\\
Let $P_1$ and $P_2$ be two distinct anisotropic points, then\vspace*{0.7mm}\\ 
\centerline{$\displaystyle \xi_{\mathbf{\varepsilon}}(P_1,P_2) = \mathbf{\varepsilon}\,(\frac{p_{21}}{p_{11}}-\frac{p_{22}}{p_{21}})^2\ne 0 = \xi_0(P_1,P_2)\,.$}\vspace*{0.5mm}\\
And if $Q = (q_1{\,:\,}q_2) \in \textrm{P}\mathbb{R}^2$ is anisotropic and $Q = x P_1{\,+\,}y{\,}P_2, {x,y\in\mathbb{R}}$, then\vspace*{0.5mm}\\
\centerline{  ${|x|\,:\,|y|} = |q_1 p_{22}-q_2 p_{21}| : |q_1 p_{12}-q_2 p_{11}|\,$.}\vspace*{0.5mm}\\
Let $Q_1, Q_2$ be two anisotropic points. We can find real numbers $s_{11},s_{12},s_{21},s_{22}$ such that $Q_1^{\;\circ} =
s_{11} P_1^{\;\circ}{\,+\,}s_{12} P_2^{\;\circ}$ and $Q_2^{\;\circ} = s_{21} P_1^{\;\circ}{\,+\,}s_{22} P_2^{\;\circ}$. Then,\vspace*{0.5mm}\\
\centerline{$\xi_{\mathbf{\varepsilon}}(Q_1,Q_2) = \xi_{\mathbf{\varepsilon}}(P_1,P_2)(s_{11}-s_{21})(s_{12}-s_{22})\,$.}\vspace*{1mm} \\
\textit{Remarks:} (1) In the following we often use the name \textit{squared distance} for $\xi_{\mathbf{\varepsilon}}(P_1,P_2)$ if $P_1,P_2$ are anisotropic points in a parabolic plane. (2) Usually, the (squared) distance of anisotropic points in a parabolic space (in a euclidean space or a Lorentz-Minkowski space, for example) is given by a real number $d$ instead of an infinitesimal number $d{\,\mathbf{\varepsilon}}$. (3) E. Study \cite[\S 23]{Stu} was presumably the first to use an infinitesimal number (an infinitesimal dual number ${\,\mathbf{\varepsilon}}$) for calculating distances between points in a euclidean space. Considerations of this kind are already presented in Klein's paper \cite[p. 612]{Kl}.
	


\subsection{} \textbf{Point reflections on a one-dimensional semi CK space $(\mathrm{P}\mathbb{R}^2, \mathfrak{A}\,).$}\; 
The reflection $\sigma_Z$ in an anisotropic point $Z\in \textrm{P}\mathbb{R}^2 $ fixes the points $Z$ and $Z^\mathfrak{A}$. Every other point $P$ is mapped to a point $\tilde{P}$ such that the points $Z, P, Z^\mathfrak{A},\tilde{P}$  form a harmonic range. It can be easily checked that $\textrm{sgn}(P)=\textrm{sgn}(\tilde{P})$. If $P$ is anisotropic, then $\tilde{P}$ is also anisotropic and $d(\tilde{P},Z) = d(P,Z)$.
The last statement applies for the elliptic and hyperbolic lines, but also for the parabolic line if we regard it as a null line.
The double point on a parabolic line is a fixed point of $\sigma_Z$; the two isotropic points on a hyperbolic line are interchanged by $\sigma_Z$. \\
If we now use a parabolic quadrance on the parabolic line given by $\xi_{\mathbf{\varepsilon}}(P,Q) = \mathbf{\varepsilon}\,(p-q)^2$ for anisotropic points $P=(1:p)$ and $Q=(1:q)$, we can write $Z=(1:z)$, $P'=(1:p')$ and get $(1:p')=P'=2Z-P=(1:2z{-}p)$. It follows $\xi_{\mathbf{\varepsilon}}(P',Z) = \xi_{\mathbf{\varepsilon}}(P,Z)$; the parabolic ${\mathbf{\varepsilon}}$-quadrance of two points is preserved by reflections.\vspace*{+1.5 mm} 

In the following we will first use a gross distance measure. In Section 3 we turn to geometries where parabolic lines represent a generic and not just a singular case; this is when fine distance measurement comes into play.


\subsection{} \textbf{Gross measurement in a semi CK space.}\;
Let $(\textrm{P}\boldsymbol{V}, \mathfrak{A}\,)$ be a semi CK space. 
A line without any isotropic point is an elliptic line. A line is hyperbolic iff it 
passes through anisotropic points of different sign; in this case it contains precisely two nonsingular isotropic points. A line that is neither hyperbolic nor elliptic, must either be totally isotropic or it must pass through exactly one isotropic point, which in this case must be a double point. A line with just one isotropic point is a parabolic line.\\
We introduce line segments with anisotropic endpoints and the length of these segments.
Lines are compact sets, such as all subspaces of $\textrm{P}V$. We give all lines, even those  entirely located in $\mathcal{Q}_{\mathfrak{A}}$, the same measure $\pi \mathbf{i}$, with imaginary unit $\mathbf{i}$. Given a line ${L}\le\textrm{P}\mathbf{V}$ and an anisotropic point $P$ on ${L}$, then there is exactly one point $Q$, denoted by $\textrm{conj}(P,{L})$, which is conjugate to $P$; it is the meet of ${L}$ and $P^\mathfrak{A}$. $Q$ is an anisotropic point, and $P = \textrm{conj}(Q,{L})$. Given two distinct anisotropic points $P$ and $Q$, we introduce two line segments $[P,Q]_+ := \{s P{\,+\,}t\,Q |\, s,t{\,\in\,}\mathbb{R}, st \geq 0\}$ and $[P,Q]_- := \{s P{\,+\,}t\,Q |\, s,t{\,\in\,}\mathbb{R}, st \leq 0\}$. Thus, $[P,Q]_+$ and $[P,Q]_-$ are the closures of the two connected components of the set $P \sqcup Q\; {\smallsetminus}{\hspace{-4.3pt}\smallsetminus}\,  \{P, Q\}$.\\
We define lengths $\mu_0([P,Q]_\pm)$ of these two segments as complex numbers with imaginary parts in the interval $[\,0,\pi\,]$, determined by the following conditions:\vspace*{0.7 mm}\\
$\bullet\;\displaystyle \cosh^2(\mu_0([P,Q]_+) = {\textrm{sign}({P})\,\textrm{sign}({Q})\,(P^\circ{{\scriptstyle{[\mathfrak{A}]}\,}} Q^\circ)^2}$ and $\mu_0(\,[P,Q]_-) = \pi \mathbf{i} - \mu(\,[P,Q]_+)$.\vspace*{1 mm}\\
$\bullet$ If $R$ is an anisotropic inner point of $[P,Q]_+$, then $\mu_0([P,R]_+) + \mu_0([R,Q]_+) = \mu_0([P,Q]_+).$\vspace*{0.6 mm}\\
\noindent\hspace*{0mm}$\bullet$ $\mu_0(\,[P,Q]_+) = \mu_0(\,[P,Q]_-) = \frac{1}{2}\pi \mathbf{i}$ precisely when $P$ and $Q$ are mutually conjugate.\vspace*{1 mm}\\
\noindent\hspace*{0mm}$\bullet$ If the line $P \sqcup Q$ is elliptic, then $(P^\circ{\scriptstyle{[\mathfrak{A}]}}\,Q^\circ)^2 < 1$ and $\mu_0(\,[P,Q]_{+})$, $\mu_0(\,[P,Q]_{-})$ are purely {\;\,imaginary}.  \vspace*{1 mm}\\
\noindent\hspace*{0mm}$\bullet$ If $\xi_0(P,Q)=0$, then there is exactly one isotropic point $R$ on $P{\,\sqcup\,}Q\,$.\\ If $R^{\,}\in [P,Q]_+$, then $\mu_0([P,Q]_+)$ $= \pi \mathbf{i}$; otherwise, $\mu([P,Q]_+) = 0$.\vspace*{1 mm}\\
\noindent\hspace*{0mm}$\bullet$ If the line $P \sqcup Q$ is hyperbolic, then it intersects $\mathcal{Q}_{\mathfrak{A}}$ transversally in two isotropic points. Here, we consider two cases:\\
Case 1: The points $P$ and $Q$ are in the same connected component of $\mathcal{Q}_\mathfrak{A}^{\;\times}$. Both isotropic points on $P \sqcup Q$ are either in $[P,Q]_+$ or in $[P,Q]_-$. Let $I$ be that of these two intervals without any isotropic point. Then $\mu_0(I)$ is one of two real numbers with the same absolute value which satisfy the equation $\cosh^2(\mu_0(I))= (P^\circ{\scriptstyle{[\mathfrak{A}]}}\,Q^\circ)^2$. Having the choice between a positive or a negative number for $\mu_0(I)$, we decide for the negative number if $P$ and $Q$ are points in $\mathcal{Q}_\mathfrak{A}^{\;-}$, and for the positive number if $P$ and $Q$ are in $\mathcal{Q}_\mathfrak{A}^{\;+}$. The length of the other interval is $\pi\,\mathbf{i} - \mu_0(I)$. \\
Case 2: 
The points $P$ and $Q$ are in different connected components of $\mathcal{Q}_\mathfrak{A}^{\;\times}$. Let us assume that ${P}\in \mathcal{Q}_\mathfrak{A}^{\;+}$ and $Q\in\mathcal{Q}_\mathfrak{A}^{\;-}$. Then  $R\!:= \text{conj}(P,P{\,\sqcup\,}Q)$ $\in \mathcal{Q}_\mathfrak{A}^{\;-}$ , and $S\!:= \text{conj}(Q,P{\,\sqcup\,}Q)\in \mathcal{Q}_\mathfrak{A}^{\;+}.$  In this case, \\
\noindent\hspace*{2mm}$\mu_0(\,[P,Q]_+) = \,\mu_0(\,[P,S]_+) + \frac{1}{2}\pi \mathbf{i}, \; \mu_0(\,[P,Q]_-) = \,-\mu_0(\,[P,S]_+) + \frac{1}{2}\pi \mathbf{i}$\, if $S \in [P,Q]_+$, and\\
\noindent\hspace*{2mm}$\mu_0(\,[P,Q]_+) = -\mu_0(\,[P,S]_+) + \frac{1}{2}\pi \mathbf{i}, \;\mu_0(\,[P,Q]_-) = \mu_0(\,[R,Q]_+) + \frac{1}{2}\pi \mathbf{i}$\, if $S \in [P,Q]_-$. \vspace*{0 mm}\vspace*{1 mm}\\ 
An analysis of the different cases shows that by knowing the number $P^\circ{{\,\scriptscriptstyle{[\mathfrak{A}]}\,}}\,Q^\circ$ and one of the two numbers ${\textrm{sgn}({P})}$, ${\textrm{sgn}({Q})}$  we can determine $\mu_0([P,Q]_\pm)$. \vspace*{1.5 mm}

We now define the lengths of the segments with one or two isotropic boundary points, guided by the principle: Isotropic boundary points are distributed equally to adjacent segments. \\
Given an anisotropic point and a double point on a parabolic line, then the length of each of the two segments with these points as boundary points is $\frac{1}{2}\pi\mathbf{i}$.\\
If $P$ is an anisotropic and $Q$ an isotropic point on a hyperbolic line, then the line splits into two segments with these points as boundary points. One contains still another isotropic point while the other does not. The length of the first segment is $ -\textrm{sgn}({P}) \infty + \frac{3}{4}\pi\mathbf{i}$, the length of the second $\,\textrm{sgn}({P})\,\infty + \frac{1}{4}\pi\mathbf{i}$.\\
A hyperbolic line ${L}$ is a disjoint union of the three sets  ${L}^0:= {L} \cap \mathcal{Q}_\mathfrak{A}, {L}^\pm:= {L} \cap \mathcal{Q}_\mathfrak{A}^{\;\pm}$. The length of the open segment  ${L}^\pm$ is $\pm\infty$, the length of the closed segment ${L}^\pm\cup {L}_0$ is $\;\pm\infty+\frac{1}{2}\pi\mathbf{i}$.\vspace*{0.5 mm}\vspace*{-2.6 mm}\\

In order to define a distance of anisotropic points, we introduce an order $\prec$ on the complex numbers by\\ 
\noindent\hspace*{0mm}$a_1 + b_1\mathbf{i} \prec a_2 + b_2\mathbf{i} \;\;\;\text{iff} \;\;\;
\begin{cases}
&b_1 < b_2 \vspace*{-2 mm}\\
\text{or}\vspace*{-2 mm}&\\
&b_1 = b_2 \;\;\text{and }\, a_1 < a_2\;.
\end{cases}$\vspace*{1.5 mm}\\
Put $\mathbb{D} := \{a{\,+\,}b{\,}\mathbf{i}\, |\, a \in {\mathbb{R}}, b \in [0,\frac{1}{2}\pi] \}$\\
and define the function $d_0{:\,}{Q}_{\mathfrak{A}}^\times\times\mathcal{Q}_{\mathfrak{A}}^\times\rightarrow \mathbb{D}$ by\vspace*{0.5 mm}\\
\centerline{$d_0(P,Q) = 
\begin{cases}
0,\;\;\text{if } P = Q\;, \\
\mu_0([P,Q]_+),\;\text{if }\, P \ne Q \;\,\textrm{and}\;\, \mu_0([P,Q]_+) \prec  \mu_0([P,Q]_-)\;, \\
\mu_0([P,Q]_-),\;\text{otherwise}\,.
\end{cases}$}\vspace*{1 mm}\\
This function $d_0$ is continuous with respect to the canonical topology on $\text{P}V$, and we call $d_0(P,Q)$ the \textit{gross distance} between the (anisotropic) points $P$ and $Q$.\vspace*{1 mm}\\
\textit{Remarks}: (1)\; On an elliptic line, the distance of two points is uniquely determined by its value under $\textrm{cosh}.$ The situation is different for a hyperbolic line $P{\,\sqcup\,}Q$; in order to determine the distance $d_0(P,Q)$, we also need to know one of the numbers $\textrm{sgn}({P}),\textrm{sgn}({Q})$.\\
(2)\, This distance function can be extended to isotropic points, as long as both points are not points on a totally isotropic line. But this extension is not continuous.\vspace*{1 mm}\\
(3) So far, no distance has been defined yet between points on the total vertex ${U}$ of $\mathcal{Q}$. This can be achieved by declaring a quadric in ${U}$ as the absolute quadric, such that $U$ becomes a CK space, see Section 3.
\subsection{} \textbf{The quadrance of two equidimensional anisotropic planes.}\;
Let $A_1$ and $A_2$ be two anisotropic planes of the same dimension $s{\,\le\,}n$ generated by points $P_1=\mathbb{R}\boldsymbol{p}_1{,\dots,}P_{s+1}=\mathbb{R}\boldsymbol{p}_{s+1}$ and $Q_1=\mathbb{R}\boldsymbol{q}_1{,\dots,}Q_{s+1}=\mathbb{R}\boldsymbol{q}_{s+1}$, respectively. We define two real numbers\\
\centerline{$\displaystyle \zeta_0(A_1,A_2){\;:=\;}\frac{(\det(\boldsymbol{p}_i{\scriptstyle{[\mathfrak{A}]}}\boldsymbol{q}_i)_{1\le i,j\le s+1})^2}{(\det(\boldsymbol{p}_i{\scriptstyle{[\mathfrak{A}]}}\boldsymbol{p}_j)_{1\le i,j\le s+1})(\det(\boldsymbol{q}_i{\scriptstyle{[\mathfrak{A}]}}\boldsymbol{q}_j)_{1\le i,j\le s+1})}\;\;\;\;(\ddagger)$}\vspace*{0.5 mm}\\

\noindent and ${\hspace{10.8 mm}}\xi_0(A_1, A_2){\;:=\;}1-\zeta_0(A_1, A_2)$.}\vspace*{0.7 mm}\\
\noindent If $A_1, A_2$ are points, $\xi_0(A_1, A_2)$ is the quadrance of the two points. We adopt this name for higher dimensional planes of the same dimension.\\

\subsection{} \textbf{Reflections on a semi CK space are isometries.}\\ Let $\phi\in \textrm{Aut}(\textrm{P}\boldsymbol{V})$ be a reflection in an anisotropic plane $U\le\textrm{P}\boldsymbol{V}$ and let $A_1$, $A_2 \le \textrm{P}\boldsymbol{V}$ be two anisotropic planes of the same dimension, then $\xi_0(\phi(A_1), \phi(A_2))=\xi_0(A_1, A_2)$.\\
\textit{Proof}. First, we give a proof for the special case $\textrm{dim}(A_1) = \textrm{dim}(A_2) = 0$.  In this case, $\phi(A_1)$ and $\phi(A_2)$ are anisotropic points. Put $M{:=\,}U\sqcap(A_1\sqcup \phi(A_1))$, $\tilde{M}{:=\,}U^{\mathfrak{A}}\sqcap(A_1\sqcup \phi(A_1))$, $N{:=\,}U\sqcap(A_2\sqcup \phi(A_2))$, $\tilde{N}{:=\,}U^{\mathfrak{A}}\sqcap(A_2\sqcup \phi(A_2))$. Only one of the following four cases can occur:\\
\noindent\hspace*{1 mm}(1) $M=A_1+\phi(A_1), \tilde{M}=A_1-\phi(A_1), N=A_2+\phi(A_2), \tilde{N}=A_2-\phi(A_2)$,\\
\noindent\hspace*{1 mm}(2) $M=A_1-\phi(A_1), \tilde{M}=A_1+\phi(A_1), N=A_2+\phi(A_2), \tilde{N}=A_2-\phi(A_2)$,\\
\noindent\hspace*{1 mm}(3) $M=A_1+\phi(A_1), \tilde{M}=A_1-\phi(A_1), N=A_2-\phi(A_2), \tilde{N}=A_2+\phi(A_2)$,\\
\noindent\hspace*{1 mm}(4) $M=A_1-\phi(A_1), \tilde{M}=A_1+\phi(A_1), N=A_2-\phi(A_2), \tilde{N}=A_2+\phi(A_2)$.\\
These cases can be treated essentially the same, therefore we only look at the first case: We have $0 = (A_1^{\,\circ}-\phi(A_1)^{\circ}){\scriptstyle{[\mathfrak{A}]}}(A_2^{\,\circ}+\phi(A_2)^{\circ})$ and $0= (A_1^{\,\circ}+\phi(A_1)^{\circ}){\scriptstyle{[\mathfrak{A}]}}(A_2^{\circ}-\phi(A_2)^{\circ})$ and therefore $\phi(A_1)^{\circ}{\scriptstyle{[\mathfrak{A}]}} A_2^{\,\circ}=\phi(A_2)^{\circ}{\scriptstyle{[\mathfrak{A}]}}A_1^{\,\circ}$. From this follows that $\xi(A_1,\phi(A_2))=\xi(A_2,\phi(A_1))$.
By interchanging the roles of $A_2$ and $\phi(A_2)$ we get $\xi(\phi(A_1),\phi(A_2))=\xi(A_1,A_2).$\\
We now assume that $A_1$ is generated by independent anisotropic points $P_1,\dots,P_s$ and $A_2$ is generated by independent anisotropic points $Q_1,\dots,Q_s$. Then $\xi_0(X,Y)=\xi_0(\phi(X),\phi(Y))$ for all $X,Y{\in\,}\{P_1,\dots P_s,Q_1,\dots,Q_s\}$ and 
 $\xi_0(A_1,A_2)=\xi_0(\phi(A_1),\phi(A_2)).$$\,\Box$\\

\subsection{} \textbf{Angles between planes.}\;
We define: Two planes $A_1, A_2\le \textrm{P}\boldsymbol{V}$ \textit{form dihedral angles} if $\textrm{dim}(A_1\sqcap A_2)=\textrm{dim}(A_1)-1=\textrm{dim}(A_2)-1$. \\
Let $A_1, A_2$ be two planes which form dihedral angles. Put $H:=A_1\sqcup A_2$ and $U:=A_1\sqcap A_2$. $H{\smallsetminus}{\hspace{-4.3pt}\smallsetminus}\; (A_1 \cup A_2)$ consists of two connected components. We call the closure of these components (\textit{dihedral}) \textit{angles}. $A_1$ and $A_2$ are the \textit{sides} of these angles, and $U$ is their \textit{vertex}.\\
\textit{Example}: Two distinct points form dihedral angles with vertex $0$.\vspace*{0.6 mm}\\
If the vertex $U$ is nonempty and anisotropic, $U^{\mathfrak{A}} \sqcap H$ is a line which meets $A_1 {\smallsetminus}{\hspace{-4.3pt}\smallsetminus}\, U$ and $A_2 {\smallsetminus}{\hspace{-4.3pt}\smallsetminus}\,U$ in points $P_1$ and $P_2$, respectively.
The segment $[P_1;P_2]_+$ belongs to one angle, the segment $[P_1,P_2]_-$ to the other. We notate the first angle by $[A_1,A_2]_+$, the second by $[A_1,A_2]_-$. If $P_1$ and $P_2$ are anisotropic, then $\xi_0(A_1,A_2) = \xi_0(P_1,P_2)$. Thus, we put $d_0(A_1,A_2) := d_0(P_1,P_2)$.\\
If the vertex is isotropic while the two sides $A_1, A_2$ are anisotropic, then $\xi_0(A_1,A_2)=0$ and $d_0(A_1,A_2)=0$.\vspace*{0.5 mm}\\
\textit{Two examples}:\\ (1) If $\boldsymbol{V}=\mathbb{R}^4,\mathfrak{A}=\textrm{diag}(1,1,-1,-1),P_1=(1{:}1{:}0{:}0),P_2=(1{:}{-}1{:}1{:}0),$ $P_3=(0{:}0{:}1{:}0)$, then
$P_1^\circ{\scriptstyle{[\mathfrak{A}]}}P_2^{\circ}=0=P_1^{\circ}{\scriptstyle{[\mathfrak{A}]}}P_3^{\circ}$ and $\zeta_0(P_1\sqcup P_2,P_1\sqcup P_3)=-1$, $\mu([P_1\sqcup P_2,P_1\sqcup P_3]_+)=\mu([P_2,P_3]_+)=\ln(\sqrt{2}+1)+\frac{1}{2}\pi\mathbf{i},$ $d_0(P_1\sqcup P_2,P_1\sqcup P_3)=\mu([P_2,P_3]_-)=\ln(\sqrt{2}-1)+\frac{1}{2}\pi\mathbf{i}$. \\
(2) If $\boldsymbol{V}=\mathbb{R}^3, \mathfrak{A}=\textrm{diag}(1,1,0), P_1=(1{:}1{:}1), P_2=(-1{:}1{:}1),P_3=(1{:}-1{:}1)$, then $\mu([P_1\sqcup P_2,P_1\sqcup P_3]_+)=\mu([P_2,P_3]_+)=\pi\mathbf{i}$ and $\xi_0(P_1\sqcup P_2,P_1\sqcup P_3)=\xi_0(P_2,P_3)=0$.
\vspace*{1.2 mm}\\
The distance $d_0(A_1,A_2)$ between almost all planes $A_1,A_2$ forming dihedral angles is nonzero as long as their dimension is less than $\textrm{rank}{(\mathfrak{A})}-1$. If their dimension exceeds $\textrm{rank}{(\mathfrak{A})}-2$, their distance is always 0.

\subsection{} \textbf{Projections and parallels.}\;
Consider some point $P = \mathbb{R}\boldsymbol{p}$ and some anisotropic plane ${U} = \mathbb{R}\boldsymbol{U}$ of dimension $k<n$ with $P \notin {U}^\mathfrak{A}$. 
The \textit{perpendicular} from $P$ to ${U}$ is the $(n{-}k)$-plane $\,\text{perp}(P,{U}) := P \sqcup {U}^\mathfrak{A}$. This plane meets ${U}$ at the point $\text{ped}(P,{U}) := \text{perp}(P,{U})\sqcap {U}$, the \textsl{projection} or \textit{pedal} of $P$ in ${U}$. If $P\in {U}$, then $\text{ped}(P,{U}) = P.$
Regardless of the choice of point $P$, at least one of the points $\textrm{ped}(P,{U}),\textrm{ped}(P,{U}^\mathfrak{A})$ is well-defined. If $P$ is a point neither on ${U}$ nor on ${U}^\mathfrak{A}$, then both are well-defined and together with $P$ they form a collinear triad.\\
If $P \in {U}^\mathfrak{A}$, then $d_0(P,{U}) = \frac{1}{2}\pi\,\mathbf{i}$, otherwise  $d_0(P,{U}) = d_0(P;\textrm{ped}(P,U))$.\\
If ${U}$ is a hyperplane with ${U}^\mathfrak{A}{=\,}Q$ and $P\ne Q$, then $\text{ped}(P,{U})  
 = ({Q^\circ}{\scriptstyle{[\mathfrak{A}]}}{P^\circ}){Q}{\,-\,}(Q^\circ{\scriptstyle{[\mathfrak{A}]}}Q^\circ){P}$.\vspace*{1 mm}\\
Suppose ${A}$ is a second plane with $\textrm{dim}\,{A}\le \textrm{dim}\,{U}$ and ${A} \sqcap {U}^\mathfrak{A} = \emptyset$. Then ped$({A},{U}) := \bigcup_{P\in {A}} \textrm{ped}(P,{U}) \le {U}$ is the \textit{pedal of} ${A}$ \textit{in} ${U}$, this is a plane of the same dimension as ${A}$.
As distance $d_0({A},{U})$ we define the distance between ${A}$ and $\textrm{ped}({A},{U})$.\vspace*{0.7 mm}\\
\textit{Example}: $n=4$, $\mathfrak{A}=\textrm{diag}(1,1,1,1,0)$.\\
Put $U=E_1\sqcup E_2\sqcup E_3,  P=(1{:}0{:}0{:}1{:}0), Q=(0{:}1{:}0{:}1{:}0), A=P\sqcup Q$. Then:\\
ped$(A,U)=E_1\sqcup E_2$ and $\displaystyle\xi_0(A,U)=\xi_0(A,E_1\sqcup E_2)={2}/{3}$.\vspace*{-2.5 mm}\\

We call the $k$-plane $\text{par}(P,{U}) := ({U} \sqcap P^\mathfrak{A}) \sqcup P$ the \textit{plane} \textit{through} $P$ \textit{parallel to} ${U}$.\\ Caution: In hyperbolic geometry  two planes  are often defined to be parallel to each other if their intersection is an isotropic set. Both concepts of parallelism have to be kept entirely apart. 

\subsection{} \textbf{The angle between a line and a plane.}\;Let $H$ be an anisotropic plane and let $L$ be an anisotropic line. We assume that $L$ meets $H$ in one point $P$. The lines $L$ and $P\sqcup\textrm{ped}(L;H)$ form two angles with vertex $P$. If $Q$ is a point on $L$ and $Q\ne P$, then  \vspace*{0.5 mm}\\ 
\centerline{$\displaystyle\sinh(\,\mu(\angle_1))= \sinh(\,\mu(\angle_2))= \frac{\sinh(d_0(Q;\textrm{ped}(Q;H))}{\sinh(d_0(Q,P))}=\frac{ \cosh(d_0(Q,H^\mathfrak{A})}{\sinh(d_0(Q,P))}\, \mathbf{i}$\;.}\vspace*{1 mm}\\
If $L$ is a line totally inside the hyperplane $H$, then the angle distance between $L$ and $H$ is put to $0$. 


\subsection{} \textbf{Symmetry points, midpoints, segment bisectors and angle bisectors.} \;
Let ${U}$ be an anisotropic plane and ${S}$ be a nonempty subset of $\textrm{P}\boldsymbol{V}$. If ${S}$ is invariant under the reflection in $U$, then ${U}$ and $U^\mathfrak{A}$ are called \textit{symmetry axes} of ${S}$. If a symmetry axis consists of only one point, then this point is called \textit{symmetry point} of ${S}$. \\
Let $P$ and $Q$ be two distinct anisotropic points in $\textrm{P}\boldsymbol{V}$. If $\textrm{sgn}(P)\ne\textrm{sgn}(Q)$, then $\{P,Q\}$ has no symmetry points. If both points have the same sign, then the points $P{\,+\,}Q$ and $P{\,-\,}Q$ are called the \textit{midpoints} of $\{P,Q\}$. $P{\,+\,}Q$ is called the \textit{midpoint} of $[P,Q]_+$ and $P{\,-\,}Q$ the \textit{midpoint} of $[P,Q]_-$. At least one of the two midpoints is anisotropic and a reflection in this point leaves both midpoints fixed. Thus, both midpoints are symmetry points of $\{P,Q\}$.\\   
If $P\sqcup Q$ is an anisotropic line and  $P{\,\pm\,}Q$ an anisotropic point,  $(P{\,\pm\,}Q)^\mathfrak{A}$ is a \textit{perpendicular bisector} of the segment $[P,Q]_\mp$.

Let $A_1, A_2$ be two distinct anisotropic planes of the same dimension, and let $\angle(A_1,A_2)$ be an angle with sides $A_1$, $A_2$ and vertex $U{:=} A_1\sqcap A_2$ situated in an anisotropic plane $A_1\sqcup A_2$.  If $U$ is anisotropic, the line $L:=U^\mathfrak{A}\sqcap A$ meets $A_1\cup A_2$ in two points. If the segment $\angle(A_1,A_2)\sqcap L$ has a midpoint $M$, $M\sqcup\, U$ is called \textit{angle bisector} of  $\angle(H_1,H_2)$.

\subsection{} \textbf{Distance of two equidimensional anisotropic planes.}\hspace*{\fill}\\
Let $A_1, A_2$ $\le\mathrm{P}\boldsymbol{V}$ be two distinct anisotropic planes of the same dimension $m>0$. 
We put $E{\,:=\,}A_1\sqcup A_2$, $A_1':=E\sqcap A_1^\mathfrak{A}$, $A_2':=E\sqcap A_2^\mathfrak{A}$ and assume that $E$ is anisotropic and that $A_1\sqcap A_2' = A_2\sqcap A_1' = 0$.\vspace*{0,8 mm}\\  Then there are points $Q_1,\dots,Q_{m+1} \in A_1$ and $R_1,\dots,R_{m+1} \in A_2$ such that\vspace*{0.5 mm}\\
\hspace*{2 mm}(1)\;\;${\textrm{cosh}(d_0(Q_i,Q_j))\; = \textrm{cosh}(d_0(R_i,R_j)) = 0, \, 1\le i<j\le m+1 }$.\\
\hspace*{2 mm}(2)\;\;${\textrm{cosh}(d_0(Q_i,R_j)) = 0\, ,\, i\ne j}$.\\
\hspace*{2 mm}(3)\;\;If $Q_i\ne R_i$, then $Q_i\sqcup R_i$ is a line that meets both planes, $A_1$ and $A_2$, perpendicularly.\\
\hspace*{2 mm}(4)\;\;$\zeta_0(A_1,A_2) = \prod_{i = 1}^{m+1} (\textrm{cosh}(d_0(Q_i,R_i)))^2$.\vspace*{1 mm}\\
\textit{Proof.}\;The proof will be given in several steps.\vspace*{0.5 mm}\\
We first show that the theorem is true if $A_1$ and $A_2$ are disjoint planes of the same dimension $m>0$. The proof is by induction on $m$.\\
Suppose $A_1$ and $A_2$ are lines that do not intersect. The proof for this initial case is rather technical and we use a CAS-system for calculations.\\ Let $P_1,P_2, P_3, P_4$ be four anisotropic points and $A_1 = P_1\sqcup P_2, A_2 = P_3\sqcup P_4$. Put $c_{ij} := P_i^\circ {\scriptstyle{[\mathfrak{A}]}} P_j^\circ, 1\le i< j\le 4$. For our purposes we may assume that $c_{12} = c_{34} = 0$. We introduce the number\hspace*{1 mm}\\
$w:= \sum\limits_{i=1}^{2}\sum\limits_{k=3}^{4}\big(c_{{ii}}^{\;2}c_{kk}^{\;2}c_{jl}^4+ 2c_{ii}^{\;2}c_{kk}c_{ll}c_{jk}^ {\;2}c_{jl}^{\;2}+2c_{kk}^{\;2}c_{ii}c_{jj}c_{il}^{\;2}c_{jl}^{\;2}\big)\\ \hspace*{55 mm}+\; 2c_{11}c_{22}c_{33}c_{44}\big( c_{13}^{\;2}c_{24}^{\;2}-4c_{13}c_{14}c_{23}c_{24}+c_{14}^{\;2}c_{23}^{\;2}\big)$, \vspace*{1 mm}\\
where the indices are chosen such that $\{i,j\}=\{1,2\}$ and $\{k,l\}=\{3,4\}$.\vspace*{1 mm} \\              
We show that this number $w$ is a non-negative real number by examining four different cases.\\ 
(1) Both lines $A_1 = P_1\sqcup P_2$ and $A_2=P_3\sqcup P_4$ are elliptic lines, and both lie in $\mathcal{Q}^+$ or both in $\mathcal{Q}^-$.  Thus, $c_{11}=c_{22}=c_{33}=c_{44}\in \{-1,1\}$ and $c_{13},c_{14},c_{23},c_{24}\in \mathbb{R}$. In this case,\vspace*{0.5 mm}\\
\centerline{ $w=\big((c_{13}+c_{23})^2+(c_{14}+c_{24})^2\big)\big((c_{13}-c_{23})^2+(c_{14}-c_{24})^2\big)\in \mathbb{R}^{\ge 0}$.}\vspace*{1 mm}\\
(2) Both lines $A_1$ and $A_2$ are elliptic lines, $P_1\sqcup P_2 \subseteq\mathcal{Q}^+$, $P_3\sqcup P_4\subseteq \mathcal{Q}^-$.\\  Here, $c_{11}=c_{22}=1, c_{33}=c_{44}=-1$ and $c_{13},c_{14},c_{23},c_{24}\in \mathbb{R}\mathbf{i}$. With $c_{13}= c_{13}^{*}\mathbf{i},c_{14}= c_{14}^{*}\mathbf{i},c_{23}= c_{23}^{*}\mathbf{i},c_{24}= c_{24}^{*}\mathbf{i}$, we get\\ 
\centerline{$w=\big((c_{13}^{*}+c_{23}^{*})^2+(c_{14}^{*}+c_{24}^{*})^2\big)\big((c_{13}^{*}-c_{23}^{*})^2+(c_{14}^{*}-c_{24}^{*})^2\big)\in \mathbb{R}^{\ge 0}$.}\vspace*{1 mm}\\
(3) $A_1$ is elliptic, $A_2$ hyperbolic; $c_{11}{\,=\,}c_{22}{\,=\,}c_{33}{\,=\,}1,c_{44}{\,=\,}-1; c_{13}, c_{23} \in \mathbb{R}$, $c_{14}, c_{24} \in \mathbb{R}\mathbf{i}$.
Then,\vspace*{-1,5 mm}\\
\centerline{$w = \big((c_{13}+c_{24}\mathbf{i})^2+(c_{23}-c_{14}\mathbf{i})^2\big)\big((c_{13}-c_{24}\mathbf{i})^2+(c_{23}+c_{14}\mathbf{i})^2\big)\in \mathbb{R}^{\ge 0}$.}\vspace*{1 mm}\\
(4) Both lines $A_1$ and $A_2$ are hyperbolic lines; $c_{11}{\,=\,}c_{33}{\,=\,}1, c_{22}{\,=\,}c_{44}{\,=\,}-1$; $c_{13},c_{24} \in \mathbb{R},$ $c_{14},c_{23}\in \mathbb{R}\mathbf{i}$. With $c_{14}= c_{14}^{*}\mathbf{i},c_{23}= c_{24}^{*}\mathbf{i}$, we get\vspace*{0,5 mm}\\
\centerline{$w=\big((c_{13}+c_{24})^2+(c_{14}^{*}+c_{23}^{*})^2\big)\big((c_{13}-c_{24})^2+(c_{14}^{*}-c_{23}^{*})^2\big)\in \mathbb{R}^{\ge 0}$.}\\

Define points $Q_1, Q_2$ on $P_1\sqcup P_2$ and $R_1, R_2$ on $P_3\sqcup P_4$  by \vspace*{1 mm}\\
$Q_1 := (-\sqrt{w}+ q + 2(c_{13}c_{23}c_{22}c_{44}+ c_{14}c_{24}c_{22}c_{33}) P_1$\\ \hspace*{55 mm} $+\;(+\sqrt{w}+ q -2(c_{13}c_{23}c_{11}c_{44}+ c_{14}c_{24}c_{11}c_{33}) P_2$,\\
$Q_2 := (+\sqrt{w}+ q + 2(c_{13}c_{23}c_{22}c_{44}+ c_{14}c_{24}c_{22}c_{33}) P_1$\\ \hspace*{55 mm} $+\;(-\sqrt{w}+ q -2(c_{13}c_{23}c_{11}c_{44}+ c_{14}c_{24}c_{11}c_{33}) P_2$,\\
 \centerline{with\;  $q:= c_{11}(c_{33}c_{24}^{\;2}+c_{44}c_{23}^{\;2}) - c_{22}(c_{33}c_{14}^{\;2}+c_{44}c_{13}^{\;2})$,}\vspace*{1 mm}\\
$R_1 := (-\sqrt{w} + r + 2(c_{13}c_{23}c_{22}c_{44}+ c_{14}c_{24}c_{22}c_{33}) P_3$\\ \hspace*{55 mm} $+\;(+\sqrt{w} + r -2(c_{13}c_{23}c_{11}c_{44}+ c_{14}c_{24}c_{11}c_{33}) P_4$,\\
$R_2 := (+\sqrt{w} + r + 2(c_{13}c_{23}c_{22}c_{44}+ c_{14}c_{24}c_{22}c_{33}) P_3$\\ \hspace*{55 mm}\,$+\; (-\sqrt{w} + r -2(c_{13}c_{23}c_{11}c_{44}+ c_{14}c_{24}c_{11}c_{33}) P_4$,\\
\centerline{with\;  $r:= c_{33}(c_{11}c_{24}^{\;2}+c_{22}c_{14}^{\;2}) - c_{44}(c_{11}c_{23}^{\;2}+c_{22}c_{13}^{\;2})$.}\vspace*{1.5 mm}\\
The mapping $(A_1{\smallsetminus}{\hspace{-4.3pt}\smallsetminus}\,\mathcal{Q})\times (A_2{\smallsetminus}{\hspace{-4.3pt}\smallsetminus}\,\mathcal{Q})\to \mathbb{C},\,(S_1,S_2)\,\mapsto\, S_1{\scriptstyle{[\mathfrak{A}]}}S_2\, $, is stationary at
$(Q_1,R_1)$ and at $(Q_2,R_2)$, \\ 
and $Q_1\sqcup R_1$ and $Q_2\sqcup R_2$ are the two lines that meet each of the sets $A_1, A_2, A_1', A_2'$. Thus,\\
\hspace*{0 mm}$Q_1{=\,}\textrm{ped}(R_1,P_1 \sqcup P_2)$, $R_1{=\,}\textrm{ped}(Q_1;P_3\sqcup P_4)$,
$ Q_2{=\,}\textrm{ped}(R_2;P_1\sqcup P_2)$,\, $R_2{=\,}\textrm{ped}(Q_2;P_3\sqcup P_4).$\vspace*{1 mm}\\ 
\hspace*{0 mm}$d_0(Q_1,Q_2) = d_0(R_1,R_2) = d_0(Q_1,R_2) = d_0(R_1,Q_2) = \frac{1}{2}\pi \mathbf{i}$. Since $A_1\cap A_2 = \emptyset$, we therefore get $Q_2 \sqcup R_2 \le (Q_1 \sqcup R_1)^\mathfrak{A}.$\vspace*{1 mm} \\
Moreover, 
\hspace*{0.5 mm} $\;\zeta_0(Q_1\sqcup Q_2,R_1\sqcup  R_2)$ \vspace*{0.5 mm}\\
\hspace*{19 mm}$\displaystyle=  \frac{\big(\det\! {\left( \begin{array}{cc}
\!\!\!\scriptstyle{Q_1^\circ} {\scriptscriptstyle{[\mathfrak{A}]}} \scriptstyle{R_1^\circ}&\!\!\!\scriptstyle{0}\!\!\!  \\ 
\!\!\!\scriptstyle{0} &\!\!\!\scriptstyle{Q_2^\circ} {\scriptscriptstyle{[\mathfrak{A}]}} \scriptstyle{R_2^\circ}\!\!\!  \\
\end{array}\right)}\big)^2 }{{\det\! {\left( \begin{array}{cc}
\!\!\!\scriptstyle{Q_1^\circ} {\scriptscriptstyle{[\mathfrak{A}]}} \scriptstyle{Q_1^\circ}&\!\!\!\scriptstyle{0}\!\!\!  \\ 
\!\!\!\scriptstyle{0} &\!\!\!\scriptstyle{Q_2^\circ} {\scriptscriptstyle{[\mathfrak{A}]}} \scriptstyle{Q_2^\circ}\!\!\!  \\
\end{array}\right)}}{\det\! {\left( \begin{array}{cc}
\!\!\!\scriptstyle{R_1^\circ} {\scriptscriptstyle{[\mathfrak{A}]}} \scriptstyle{R_1^\circ}&\!\!\!\scriptstyle{0}\!\!\!  \\ 
\!\!\!\scriptstyle{0} &\!\!\!\scriptstyle{R_2^\circ} {\scriptscriptstyle{[\mathfrak{A}]}} \scriptstyle{R_2^\circ}\!\!\!  \\
\end{array}\right)}}}
$ \vspace*{1.5 mm}\\ 
\hspace*{19 mm}= $\big(\cosh(d_0(Q_1,R_1))\big)^2\big(\cosh(d_0(Q_2,R_2))\big)^2$.\vspace*{2.5 mm}\\
Let us now assume that $m>1$ and that the statement holds for $m-1$.\\ 
First of all, there exists (at least) one line $L$ which meets the planes $A_1, A_2, A_1', A_2'$. This is a consequence of a lemma in enumerative geometry (Schubert calculus) which states that, if four $m$-planes in a projective space of dimension $2m+1$ are in a "general position", then there exist $m+1$ lines which meet all these planes, see \cite{Schu} for a proof.\\
If $L$ meets $A_1$ at $Q_1$, then there exist anisotropic points $Q_2,\dots,Q_{m+1}\in A_1$ such that $d_0(Q_i,Q_j)=\frac{1}{2}\pi\mathbf{i}, 1{\le}i{<}j{\le\,}m{+}1$, and if
$L$ meets $A_2$ at $R_1$, then there are anisotropic points $R_2,\dots,R_{m+1}\in A_2$  with $d_0(R_i,R_j)=\frac{1}{2}\pi\mathbf{i}, 1{\le}i{<}j{\le\,}m{+}1$.\\
The proof for $m=1$ shows that $d_0(Q_i,R_j)=\frac{1}{2}\pi \mathbf{i}$ for $ \,i\ne j$. $U_1{:=\,}\textrm{span}(Q_2,\dots,Q_{m+1})$ and $U_2{:=\,}\textrm{span}(R_2,\dots,R_{m+1})$ are disjoint $(m{-}1)$-planes in $(P_1 \sqcup Q_1)^\mathfrak{A}$. We apply formula $(\ddagger)$ getting
{$\zeta_0(A_1,A_2)=\big(\cosh(d_0(P_1,Q_1))\big)^2\,\zeta_0(U_1,U_2)$}. We complete the proof for disjoint planes $A_1$ and $A_2$ by including the induction assumption.\vspace*{1.5 mm}\\
\hspace*{2 mm}We now assume that $A_1$ and $A_2$ intersect in a non-empty subspace $U$ of dimension $s$. We can find points $Q_1,\dots,Q_{m+1}$ which span $A_1$ such that $\textrm{span}(Q_1,\dots,Q_{s+1}) = U$ and $\textrm{cosh}(d_0(Q_i,Q_j)) = 0$ for ${1{\,\le\,}i{\,<\,}j{\le\,}}m{+}1$. Put $R_i = Q_i,\, 1\le\! i\le\! s{+}1$. We can find points $R_{s+2},\dots, R_{m+1}$ such that $\textrm{cosh}(d_0(R_i,R_j)) = 0$ for ${1{\,\le\,}i{\,<\,}j{\le\,}}m{+}1$. The intersections $A_1 \sqcap\,U^\mathfrak{A}$ and $A_2 \sqcap\,U^\mathfrak{A}$ are two disjoint subspaces of equal dimension in $U^\mathfrak{A}$. With formula $(\ddagger)$ we  get {$\zeta_0(A_1,A_2) = \zeta_0(A_1 \sqcap U^\mathfrak{A},A_2 \sqcap U^\mathfrak{A}\,)$}.  We can apply the theorem to the subspaces $A_1 \sqcap U^\mathfrak{A}$ and $A_2 \sqcap U^\mathfrak{A}$ and finish the proof. $\scriptstyle{\Box}$ \vspace*{0.7 mm}
\subsection{} \textbf{The function $\psi$ and staudtian of a finite set of points.}\hspace*{\fill} \\ 
We introduce a function $\psi$ which is defined on the powerset of $\textrm{P}\boldsymbol{V}$ and has values in $\mathbb{R}$ by: \\
$\psi(\mathcal{S}) = 0$, unless $\mathcal{S}$ consists of finitely many anisotropic points. If $\mathcal{S}$ consists of $s\in\mathbb{N}^+$ anisotropic points  ${P_1,\dots,P_s}$, then $\displaystyle\psi(\mathcal{S}) = {\det((P_i^{\circ}{\scriptstyle{[\mathfrak{A}]\,}}P_j^{\circ})_{1\leq\,i,j\leq\,s}})\,$.\\
If $\psi(\mathcal{S})\ne 0$, we call the number $\frac{1}{(s{-}1)!}\sqrt{|\psi(\mathcal{S})|}$ the \textit{staudtian} of $\mathcal{S}$, cf. \cite{H1, Ev2} for the case $n{\,=\,}2$. \vspace*{1 mm}\\
Given two anisotropic points $P_1, P_2$, then $\psi(\{P_1,P_2\})=\xi_0(P_1,P_2)=-\sinh^2(d_0(P_1,P_2))$.
Let $P_1{,\dots,}P_s$ be independent anisotropic points and $\mathcal{S}{=}\{P_1{,\dots,}P_s\}$. {If the plane span$(\mathcal{S})$ is isotropic, then $\psi(\mathcal{S}) = 0.$\newpage 
\noindent If the planes span$(\mathcal{S})$ and $P_1{\sqcup\dots\sqcup}P_{s-1}$ are anisotropic, then\\
\centerline{$-\psi(\mathcal{S}{\smallsetminus}{\hspace{-4.3pt}\smallsetminus}\{P_s\})\sinh^2(d_0(P_s,\textrm{span}(\mathcal{S}{\smallsetminus}{\hspace{-4.3pt}\smallsetminus}\{P_s\}))){\,=\,}\,\psi(\mathcal{S})$.}\vspace*{0.5 mm}}
If ${P_1,\dots,P_s}$ are anisotropic and $Q = q_1 P_1{\,+\dots+\,}q_s P_s$ with $q_1q_2{\cdots}q_s{\,\neq\,}0$ is another anisotropic point, dependent on $P_1{\,,\dots,\,}P_s$, then for $0<i,j\le s$ \vspace*{0.0 mm}\\
\centerline{$\displaystyle{\psi(\mathcal{S}{\smallsetminus}{\hspace{-4.3pt}\smallsetminus}\{P_i\}\cup \{Q\})}{\,q_j^{\,2}}{\,=\,}{\psi(\mathcal{S}{\smallsetminus}{\hspace{-4.3pt}\smallsetminus}\{P_j\}\cup \{Q\})}{\,q_i^{\;2}}\,$.} 

\subsection{} \textbf{Trigonometry.}\hspace*{\fill} \\ 
Three independent points $P_1, P_2, P_3$ in $\textrm{P}\boldsymbol{V}$ form a triangle. We assume that the plane ${P_1\sqcup P_2\sqcup P_3}$ and the three sidelines $P_2\sqcup P_3, P_3\sqcup P_1,P_1\sqcup P_2 $ of this triangle are anisotropic.
The sidelengths are $a_1:=\mu([P_2,P_3]_+), a_2:=\mu([P_3,P_1]_+), a_3:=\mu([P_1,P_2]_+)$; the measures of the interior angles are $\alpha_1:=\mu([P_1\sqcup P_2,P_1\sqcup P_3]_+),\alpha_2:=\mu([P_2\sqcup P_3,P_2\sqcup P_1]_+), \alpha_3:=\mu([P_3\sqcup P_1,P_1\sqcup P_2]_+)$.
We mention that \vspace*{1 mm}\\
\centerline{$\hspace{3pt}\sinh(a_1) = \sinh(d_0(P_2,P_3)),\dots ,\hspace{20pt} \sinh(\alpha_1)\;=\sinh(d_0(P_1\sqcup P_2,P_1\sqcup P_3)),\dots,$\hspace{0pt}}\\
\centerline{$\;\,\cosh^2(a_1) = \cosh^2(d_0(P_2,P_3)),\dots, \hspace{10pt} \cosh^2(\alpha_1)=\cosh^2(d_0(P_1\sqcup P_2,P_1\sqcup P_3)),\dots$.}\vspace*{0.7 mm}\\
The law of sines: \vspace*{1 mm}\\
\centerline{$\displaystyle\sinh^2(a_1) : \sinh^2(\alpha_1)\; =\; \sinh^2(a_2) : \sinh^2(\alpha_2)\; =\; \sinh^2(a_3) : \sinh^2(\alpha_3)\,$.}\vspace*{1 mm} 

\noindent Two laws of cosines:\\
\centerline{$\displaystyle\cosh^2(\alpha_1) = \big(\frac{\cosh(a_1)-\cosh(a_2)\cosh(a_3)}{\sinh(a_2)\sinh(a_3)}\big)^2$}\hspace{10pt}\,\\
\centerline{$\hspace{6pt}\displaystyle\cosh^2(a_1) = \big(\frac{\cosh(\alpha_1)+\cosh(\alpha_2)\cosh(\alpha_3)}{\sinh(\alpha_2)\sinh(\alpha_3)}\big)^2\,.$}\vspace*{0.5 mm}\\
\begin{figure}[!htbp]
\includegraphics[height=6cm]{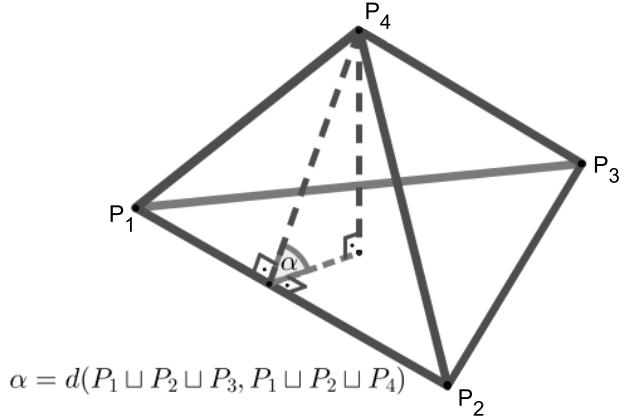}
\caption{All figures were created with the software program GeoGebra \cite{GG}.}
\vspace*{1 mm}\end{figure}
\noindent Two generalized laws of sines:\vspace*{0.5 mm}\\
Let $U\le \textrm{P}\boldsymbol{V}$ be an anisotropic plane of dimension $s\geq 2$ and let $S{\,:=\,}\{P_1,\dots,P_{s+1}\}$ be a set of independent anisotropic points in $U$. Hence $S$ is the set of vertices of an $s$-simplex in $U$ and the sets $S_i{\,:=\,}S{\smallsetminus}{\hspace{-4.3pt}\smallsetminus}\{P_i\}$, $i{\,=\,}1,\dots,s{+}1$, are the sets of vertices of its $(s{-}1)$-faces. 
Define for $i\ne j$ points $P_{ij} := (P_i\sqcup P_j)\sqcap P_i^\mathfrak{A}$ and put $\alpha_i := \big\{P_{ij}|1\le j\le s{+}1\; \textrm{and}\; j\ne i\big\}$. $\alpha_i$ can be interpreted as the \textit{interior angle} of the simplex \textit{at the vertex} $P_i$. Then, cf. \cite{Er}:\vspace*{1 mm}\\
(1)$\;\;$ If $\textrm{span}(S_1),\dots,\textrm{span}(S_{s+1})$ are anisotropic and $i \in \{1,\dots,s{+}1\}$, \vspace*{1 mm}\\
\centerline{$\displaystyle\frac{\psi(\alpha_i)}{\psi(S_i)} = \frac{\psi^{s-1}(S)}{\prod_{j= 1}^{s{+}1}\psi(S_j)}\;.$} \vspace*{0.6 mm}\\
(2)$\;\;$ If for distinct numbers $i,j \in \{1,\dots,s{+}1\}$  the three planes generated by $S_i, S_j$ and \\$\hspace*{6.5mm} S{\smallsetminus}{\hspace{-4.3pt}\smallsetminus}\{P_i,P_j\})$ are anisotropic, then \hspace{-4.3pt}  \\
\centerline{$\xi_0(\textrm{span}(S_i),\textrm{span}(S_j)\big)=-\displaystyle{\sinh^2}{\big(d\big(\textrm{span}(S_i),\textrm{span}(S_j)\big)\big)} = \frac{\psi(\mathcal{S})\,\psi(\mathcal{S}{\smallsetminus}{\hspace{-4.3pt}\smallsetminus}\{P_i,P_j\})}{\psi({S_i})\,\psi({S_j})} $\,.}\vspace*{1 mm}\\
\subsection{} \textbf{Quadrics and spheres.}\; Let $(P_1,\dots,P_{s+1})$ be a tuple of anisotropic points generating an anisotropic $s$-plane $A$.\,Let
$\mathfrak{N}{\,=\,}(\mathfrak{n}_{ij}){\,\in\,}\mathbb{R}^{(s{+}1)\times(s{+}1)}$ be an indefinite, regular symmetric matrix.\,Then $\mathcal{Q}_{\mathfrak{N}}{\,=\,}\mathcal{Q}(\mathfrak{N},(P_1{,\dots,}P_{\!{s+1}})){\,:=\,}\{r_1P_1{+\dots+\,}r_{\!{s+1}}P_{s+1}|\!\displaystyle\sum_{1\le i,j\le s{+}1}\!\!\!\mathfrak{n}_{ij}r_ir_j{=\,}0\}$ is a non-empty quadric in $A$. 
The polar of a point $R{\,=}\sum_{j}r_j P_j$ with respect to $\mathcal{Q}_{\mathfrak{N}}$  , polar$(R,\mathcal{Q}_{\mathfrak{N}})$ $:=\{\sum_{j} t_j P_j\,|\!\displaystyle\sum_{1\le i,j\le s{+}1}\!\!\!\mathfrak{n}_{ij}t_ir_j{\,=\,}0\}$, is a hyperplane of A.\vspace*{1 mm}\\
An anisotropic point $R$ is a symmetry point of $\mathcal{Q}_{\mathfrak{N}}$ iff\, polar$(R,\mathcal{Q}_{\mathfrak{N}}) = R^{\mathfrak{A}}\sqcap A$.\\
\textit{Proof}: Take any line through $R$ that meets $\mathcal{Q}_{\mathfrak{N}}$ in two points $T_1$ and $T_2$. This line meets the polar$(R,\mathcal{Q}_{\mathfrak{N}})$  at the harmonic conjugate  $R'$ of $R$ with respect to $T_1$ and $T_2$. Precisely when $R$ is an anisotropic midpoint of $\{T_1, T_2\}$, the point $R'$ is a (second) anisotropic midpoint of $\{T_1, T_1\}$ and lies on $R^\mathfrak{A}.\;\;\scriptstyle{\Box}$ \vspace*{0.5 mm}\\
\begin{figure}[!b]
\includegraphics[height=7cm]{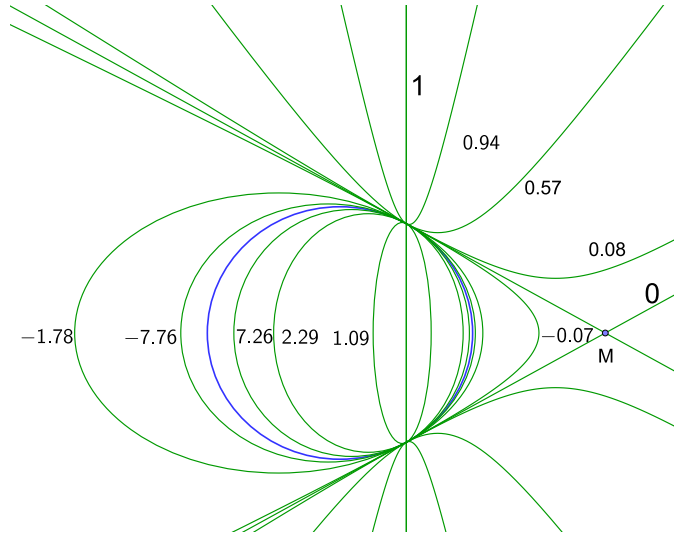}
\caption{The blue conic is the absolute in a hyperbolic plane. The green conics are 1-spheres (circles) with center $M$. The numbers indicate their quadrances.}
\end{figure}
$\mathcal{Q}_{\mathfrak{N}}$ is a \textit{hypersphere} of $A$ if there exists an anisotropic hyperplane $H$ of $A$ whose anisotropic points are symmetry points of $\mathcal{Q}_{\mathfrak{N}}$. In this case, the point $H^\mathfrak{A}\sqcap A$ is called the \textit{center} of the hypersphere.\\
Let us assume that $\mathcal{Q}_{\mathfrak{N}}$ is a \textit{hypersphere} of $A$ and $M$ is its center. Then $M^\mathfrak{A}\cap \mathcal{Q}_{\mathfrak{N}}$ is the set of isotropic points of this sphere. A sphere can have real isotropic points; but these form a thin set within this sphere.\\
Suppose the center $M$ of $\mathcal{Q}_{\mathfrak{N}}$ is anisotropic. Then $M$ is also a symmetry point of $\mathcal{Q}_{\mathfrak{N}}$. If $T_1$ and $T_2$ are two distinct anisotropic points on $\mathcal{Q}_{\mathfrak{N}}$, then $d_0(T_1,M)$ $= d_0(T_2,M)$, because: 
Let $R$ be the intersection of the line $T_1\sqcup T_2$ with the hyperplane $M^\mathfrak{A}$ of $\textrm{P}\boldsymbol{V}$, then the triple $(T_1,M,T_2)$ is mapped onto the triple $(T_2,M,T_1)$ by the (distance preserving) reflection in the point ${R}$.\\
If $M$ is anisotropic, we can define the \textit{radius} $r$ of a hypersphere as the distance between its center and any of its anisotropic points, and we call $\xi_0(\mathcal{Q}_{\mathfrak{N}}){\,:=\,}1-\cosh^2(r)$ the quadrance of $\mathcal{Q}_{\mathfrak{N}}$.\\
If the center $M$ of $\mathcal{Q}_{\mathfrak{N}}$ is isotropic, $M$ is a touchpoint of the two quadrics $\mathcal{Q}_{\mathfrak{N}}$ and $\mathcal{Q}_{\mathfrak{A}}$, and $\mathcal{Q}_{\mathfrak{N}}$ is called a \textit{horosphere}.\\
If $\mathcal{Q}_{\mathfrak{N}}$ has several (more than one) isotropic points, its center is a point outside that sphere and the lines connecting an isotropic point of $\mathcal{Q}_{\mathfrak{N}}$ with the center are \textit{tangents} of the sphere. The union of these tangents through $M$ is called the \textit{tangent}-\textit{cone} of $\mathcal{Q}_{\mathfrak{N}}$ with vertex $M$.  

Suppose $\mathcal{Q}_{\mathfrak{N}}$ has an anisotropic center $M=m_1P_1{+\dots+}m_{s+1}P_{s+1}$ and radius $r\in \mathbb{D}$, then\vspace*{0.6 mm}\\
\centerline{$\mathcal{Q}_{\mathfrak{N}}{=}\big\{r_1\!P_1{+\dots+\,}r_{s+1}\!P_{s+1}|\displaystyle (\sum_{1\le i,j\le s{+}1}\!\!\!\!\mathfrak{n}_{ij}r_im_j)^2{=} \cosh^2(r)\,(\sum_{1\le i,j\le s{+}1}\!\!\!\mathfrak{n}_{ij}m_im_j)(\sum_{1\le i,j\le s{+}1}\!\!\!\mathfrak{n}_{ij}r_ir_j)\big\}.$}\vspace*{1 mm}\\
If $\mathcal{Q}_{\mathfrak{N}}$  passes through a point $T=t_1P_1{+\dots}$ $+\,t_{s+1}P_{s+1}$,\, $\mathcal{Q}_{\mathfrak{N}}$ consists of all points ${Q}\,=r_1\!P_1{+\dots+\,}r_{s+1}\!P_{s+1}$ satisfying the equation\vspace*{1.6 mm}\\
\centerline{$\displaystyle(\sum_{1\le i,j\le s{+}1}\!\!\!\mathfrak{n}_{ij}r_im_j)^2(\sum_{1\le i,j\le s{+}1}\!\!\!\mathfrak{n}_{ij}t_it_j){\;=\;}(\sum_{1\le i,j\le s{+}1}\!\!\!\mathfrak{n}_{ij}t_im_j)^2(\sum_{1\le i,j\le s{+}1}\!\!\!\mathfrak{n}_{ij}r_ir_j).$}\vspace*{1.5 mm}\\
Let $H_1$ and $H_2$ be two hyperplanes of $A$ touching $\mathcal{Q}_{\mathfrak{N}}$ in anisotropic points $Q_1, Q_2$, respectively. If the center $M$ of $\mathcal{Q}_{\mathfrak{N}}$ is anisotropic, then $\zeta_0(H_1,M^\mathfrak{A})=\zeta_0(H_2,M^\mathfrak{A})$.\\

Here are some \textit{examples} of hyperspheres: (1) Suppose $A$ is an anisotropic line. Then two distinct anisotropic points $Q_1, Q_2\in A$ form a hypersphere iff they have a midpoint. A single anisotropic point $Q\in A$ (when counted with multiplicity $2$)  is a hypersphere in $A$ with radius $0$. (2) A (non-empty) tangent-cone of a hypersphere in an anisotropic plane $A$ with dim$(A)>1$ is a hypersphere in $A$ with radius $r=0$.  (3) If $M$ is an anisotropic point in an anisotropic plane $A$, then $A\cap M^{\mathfrak{A}}$ is (doubly counted) a hypersphere in $A$; $M$ is its center and its radius is $\frac{1}{2}\pi{\mathbf{i}}$.\vspace*{-0.5 mm}
\subsection{} \textbf{Simplices and their classical centers\footnote{$^)$ For a definition of the term \textit{triangle center} see \cite{Ki}. See also \cite{Ev1, Ev3,W2,Y} .}$^)$.}\;
Let $(P_1,\dots,P_{s+1})$ be a tuple of anisotropic points generating an anisotropic $s$-plane $A$. We denote a point $R{\,=\,}r_1P_1{\,+\dots+\,}r_{s+1}P_{s+1}\in A$ by $[r_1:\dots:r_{s+1}]$. The set $\{[r_1:\dots:r_{s+1}]\,|\, r_1r_2{\cdots}r_{s+1}$ $\neq 0\,\}$ consists of $2^s$ connected components. The closure of each component is a \textit{simplex} with vertices $P_1,\dots,P_{s+1}$. We will especially focus on the simplex $\Sigma = \Sigma_0{\,:=\,}\{[r_1{\,:\dots:\,}r_{s+1}]|\, p_1{,\dots,}\,p_{s+1} \geq 0\,\}$; it will be our reference simplex. The other simplices with these vertices will be called the \textit{mates} of $\Sigma$. \\
A plane of dimension $k\le s$ through $k{+}1$ vertices of $\Sigma$ is called a $k$-\textit{sideplane}, the intersection of $\Sigma$ with one of its $k$-sideplanes a $k$-\textit{face} of $\,\Sigma$. Instead of $1$-sideplanes we usually speak of \textit{sidelines}, and the $1$-faces are also called \textit{edges}, the $(n-1)\!$ -faces \textit{facets} of $\Sigma$. Even the vertices $P_i$ are accepted as sideplanes of $\Sigma$.\vspace*{1 mm}\\
The classical centers of an $n$-simplex are the \textit{centroid} $G$, the \textit{circumcenter} $O$ and the \textit{incenter} $I$.\vspace*{0.2 mm} \\  
We want all the edges of $\Sigma$ to have midpoints and demand that all sidelines of $\Sigma$ are anisotropic and all vertices $P_1,\dots,P_{s+1}$ have the same sign. In this case, the \textit{centroid} of $\Sigma$ is the point $ G = [1:\dots:1]$ and, more generally, the centroid of a $k$-face $\Sigma\cap (P_{i_1}\sqcup\cdots\sqcup P_{i_{k+1}})$ is the point $P_{i_1}+\dots+P_{i_{k+1}}$. The line joining the centroid of a $k$-face of $\Sigma$ with the centroid of the opposite $(s{-}k{-}1)$-face passes through $G$. 
The centroids of the mates of $\Sigma$ are the points $R{\,=\,}[r_1:\dots:r_{s+1}]$ with $P\ne G$ and  $r_j \in \{-1,1\}$ for $1\le j \le s{ +}1$. Without loss of generality we can put $r_1$ to $1$. We label the centroids of $\Sigma$ and its mates by integers from $0$ to $2^{s}-1$ as follows.  Given any centroid $[1:r_2:\dots:r_{s+1}]$, define numbers $r_j^{\star} \in \{0,1\}$ by $r_j^{\star} = 1$ iff $r_{s+1-j} = -1$. Then $G_k := [1:r_2:\cdots:r_{s+1}]$ iff $k=\sum_{\substack{0\,\leq j \leq {s}}} r_j^{\star} 2^j$. Here are two examples: $G_0 = G$ and $G_5 = [1:\cdots:1:-1:1:-1]$.\vspace*{1 mm}\\  The simplices with vertices $P_1,\dots,P_{s+1}$ will now be labeled such that $\Sigma_k$ is the simplex with $G_k\in \Sigma_k$ for $k = 0,\dots,2^{s}{-}1$. We accept each point $G_k, 0 \le k \le 2^{s}{-}1$, as a centroid of each of the simplices $\Sigma_j, 0 \le j \le 2^{s}{-}1$, but  define $G_k$ to be the \textit{proper} centroid of $\Sigma_j$ iff $k=j$.\vspace*{1 mm}\\ 
We specify the edge lengths of simplex $\Sigma_k$. The edge of $\Sigma_k$, $k=\sum_{\substack{0\,\leq j \leq {s}}} r_j^{\star} 2^j$, with vertices $P_j$ and $P_l$ has the length\vspace*{0 mm}\\
\centerline{\vspace*{0 mm}$\mu_{k;jl}\,=\begin{cases}
 \mu_{0;jl}= \mu([P_j,P_l]_+)\,, & \textrm{if}\;r^\star_{j}r^\star_{l}=1,   \\  
\mu([P_j,P_l]_-)\,, & \textrm{if}\;r^\star_{j}r^\star_{l}=-1.\\
\end{cases}$}\vspace*{1 mm}\\

A \textit{circumsphere} of $\Sigma$ is a hypersphere in $A$ passing through all vertices of $\Sigma$. There are $2^s$ circumspheres $\mathcal{S}_{0^{\,}},\dots,\mathcal{S}_{2^s-1}$ of $\Sigma$, and all of them are also hyperspheres of the mates of $\Sigma$. But as for the centroids, there is a bijective correspondence between these hyperspheres and simplices in such a way that we say that $\mathcal{S}_k$ is the \textit{{proper} circumsphere} of $\Sigma_k$. All inner points of $\Sigma_k$ are points inside the proper circumsphere  of $\Sigma_k$, but the center of this circumsphere need not lie inside $\Sigma_k$.\vspace*{-1,5 mm}  \\
The proper circumsphere of $\Sigma$ is $\mathcal{S}_0=\big\{[r_1:\dots:r_{s+1}]\,| \!\displaystyle\sum_{1\le i,j\le s{+}1}\!\!\!\sinh^2(\frac{1}{2}d_0(P_i,P_j))\,r_ir_j{\,=\,}0\big\}$. \vspace*{1 mm} \\
\textit{Proof}: It is obvious that if we intersect a circumsphere $\mathcal{S}_{i}$ with $k$-sideplanes of $\Sigma$, we get spheres of dimension $k-1$, each passing through $k+1$ of the vertices of $\Sigma$. For $k=2$ we get circumcircles of triangles. The validity of the equation is proven for these one-dimensional subsets of $S$ (see \cite{Ev2,Ev3}). From this follows the correctness of the equation. $\scriptstyle{\Box}$\vspace*{0 mm} \\

Define $\mathfrak{C}=(\mathfrak{c}_{ij})\in \mathbb{R}^{(s{+}1)\times(s{+}1)}$ by $\mathfrak{c}_{ij}:=P_i^{\,\circ}{\scriptstyle{[\mathfrak{A}]}}P_j^{\,\circ}$. The proper circumcenter of $\Sigma_0$ is $O=O_0=[r_1:\dots:r_{s+1}]$ with $(r_1,\dots,r_{s+1})\,\mathfrak{C}=(1,1,\dots,1)$, the proper circumcenter of $\Sigma_1$ is $O_1=[r_1:\dots:r_{s+1}]$ with $(r_1,\dots,r_{s+1})\,\mathfrak{C}=(1,\dots,1,-1)$ etc., cf. \cite{H1,Ev3}.\vspace*{1 mm}  \\
\textit{Proof}: We give a proof for $O=O_0$: $O$ is a point on all bisectors of segments $[P_i,P_j]_+, 0\le i{\,<\,}j{\,\le\,}s{+}1$. Therefore, $\xi_0(O,P_i-P_j)=1\;\,(0{\le}i{<}j{\le}s{+}1)$ and $O^\mathfrak{A}\sqcap A = \{[x_1\dots x_{n+1}]\,|\,x_1+\dots+x_{n+1}=0\}.\;\Box$\vspace*{2 mm}\\
\textit{Remark}: If $Q=[q_1:\dots,q_{s+1}]$ is a point in $A$ with $q_1q_2\dots q_{s+1}\ne 0$, we call the $(s{-}1)$-plane $\displaystyle Q^\Sigma{:=\,}\big\{[x_1{\,:}\dots{:\,}x_{s+1}]\,|\,\frac{x_1}{q_1}+\dots+\frac{x_{s+1}}{q_{s+1}}=0\big\}$ the $\Sigma$-polar of $Q$. For $i=0,\dots,2^{s}{-}1$ the equation $O_i:=A\sqcap (G_i^{\,\Sigma})^\mathfrak{A}$ applies.\vspace*{2 mm}\\
An equation for the radius $r$ of $\mathcal{S}_{0^{\,}}$ is\\
\centerline{ $\tanh^2(r) = \tanh^2(d_0(O,P_1)) = 1-\textrm{sgn}(P_1)\big((1,\dots,1){\scriptstyle{[\mathfrak{C}^{-1}]}}(1,\dots,1)\big)$.}\vspace*{1.6 mm}  \\
\textit{Example}: Let $P_1, P_2, P_3$ be three distinct points in the elliptic plane $(\textrm{P}\mathbb{R}^3,\mathfrak{A}=\textrm{diag}(1,1,1))$. Put $c_{ij}:=\cosh(d_0(P_i,P_j))\le 1$. Then $\mathfrak{C} =\left(\begin{array}{rrr} 
1_{\;\,} & c_{12} & c_{13} \\ 
c_{12} & 1_{\;\,} & c_{23} \\ 
c_{13} & c_{23} & 1_{\;\,}
\end{array}\right)$.\vspace*{1 mm}\\ We calculate $O$ and $\tanh^2(r)$:\vspace*{1 mm}\\ \hspace*{14.7 mm}$O=(1{\,-\,}c_{23})(1{+}c_{23}{-}c_{12}{-}c_{13})\,P_1$${\,+\,}(1{\,-\,}c_{13})(1{+}c_{13}{-}c_{12}{-}c_{23})\,P_2$ \\ \vspace*{0.8 mm}\hspace*{70 mm}${\,+\,}(1{\,-\,}c_{12})(1{+}c_{12}{-}c_{13}{-}c_{23})\,P_3$,\vspace*{1.1 mm}\\
 \hspace*{5 mm} $\displaystyle\tanh^2(r) = \frac{2 (c_{12}-1) (c_{13}-1) (c_{23}-1)}{\det(\mathfrak{C})}= \frac{2 (c_{12}-1) (c_{13}-1) (c_{23}-1)}{2 c_{12} c_{13} c_{23} - c_{12}^2  - c_{13}^2 - c_{23}^2 + 1}$.\vspace*{1.5 mm}\\
\hspace*{2 mm}
We denote the $(s{-}1)$-sideplane of $\Sigma$ opposite the vertex $P_k$ by  $H_k$ and put ${Q_k{:=\,}H_k^{\,\mathfrak{A}}{\,\sqcap\,}\,A}$. There are $2^s$ simplices $\Sigma^0,\dots,\Sigma^{2^s-1}$ with vertices $Q_1,\dots,Q_{s+1}$ and there is a $1{:}1$-corres-pondence $\Sigma_k\mapsto\Sigma^k$ given by $\Sigma^k :=  A\, {\smallsetminus}{\hspace{-4.3pt}\smallsetminus} \bigcup\limits_{R \,\textrm{inner point of}\, \Sigma_k} R^{\mathfrak{A}}$.
We call $\Sigma^k$ the \textit{dual} of $\Sigma_k$.\vspace*{1 mm}

An \textit{insphere} of $\Sigma$ is a hypersphere in $A$ which touches each $(s{-}1)$-dimensional sideplane $H_k$ of $\Sigma$. Thus, its center has the same distance from all these planes and hence the same distance from the vertices $Q_k$ of  $\Sigma^0$. In other words, the incenters of $\Sigma$ are the circumcenters of $\Sigma^0$. Since circumcenters of a simplex exist only if all its vertices have the same sign, a necessary condition for the existence of incenters of $\Sigma$ is that all vertices of $\Sigma^0$ have the same sign. We assume now, that $\Sigma_0$ has got an insphere and that the facet of $\Sigma_0$ opposite vertex $P_i$ has the staudtian $s_i$. Then $s_i/s_j\in \mathbb{R}^+$ for $1{\,\le\,}i,j{\,\le\,}n{+}1$. The proper incenter of $\Sigma_0$ is $I_0 = [s_1:{\dots}:s_{n+1}] $
 $= [1/\sinh(d_0(P_1,H_1)):\dots:1/\sinh(d_0(P_{s+1},H_{s+1}))]$. $I_0$ can also be written as a function of the interior angles $\alpha_i$ at the vertices of $\Sigma_0$: $I_0 = [\sqrt{\psi(\alpha_1)}:\dots:\sqrt{\psi(\alpha_{s+1})}\,]$.
The proper circumcenter of $\Sigma^k$ is the proper incenter of $\Sigma_k$, and the proper incenter of $\Sigma^k$ is the proper circumcenter of $\Sigma_k$. We denote the proper incenter of $\Sigma_k$ by $I_k$. This incenter $I_k$ is always a point inside the simplex $\Sigma_k$, but all the points inside this insphere may lie outside the simplex, see \cite{Ev3}. Since $I_k\in \Sigma_k$, $I_k$ is the barycentric product of  $I_0$ with the centroid $G_k$. For example, if $s=3$ and $k=2$, then $I_k=[s_1:s_2:-s_3 :s_4]$.\vspace*{1.5 mm}\\
\hspace*{2.5 mm}Let $H_i$, as before, be the $(s-1)$-sideplane of $\Sigma$ opposite vertex $P_i$. {\!}The points $R_1,\dots,R_{s+1}$, $R_i:= (P_i\sqcup H_i^{\,\mathfrak{A}})\sqcap H_i$, are called the \textit{pedals} of $\Sigma$. The tuples $P_1,\dots,P_{s+1}$ and $R_1,\dots,R_{s+1}$ are, in general, not perspective. If they are, the perspector is called \textit{orthocenter }of $\Sigma$.\vspace*{-1 mm}
\subsection{} \textbf{Radical centers of spheres.} We adopt names and assumptions from the last subsection. Now we put around each point $P_i$ an $(s{-}1)$-sphere $\tilde{\mathcal{S}}_i$ (a hypersphere in $A$) in such a way that all anisotropic points on these spheres have the same sign. We denote the radius of $\tilde{\mathcal{S}}_i$ by  $\rho_i\, (\ne \frac{1}{2}\pi\textrm{i})$ and put $r_i{\,:=\,}\cosh(\rho_i), 1\le i\le s{+}1$.\\ 
Given any point $Q\in A$, we will call the number $\tau_i(Q):=\cosh(Q,P_i)/r_i$ the \textit{power} of $Q$ with respect to the sphere $\tilde{\mathcal{S}}_i, i=1,\dots,s{+}1$.\vspace*{1 mm}\\
\begin{figure}[!htbp]
\includegraphics[height=10cm]{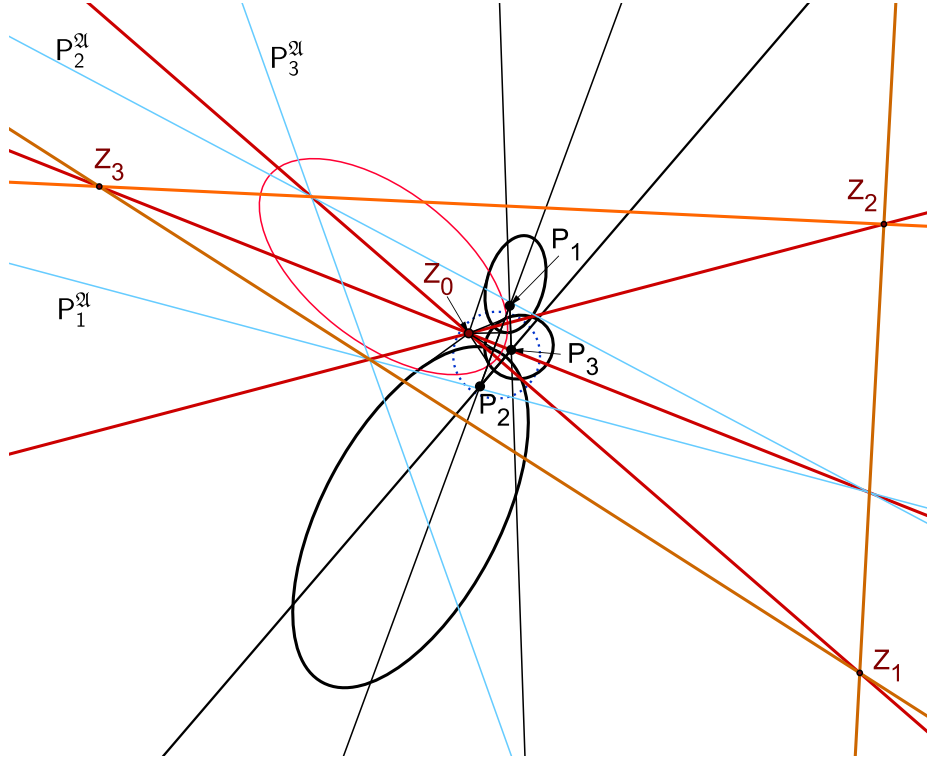}
\caption{A 3-simplex (triangle) $(P_1,P_2,P_3)$ is situated in an elliptic plane. The blue dotted circle $\{(q_1{:}q_2{:}q_3)\,|\,q_1^{\,2}+q_2^{\,2}=q_3^{\,2}\}$ is used as a construction aid and serves for orientation. The points $Z_0, Z_1,Z_2,Z_3$ are the radical centers of the black circles around the vertices of the triangle. The red and orange lines are the radical axes. The red circle is the radical circle with center $Z_0$.}\vspace*{0.8 mm}
\end{figure} 
There is a unique point $Z_0=[z_1,\dots,z_{s+1}]\in A$ with $(z_1,\dots,z_{s+1})\,\mathfrak{C}=(r_1,\dots,r_{s+1})$. $Z_0$ has the same power with respect to all spheres $\tilde{\mathcal{S}}_1,\dots,\tilde{\mathcal{S}}_{s+1}$ and is called \textit{radical center} of these spheres. $Z_0$ is a point inside the simplex $\Sigma^0$; but in each simplex $\Sigma^k,\, k\in\{1,\dots,s+1\}$, there also exists one radical center, a point $Z_k$ with $\tau_i(Z_k)=\tau_{j}(Z_k)$ for $1{\le\,}i{<}j{\le}\,s{+}1$ (cf. \cite{Ev3}). The points $Z_1,\dots,Z_{s+1}$ are the vertices of the anticevian simplex of the pair $(\Sigma^0,Z_0)$. \\
If we now assume that $\textrm{dim}(A)>1$ and $\tau_1(Z_i)>1$ for a number $i\in \{0,\dots,s+1\}$, then a hypersphere in $A$ can be drawn around $Z_i$, that meets all the spheres $\tilde{\mathcal{S}}_1,\dots,\tilde{\mathcal{S}}_{n+1}$ ortho\-gonally.\\
\begin{figure}[!t]
\includegraphics[height=11cm]{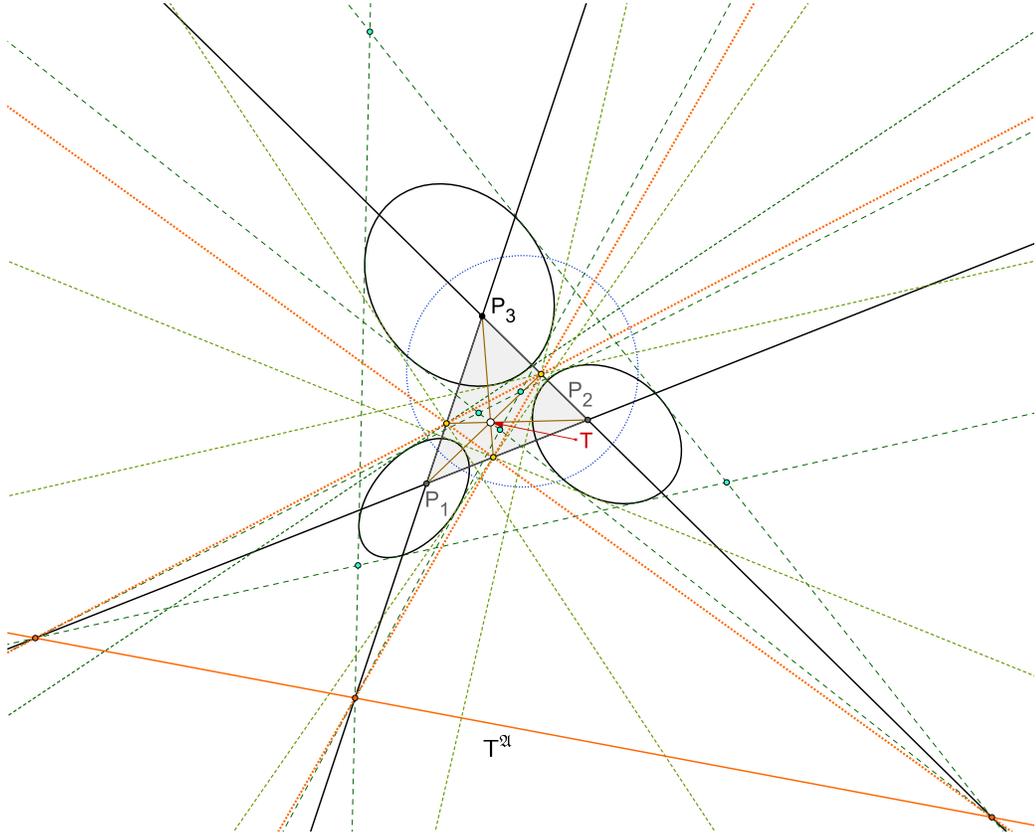}
\caption{$(P_1P_2P_3)$ is a triangle in an elliptic plane. The point $T$ is the inner center of similitude of the black circles around the vertices of this triangle.} 
\end{figure} 
A special case:  If all the radii $\rho_k$ agree, then the radical centers are the circumcenters of $\Sigma$ and $Z_0$ is the proper circumcenter $O$ of $\Sigma_0$. \vspace*{-0.2 mm} \\\vspace*{-2 mm}

Exactly two points on the line $P_1 \sqcup P_2$ have the same power with respect to the spheres $\tilde{\mathcal{S}}_1$ and $\tilde{\mathcal{S}}_2$:
{$M_1{\,=\,}(c_{12} r_2 -r_1)P_1 + (c_{12} r_1-r_2)P_2$ and 
$M_2{\,=\,}(c_{12} r_2 + r_1)P_1 - (c_{12} r_1+r_2)P_2$.}\\
If we now assume dim$(A)> 1$, then the set of all radical centers of these two spheres is the union of the two hyperplanes $H_{12} :=(M_1\sqcup (P_1\sqcup P_2)^\mathfrak{A})\sqcap A$ and $\hat{H}_{12} :=(M_2\sqcup (P_1\sqcup P_2)^\mathfrak{A})\sqcap A$ of $A$.
The four planes $P_1^{\,\mathfrak{A}}\sqcap A, {H}_{12} , P_2^{\,\mathfrak{A}}\sqcap A ,\hat{H}_{12}$ form a harmonic pencil, see Figure 3.\vspace*{0.8 mm}\\
To each sphere $\tilde{\mathcal{S}}_k$ is assigned a dual $(s-1)$-sphere $\tilde{\mathcal{S}}_k^\delta = \{ T^\mathfrak{A} \sqcap A |\; T\, \textrm{is a tangent hyperplane}$ $\textrm{of}\,\tilde{\mathcal{S}}_k \textrm{in\,} A\}$. Let us look at the two spheres $\tilde{\mathcal{S}}_1$ and $\tilde{\mathcal{S}}_2$. Their duals have together two radical hyperplanes, $H'_{12}$ and $\hat{H}'_{12}$. The duals in $A$ of these hyperplanes are the points $R_{12}^{\,\pm} = \sinh(\rho_2) P_1 \pm \sinh(\rho_1) P_2$ on the line $P_1\sqcup P_2$. These two points are called \textit{centers of similitude} of $S_1$ and $S_2$. Together with $P_1$ and $P_2$ they form a harmonic range $P_1, R_{12}^{\,+}, P_2,R_{12}^{\,-}$. If the spheres $S_1$, $S_2$ have a common tangent plane, it passes either through $R_{12}^{\,+}$ or through $R_{12}^{\,-}$. \\
We call the point $[1/\sinh(\rho_1):\cdots:1/\sinh(\rho_{n+1})]$ the \textit{inner center of similitude} of the collection $\{S_k | 1\le k\le s+1\}$, and it is obvious that subsets of this collection have their own inner center of similitude. The barycentric product of an inner center of similitude with a centroid  of $\Sigma$ gives a \textit{center of similitude}. \\
Two special cases:\, (1) If all radii $\rho_k$ are the same, then the inner center of similitude is $ G_0$.\; (2) If $\rho_k$ equals the distance between the vertex $P_k$ and its opposite $(n{-}1)$-sideplane $H_k$, $1\le k\le s{+}1$, then $I_0$ is the inner center of similitude.\vspace*{2 mm}

\section{Geometry on Cayley-Klein spaces}\vspace*{0.6 mm}
\subsection{} \textbf{Cayley-Klein spaces.}\; A \textit{Cayley-Klein space} with underlying projective space $\textrm{P}\boldsymbol{V}$ is a nested sequence of finitely many, say $\rho+1$, semi CK spaces $(A_i,{{\beta_i}})$, such that $A_0 = \textrm{P}\boldsymbol{V}$, $A_{i+1}$ is a nonempty radical rad(${\beta_i}$), $0\le i\le \rho-1$, and $A_{\rho+1}{:=\textrm{rad}(\beta_\rho)}=\emptyset.$ Since the subspaces $A_i$ can be read off from the functions $\beta_{i-1}$, we denote this CK-space by $\big(\textrm{P}\boldsymbol{V}; {\beta_0}, {\beta_1},\dots, {\beta_\rho}\big)$. \\
If $\beta_0$ is nondegenerate\,(i.e. $\rho=0$), then $(\textrm{P}\boldsymbol{V}, {\beta_0})$ is a CK-space.\vspace*{0 mm}\\
We call a point $P\in \textrm{P}\boldsymbol{V}$ \textit{anisotropic of degree} $k\in \{0,\dots,\rho\}$ if $P\in A_k $ $\textrm{\,and\,} P\not\in \mathcal{Q}_{\beta_k}$, and we call a point $P$ \textit{isotropic of degree} $k$ (or shorter $k$-\textit{isotropic}) if $P\in \mathcal{Q}_{\beta_k}\!{\smallsetminus}{\hspace{-4.3pt}\smallsetminus}\,A_{k+1}$. 
A point $P$ is called \textit{isotropic} if it is $k$-isotropic for some $k\in\{0,\dots,\rho\}$, otherwise it is called \textit{anisotropic}.\\
A mapping $\phi\in \textrm{Aut}(\textrm{P}\boldsymbol{V})$ is an \textit{automorphism} on $\big(\textrm{P}\boldsymbol{V}; {\beta_0}, {\beta_1},\dots, {\beta_\rho}\big)$ if
the restriction of $\phi$ to $A_i$ is an automorphism on $(A_i,\beta_i), 0\le i< \rho$. See \cite{Gi} for a more detailed study of these mappings.
\vspace*{-4 mm}\\
\subsection{} \textbf{The $\mathbf{\varepsilon}$-quadrance  as a distance-function on Cayley-Klein spaces.}
Let $(n^\downarrow,n^+,n^-,n^{\uparrow})$ be a quadruple of nonnegative integers with $n^\downarrow+n^++n^-+n^\uparrow=n{+}1$. A matrix $\mathfrak{A}=(\mathfrak{a})_{ij}\in \mathbb{R}^{(n+1)\times(n+1)}$ we denote by $\mathfrak{A}^{(n^\downarrow,\,n^+,\,n^-,\,n^{\uparrow})}$ if $\mathfrak{a}_{ii}{\,=\,}1$ for $n_\downarrow\le i\le n_\downarrow+n_+$, $\mathfrak{a}_{ii}{\,=\,}{-}1$ for $n^\downarrow+n^+{\,<\,}i\le n^\downarrow+n^++n^-$, and $\mathfrak{a}_{ij}{\,=\,}0$ otherwise.\vspace*{+0.5 mm}    
We now refer to the CK space introduced in the last subsection. We put $\beta_i^{\,\downarrow}{:=\,}\sum_{k<i}\,(\beta_k^{\;+}+\beta_k^{\;-}),\beta_i^{\,\uparrow}{:=\,}\sum_{k>i}\,(\beta_k^{\;+}+\beta_k^{\;-})$, $\beta_0^{\,\downarrow}{:=\,}0$,\;$\beta_\rho^{\,\uparrow}{:=\,}0$ and $\mathfrak{A}_i:=\mathfrak{A}^{(\beta_i^{\;\downarrow},\,\beta_i^{\;+},\,\beta_i^{\;-},\,\beta_i^{\;\uparrow})}$.\\
Then $\big(\textrm{P}\boldsymbol{V}; {\beta_0}, {\beta_1},\dots, {\beta_\rho}\big)$ and $\big(\textrm{P}\mathbb{R}^{n+1}; {\mathfrak{A}_0},\dots,{\mathfrak{A}_\rho}\big)$
are isomorphic CK spaces.

Let ${\mathbf{\varepsilon}}{\,\in\,}\mathbb{R}^\ast$ be an infinitesimal. We put $\displaystyle\mathfrak{A}_{\mathbf{\varepsilon}}{\,:=\,}{{\mathfrak{A}_0}}+\!\!\sum_{1\le i\le \rho}\!\!{\mathfrak{A}_i}\mathbf{\varepsilon}^{i}\,$. 
$\mathfrak{A}_{\mathbf{\varepsilon}}$ holds the complete information about the CK space $\big(\textrm{P}\mathbb{R}^{n+1}; {\mathfrak{A}_0},\dots,{\mathfrak{A}_\rho}\big)$. In the following we study the CK space $\big(\textrm{P}\mathbb{R}^{n+1}; \mathfrak{A}_\mathbf{\varepsilon}\big)$ more closely.\vspace*{0.8 mm}\\
It can be easily checked that a point $P=\mathbb{R}\boldsymbol{p}\in \textrm{P}\boldsymbol{V}$ is isotropic precisely when $\boldsymbol{p}{\;\scriptstyle{[\mathfrak{A}_{\mathbf{\varepsilon}}]}\;}\boldsymbol{p}=0.$\vspace*{0.5 mm} \\
Let $U$ be a plane of dimension $s$ generated by points $P_1=\mathbb{R}\boldsymbol{p}_1,\dots,P_{s+1}=\mathbb{R}\boldsymbol{p}_{s+1}$.
Then $U$ is called \textit{isotropic} if the matrix $(\boldsymbol{p}_i{\scriptstyle{\,[\mathfrak{A}_{\mathbf{\varepsilon}}]\,}}\boldsymbol{p}_j )_{_{ij}}$ is degenerate.
We now assume that $U$ is anisotropic. If $\tilde{U}$ is a second anisotropic plane of the same dimension $s$, generated by points $\tilde{P}_1=\mathbb{R}\tilde{\boldsymbol{p}}_1,\dots,\tilde{P}_{s+1}=\mathbb{R}\tilde{\boldsymbol{p}}_{s+1}$, we define the $\mathbf{\varepsilon}$-quadrance of $U$ and $\tilde{U}$ by \\
\centerline{$\,\displaystyle \xi_{\mathbf{\varepsilon}}(U,\tilde{U}):= \star_{\mathbf{\varepsilon}}\,\big(1-\frac{\big(\det\,({\boldsymbol{p}}_i\,{\scriptstyle{[\mathfrak{A}_{\mathbf{\varepsilon}}]}}\,\tilde{\boldsymbol{p}}_j)_{_{ij}}\big)^2}{\big(\det\,({\boldsymbol{p}}_i\,{\scriptstyle{[\mathfrak{A}_{\mathbf{\varepsilon}}]}}\,{\boldsymbol{p}}_j)_{_{ij}}\big)\big(\det\,({\tilde{\boldsymbol{p}}}_i\,{\scriptstyle{[\mathfrak{A}_{\mathbf{\varepsilon}}]}}\,{\tilde{\boldsymbol{p}}}_j)_{_{ij}}\big)}\big)$\,.}\vspace*{1.2 mm}

We call a plane $U{\,\le\,}\textrm{P}\boldsymbol{V}$ regular if $U{\,\smallsetminus}{\hspace{-4.3pt}\smallsetminus\,}A_0\ne \emptyset$, and we call a regular plane $U\le \textrm{P}\boldsymbol{V}$ a \textit{light-plane} if $\xi_{\mathbf{\varepsilon}}(Q,R)\in \mathbb{R}\,\mathbf{\varepsilon}^2$ for all regular points $P,Q\in U$. (According to this definition regular points are light-planes of dimension 0.)\vspace*{0 mm}
\subsection{} \textbf{Total polars and polar varieties.}\hspace*{\fill}\\
There are various definitions for a \textit{(total) polar of a plane in a CK-space} in the literature.
We adopt the definition given in \cite{HH}. If $P$ is an anisotropic point of degree $s\le\rho$ in $(\textrm{P}V,\mathfrak{A}_0,\dots,\mathfrak{A}_\rho)$, then any hyperplane in $\textrm{P}V$ that contains $P^{\mathfrak{A}_{s}}$ is a \textit{total polar} of $P$. Let $U\le \textrm{P}V$ be a $k$-plane. We can find an independent set $B=\{P_0,\dots,P_{k}\}$ of points that generate $U$ such that $\#(B\cap A_j)=\dim(U\cap A_j){\,+\,}1$, $0\le j\le \rho$. If $H_i$ is a total polar of $P_i$ for $0\le i\le k$ and $H_0,\dots,H_{k}$ are independent, the $(n{-}k{-}1)$-plane $\bigcap_{\,0\le i\le k}H_i$ is a
\textit{total polar} of U.\\ If $\tilde{U}$ is a total polar of $U$, then $U$ is a total polar of $\tilde{U}$. If $P$ is an anisotropic point in $U$ and $Q$ an anisotropic point in a total polar of $U$, then $\xi_{\mathbf{\varepsilon}}(P,Q)=1$.\\
The set of all total polars of a $k$-plane $U$ is called \textit{polar variety} of $U$, denoted by $U^\perp$.\vspace*{0 mm}
\subsection{} \textbf{Examples}.
(1) $\big(\textrm{P}\mathbb{R}^1;{\mathfrak{A}_\mathbf{\varepsilon}= {(1)}}\big)$ is a CK space of dimension $0$; it consists of one point, $E_1=\mathbb{R}\boldsymbol{e}_1$. $\xi_{\mathbf{\varepsilon}}(E_1,E_1)=0$. The only total polar of $E_1$ is the empty set.\\
(In all CK spaces with underlying projective space $\textrm{P}\boldsymbol{V}$ the total polar of $0$ is $\textrm{P}\boldsymbol{V}$ and vice versa.)\\
(2) The isotropic points on the hyperbolic line $\big(\textrm{P}\mathbb{R}^2;{\mathfrak{A}_\mathbf{\varepsilon}=\textrm{diag}{(1,-1)}}\big)$ are the 
points $P_1=(1:1)$ and $P_2=(1:-1)$. (They are isotropic of degree 0.) The quadrance of two distinct points on this line is not defined if one of these is $P_1=(1:1)$ or $P_2=(1:-1)$. The $\mathbf{\varepsilon}$-quadrance of two anisotropic points $Q=(q_1{:}q_2)$ and $R=(r_1{:}r_2)$ is $\displaystyle \xi_{\mathbf{\varepsilon}}(Q;R)=\frac{(q_1r_2-r_1q_2)^2}{(q_1^{\,2}-q_2^{\,2})(r_1^{\,2}-r_2^{\,2})}$. As ${\mathfrak{A}_\mathbf{\varepsilon}}$, $\xi_{\mathbf{\varepsilon}}(Q;R)$ is independent of $\mathbf{\varepsilon}$.\vspace*{+0.5 mm}\\
$Q^\perp = \{(q_2{:}q_1)\}$, $P_1^\perp = \{P_1\}$ and $P_2^\perp = \{P_2\}$.\\
(3) On the parabolic line $\big(\textrm{P}\mathbb{R}^2;{{\mathfrak{A}_\mathbf{\varepsilon}}}\big)$, $\mathfrak{A}_\mathbf{\varepsilon}=\textrm{diag}(1,\mathbf{\varepsilon})$, the point $P=(0{:}1)$ is anisotropic of degree 1, all other points are anisotropic of degree $0$.
The $\mathbf{\varepsilon}$-quadrance between two points $Q=(1{:}q_2)$ and $R=(1{:}r_2)$ is $\mathbf{\varepsilon} (q_2-r_2)^2$, the quadrance between $P$ and $Q$ is 1. The only total polar of $Q$ is the point $P$, while $P^\perp = \textrm{P}\mathbb{R}^2$.\vspace*{0.7 mm}\\ 
(4) Consider the CK space $(\textrm{P}\mathbb{R}^3; \mathfrak{A}_\mathbf{\varepsilon})$, $\mathfrak{A}_\mathbf{\varepsilon}=\textrm{diag}(1,1,\mathbf{\varepsilon})$ . It is called the \textit{co-euclidean plane} or \textit{polar-euclidean plane}. The point $E_3$ is anisotropic of degree 1; it is often called the \textit{absolute midpoint}. All other points are anisotropic of degree 0. $E_3^{\;\perp}$ is the set of all lines in $\textrm{P}\mathbb{R}^3$. $E_1^{\;\perp}$ consists of only one element, the line $E_2\sqcup E_3$. The polar variety of the line $E_1\sqcup E_2$ consists of one point, the point $E_3$, while the polar variety of the line $E_1\sqcup E_3$ is the set (the range) of points on the line $E_1\sqcup E_2$, thus $(E_1\sqcup E_3)^\perp=E_1\sqcup E_2$.\\
Let $Q=(q_1{:}q_2{:}q_3)$ and $R=(r_1{:}r_2{:}r_3)$ be anisotropic points of degree 0. If $Q$ and $R$ are not on one line with $E_3$, then $\displaystyle \xi_{\mathbf{\varepsilon}}(Q,R)= \frac{(q_1r_2-r_1q_2)^2}{(q_1^{\,2}-q_2^{\,2})(r_1^{\,2}-r_2^{\,2})}$ . If $Q,R,E_3$ are collinear, then $\displaystyle \xi_{\mathbf{\varepsilon}}(Q,R)= \mathbf{\varepsilon}\,\frac{(q_1r_3 -r_1q_3)^2 + (q_2r_3-r_2q_3)^2}{(q_1^{\,2}-q_2^{\,2})(r_1^{\,2}-r_2^{\,2})}$. The quadrance between $E_3$ and an anisotropic point $P^{{\,}^{\,}}\ne E_3$ is $\,\xi_0(E_3,P)=\xi_{\mathbf{\varepsilon}}(E_3,P) = 1$. (The quadrance between two points which are anisotropic of different degree takes always the value 1.)\vspace*{+1 mm}\\
Let us calculate the quadrance between two lines in this plane. All lines are anisotropic, and for all lines $L$ we find $\xi_{\mathbf{\varepsilon}}(L,L)=0$. Let $L_1, L_2$ be two distinct lines. If $L_1, L_2$ both pass through $E_3$, then $\xi_{\mathbf{\varepsilon}}(L_1,L_2)=\xi_{\mathbf{\varepsilon}}(S,T)$, where $S \in L_1$ and $T\in L_2$ are any two points different from $E_3$.  If just one of these lines passes through $E_3$, then  $\xi_{\mathbf{\varepsilon}}(L_1,L_2)=1$. We now assume that neither $L_1$ nor $L_2$ passes through $E_3$
and that they meet at a point $P=(p_1{:}p_2{:}p_3)$. If $Q=(q_1{:}q_2{:}q_3)\ne P$ is another point on $L_1$ and $R=(r_1{:}r_2{:}r_3)\ne P$
a point on $L_2$ we calculate $\xi_{\mathbf{\varepsilon}}(L_1,L_2)$, getting\vspace*{-0 mm}\\ 
\centerline{$\displaystyle \xi_{\mathbf{\varepsilon}}(L_1,L_2)=\mathbf{\varepsilon}\, \frac{(p_1^2+p_2^2)(\det\left(\begin{array}{rrr} 
p_1 & p_2 & p_3 \\ 
q_1 & q_2 & q_3 \\ 
r_1 & r_2 & r_3
\end{array}\right))^2}{(p_1q_2-p_2q_1)^2 (p_1r_2-p_2r_1)^2} \,$.} \vspace*{1.7 mm}\\
When we choose $Q=L_1\sqcap P^{\mathfrak{A}_{0}}$ and $R=L_2\sqcap P^{\mathfrak{A}_{0}}$, then $\xi_{\mathbf{\varepsilon}}(L_1,L_2)=\xi_{\mathbf{\varepsilon}}(Q,R)$.\\
(5) We calculate the $\boldsymbol{\varepsilon}$-quadrance of two lines in the euclidean space  $(\textrm{P}\mathbb{R}^{4}$, $\mathfrak{A}_{\boldsymbol{\varepsilon}}), \mathfrak{A}_{\boldsymbol{\varepsilon}}=\textrm{diag}(1,\boldsymbol{\varepsilon},\boldsymbol{\varepsilon},\boldsymbol{\varepsilon})$. Let $P:=(1{:}p_2{:}p_3{:}p_4)$ and $Q:=(1{:}q_2{:}q_3{:}q_4)$ be two regular points in $\textrm{P}\mathbb{R}^{4}$, and let $R:=(0{:}r_2{:}r_3{:}r_4)$, $S:=(0{:}s_2{:}s_3{:}s_4)\in A_1$ be points "at infinity". Put $\boldsymbol{p}:=(p_2,p_3,p_4)$, $\boldsymbol{q}:=(q_2,q_3,q_4)$, $\boldsymbol{r}:=(r_2,r_3,r_4)$  and $\boldsymbol{s}:=(s_2,s_3,s_4)$. Without loss of generality we may assume that $|\boldsymbol{r}|^2=r_2^{\,2}{\,+\,}r_3^{\,2}{\,+\,}r_4^{\,2}=1=|\boldsymbol{s}|^2$.\\
{If $R\ne S$, then $\xi_{\boldsymbol{\varepsilon}}(P\sqcup R,Q\sqcup S)=|\boldsymbol{r}\times\boldsymbol{s}|^2=\xi_{\boldsymbol{\varepsilon}}(R,S)$. }
If $R=S$, then $\xi_{\boldsymbol{\varepsilon}}(P\sqcup R,Q\sqcup S)=\boldsymbol{\varepsilon}|\boldsymbol{r}\times(\boldsymbol{p}-\boldsymbol{q})|^2$. This number can be interpreted as the euclidean squared distance of the parallel lines $P\sqcup R$, $Q\sqcup S$. There is one total polar of  $ R$ passing through $P$. This plane meets the line $Q\sqcup R$ at a point $Q'$ with $\xi_{\boldsymbol{\varepsilon}}(P,Q')=\boldsymbol{\varepsilon}|\boldsymbol{r}\times(\boldsymbol{p}-\boldsymbol{q})|^2.$\\
\textit{Remark}: E. Study \cite[p. 205]{Stu} introduced an angular distance of two skew lines which is a dual number made up of the angle and the euclidean distance of these lines. \vspace*{-2.0 mm}
\subsection{} \textbf{Reflections on $\big(\textrm{P}\mathbb{R}^{n+1}; \mathfrak{A}_\mathbf{\varepsilon}\big)$.}\hspace*{\fill} \\
Let $Q,R,S\in \textrm{P}\mathbb{R}^{n+1}$ be anisotropic of degree 0.\vspace*{0.4 mm}\\
(1) Assume that $U$ is a plane in $\textrm{P}\mathbb{R}^{n+1}$ and $R\in U$. The image of $Q$ under the reflection in $R$ has the same distance from $U$ as $Q$.\vspace*{0.4 mm}\\
(2) Let $Q'$ be the image of $Q$ under a reflection in $R$ and $Q''$ the image of $Q'$ under a reflection in $S$. Put $E:=Q\sqcup R\sqcup S$. If $\textrm{dim}(E){\,<\,}2$, then the distance of $Q$ and $Q''$ is twice distance of $R$ and $S$. Otherwise there is uniquely determined point $P \in E$, the pole of $R\sqcup S$ in $E$, and angle between $P\sqcup Q$ and $P\sqcup Q''$ is twice the angle between $P\sqcup R$ and $P\sqcup S$.\vspace*{0.6 mm} \\
\textit{Proof} of (1). Put $m:= \textrm{dim}(U)$. We may assume that $Q\notin U$. We can find points $P_1=\mathbb{R}(\boldsymbol{p}_1),\dots,P_{m+2}=\mathbb{R}(\boldsymbol{p}_{m+2})\in \textrm{P}\mathbb{R}^{n+1}$ such that $U=P_1\sqcup\dots\sqcup P_{m+1}$,\, $Q\in U\sqcup P_{m+2}$ and
$\boldsymbol{p}_i {\scriptstyle{[\mathfrak{A}_0]\,}} \boldsymbol{p}_j = 0$ for $1\le i<j\le m{+}2$. Put 
$\tau_i{:=\,}\boldsymbol{p}_i {\,\scriptstyle{[\mathfrak{A}_0]\,}\,} \boldsymbol{p}_i,\, i=1,\dots,n{+}2$. 
We can find real numbers $q_1,\dots,q_{m+2},r_1,\dots,r_{m+1}$ with
$Q{\,=\,}\mathbb{R}(\sum_{1\le i\le m{+}2} q_i \boldsymbol{p}_i),\, R{\,=\,}\mathbb{R}(\sum_{1\le i\le m{+}1} r_i \boldsymbol{p}_i)$.\vspace*{0.6 mm} The image of $Q$ under the reflection in $R$ is the point\\
\centerline{$Q'= \mathbb{R}\big(q_{m+2}\boldsymbol{p}_{m+2} + \sum_{1\le i\le m{+}1}(q_i+\lambda r_i)\boldsymbol{p}_i\big)$ with $\displaystyle\lambda=-2\frac{\sum_{1\le i\le m{+}1} q_i r_i \tau_i}{\sum_{1\le i\le m{+}1} r_i^{\;2} \tau_i}$.}\vspace*{1.2 mm}\\The distance between $Q$ and $U$ is the same as the distance between the points $Q$ and $\tilde{Q}:= \mathbb{R}(\sum_{1\le i\le m{+}1} q_i \boldsymbol{p}_i) = (P_{m+2}\sqcup Q)\sqcap U$ and the distance between ${Q}'$ and $U$ agrees with the distance between $Q'$ and $\tilde{Q}':= \mathbb{R}\big(\sum_{1\le i\le m{+}1}(q_i+\lambda r_i)\boldsymbol{p}_i\big) = (P_{m+2}\sqcup Q')\sqcap U$.
So we calculate $\xi_\mathbf{\varepsilon}(Q,\tilde{Q})$ and $\xi_\mathbf{\varepsilon}(Q',\tilde{Q}')$ and get\vspace*{0.7 mm}\\
\centerline{$\displaystyle\xi_\mathbf{\varepsilon}(Q,\tilde{Q}) = \xi_\mathbf{\varepsilon}(Q',\tilde{Q}') = \frac{q_{m+2}(\boldsymbol{p}_{m+2} {\,\scriptstyle{[\mathfrak{A}_\mathbf{\varepsilon}]\,}\,} \boldsymbol{p}_{m+2})}{\sum_{{1\le i\le m{+}2}}q_i^{\,2}\tau_i}\; .$} \vspace*{+2 mm}\\
Proof of (2). If $Q=R=S$, the statement is obviously true. Let us assume now that $\textrm{dim}(E)\ge 1$. Let $P_1=\mathbb{R}\boldsymbol{p}_1,P_2=\mathbb{R}\boldsymbol{p}_2,P_3=\mathbb{R}\boldsymbol{p}_3 \le \textrm{P}\mathbb{R}^{n+1}$ be independent points with $R{\,=\,}P_1, R\,\sqcup S \le P_1\sqcup P_2, R\,\sqcup S\sqcup Q\le P_1\sqcup P_2\sqcup P_3$ and $\boldsymbol{p}_1{\scriptstyle{[\mathfrak{A}_0]\,}}\boldsymbol{p}_2{\,=\,}\boldsymbol{p}_1{\scriptstyle{[\mathfrak{A}_0]\,}}\boldsymbol{p}_3{\,=\,}\boldsymbol{p}_2{\scriptstyle{[\mathfrak{A}_0]\,}}\boldsymbol{p}_3{\,=\,}0$.
(If n=1, we imagine $\textrm{P}\mathbb{R}^{2}$ embedded isometrically in a plane spanned by $P_1,P_2,P_3$.)  Put 
$\tau_i{:=\,}\boldsymbol{p}_i {\,\scriptstyle{[\mathfrak{A}_0]\,}\,} \boldsymbol{p}_i,\, i=1,2,3$. Then, $Q'{\,=\,}\mathbb{R}(-q_1\boldsymbol{p}_1{\,+\,}q_2\boldsymbol{p}_2+q_3\boldsymbol{p}_3)$\vspace*{0.4 mm}\\ and $Q''{\,=\,}\mathbb{R}((-q_1+\lambda s_1)\boldsymbol{p}_1{\,+\,}(q_2+\lambda s_2)\boldsymbol{p}_2+q_3\boldsymbol{p}_3)$ with $\displaystyle\lambda=2\frac{q_1s_1\tau_1-q_2s_2\tau_2}{s_1^{\;2}\tau_1+s_2^{\;2}\tau_2}$.\\
We first consider the case $Q\in P_1\sqcup P_2$. In this case,\vspace*{0.7 mm}\\
\centerline{$\displaystyle\xi_\mathbf{\varepsilon}(Q,Q'')  = 4\,\frac{s_2^{\;2}(\boldsymbol{p}_{2} {\,\scriptstyle{[\mathfrak{A}_\mathbf{\varepsilon}]\,}\,} \boldsymbol{p}_{2})}{s_1^{\;2}\tau_1} = 4\,\xi_\mathbf{\varepsilon}(R,S)\; $ if $\tau_2=0$}\vspace*{-0.8 mm}\\ 
and\vspace*{-0.8 mm}\\
\centerline{$\displaystyle\xi_\mathbf{\varepsilon}(Q,Q'')  = 4\,\frac{s_1^{\;2}s_2^{\;2}\tau_1\tau_2}{(s_1^{\;2}\tau_1+s_2^{\;2}\tau_2)^2} = 4\,\xi_\mathbf{\varepsilon}(R,S)\big(1-\xi_\mathbf{\varepsilon}(R,S)\big)\; $ if $\tau_2\ne 0$\;.}\vspace*{0.7 mm}\\
This proves the statement if $Q\in P_1\sqcup P_2$. (The result is independent of the choice of point $P_3$.)\vspace*{0.2 mm}\\
We now omit the restriction $Q\in P_1\sqcup P_2$ and calculate $\xi_\mathbf{\varepsilon}(P_3\sqcup R,P_3\sqcup {S})$ and $\xi_\mathbf{\varepsilon}({P_3\sqcup Q,}$ ${P_3\sqcup {Q}''})$. We have $P=P_3$ and get\vspace*{0.6 mm}\\
\centerline{$\displaystyle\xi_\mathbf{\varepsilon}(P\sqcup Q,P\sqcup Q'')  = 4\,\frac{s_2^{\;2}(\boldsymbol{p}_{2} {\,\scriptstyle{[\mathfrak{A}_\mathbf{\varepsilon}]\,}\,} \boldsymbol{p}_{2})}{s_1^{\;2}\tau_1} = 4\,\xi_\mathbf{\varepsilon}(P\sqcup R,P\sqcup S)\; $ if $\tau_2=0$\hspace*{10 mm}}\vspace*{-0.8 mm}\\ 
and\vspace*{-0.8 mm}\\
{\hspace*{9 mm}$\displaystyle\xi_\mathbf{\varepsilon}(P\sqcup Q,P\sqcup Q'')  = 4\,\frac{s_1^{\;2}s_2^{\;2}\tau_1\tau_2}{(s_1^{\;2}\tau_1+s_2^{\;2}\tau_2)^2}$}\vspace*{0.8 mm}\\{\hspace*{36 mm} = $4\,\xi_\mathbf{\varepsilon}(P\sqcup R,P\sqcup S)\big(1-\xi_\mathbf{\varepsilon}(P\sqcup R,P\sqcup S)\big)\; $ if $\tau_2\ne 0\;.\;\Box$\vspace*{0.7 mm}\\
\subsection{} \textbf{Metric affine spaces.}\hspace*{\fill} \\
We next consider a CK space $(\textrm{P}\mathbb{R}^{n+1},\mathfrak{A}_0,\dots,\mathfrak{A}_\rho)$, with $\mathfrak{A}_0=\textrm{diag}(1,0,\dots,0)$ and $\rho>0$. 
The geometry on $\textrm{P}\mathbb{R}^{n+1}{\,\smallsetminus}{\hspace{-4.3pt}\smallsetminus\,} A_1$ is called metric affine. Therefore we call $(\textrm{P}\boldsymbol{\mathbb{R}^{n+1}}; \mathfrak{A}_0,\dots,\mathfrak{A}_\rho)$ a metric affine CK space.\\
If $Q=(q_1{\,:\dots:\,}q_{n+1})\in\textrm{P}\mathbb{R}^{n+1}{\,\smallsetminus}{\hspace{-4.3pt}\smallsetminus\,} A_1$, we put $Q^\circ := (1,{q_2/q_1{\,,\dots,\,}q_{n+1}/q_1})$. Given two regular points $Q$ and $R$, then
{  $\xi_{\mathbf{\varepsilon}}(Q,R){\,=\,}\star_{\mathbf{\varepsilon}}\big(({Q^\circ-R^\circ}){\scriptstyle{\,[\mathfrak{A}_{\mathbf{\varepsilon}}]\,}}({Q^\circ-R^\circ})\big)$}\vspace*{0.4 mm}\\
This number can be interpreted as the \textit{squared distance} of the two points. (This squared distance can be negative.) For anisotropic points $P$ and $Q$ of different degree we have $\xi_{\mathbf{\varepsilon}}(P,Q)=1$.\\
Special cases: (1) If $\rho=1$ and $(A_1; {\mathfrak{A}_{1}})$ is an elliptic space, the geometry on $\textrm{P}\boldsymbol{V}{\,\smallsetminus}{\hspace{-4.3pt}\smallsetminus\,} A_1$ is called \textit{euclidean}. (2) If $n=4$ and $(A_1; {\mathfrak{A}_{1}})$ is a hyperbolic space, the geometry on $\textrm{P}\boldsymbol{V}{\,\smallsetminus}{\hspace{-4.3pt}\smallsetminus\,} A_1$ is called \textit{Minkowski-geometry (of space-time)}. (3) If $\rho=2$ and $(A_1; \mathfrak{A}_{2},\mathfrak{A}_{3})$ is a polar-euclidean space, the geometry on $\textrm{P}\boldsymbol{V}{\,\smallsetminus}{\hspace{-4.3pt}\smallsetminus\,} A_1$ is called \textit{galilean}. (4) A \textit{flag-space} can be characterized by $\rho = n$.

Let $M$ be a regular point and $r \in \mathbb{R}\,\mathbf{\varepsilon}$ be a real infinitesimal. The \textit{hypersphere} \textit{with} \textit{{center}} $M$ \textit{{and}} \textit{squared radius} $r$ is a quadric $\mathcal{S}(M,r)$ consisting of all points $Q$ which fulfill the condition $\xi_{\varepsilon}(Q,M)=r$ if $Q$ is regular, and the condition ${Q{\scriptstyle{\,[\mathfrak{A}_{1}}]\,}}Q = 0$, otherwise. 
Together with $M$ all points of $A_1$ are symmetry points of $\mathcal{S}(M,r)$, but $\mathcal{S}(M,r)$ has more symmetry points if $\rho{\,>\,}1$. Let us denote the set of symmetry points of $\mathcal{S}(M,r)$ by $\mathfrak{S}(\mathcal{S}(M,r))$.\vspace*{1.8 mm} \\
\textit{Examples}: (1) We determine the sphere $\mathcal{S}(M,r)$ with regular center $M=E_1$ and squared radius $r=1\mathbf{\varepsilon}$ for the spaces listed above, together with its symmetry points.\\
(1a) In the euclidean space of dimension $n$ this sphere is the set $\{(1{:}p_2{:}\dots{:}p_{n+1})$ $|$ $\,p_2^2+\dots+p_{n+1}^2=1\}$, and
$\mathfrak{S}(\mathcal{S}(M,1))=\{\!M\!\}\cup A_1$.\\
(1b) In the $4$-dimensional Minkowski space-time we have $\mathcal{S}(M,1)=\{(p_1{:}\dots{:}p_{5})$ $|\,p_1^2=p_2^2{+}p_3^2{+}p_4^2{-}p_{5}^2\}$ and $\mathfrak{S}(\mathcal{S}(M,r))=\{M\}\cup A_1$.\\(1c) In the galilean space of dimension $n{\,>\,}1$, this sphere is $\{(p_1{:}\dots{:}p_{n+1})$ $|\,p_1^2=p_2^2+\dots+p_{n}^2\}$ and the set of its symmetry-points is $(M\sqcup E_{n+1}) \cup A_1$. (4) In the flag-space of dimension $n{\,>\,}1$  we have $\mathcal{S}(M,1)=\{(p_1{:}\dots{:}p_{n+1})$ $\,|\,p_1^2=p_2^2\}$ and $\mathfrak{S}(\mathcal{S}(M,1))=(M\sqcup E_3\sqcup\dots\sqcup E_{n+1}) \cup A_1$. \\
(2) Now we consider spheres with center $M=E_1$ and radius $r=0$.\\ (2a) In the euclidean space of dimension $n$ this sphere is the set $\{(1{:}p_2{:}\dots{:}p_{n+1})$ $|\,p_2^2+\dots+p_{n+1}^2=0\}$. So it is the point $E_1$ with multiplicity $2$ (the double point $E_1$) and $\mathfrak{S}(\mathcal{S}(M,1))=\{M\}\cup A_1$.\\ (2b) In the $4$-dimensional Minkowski space-time $\mathcal{S}(M,0  )$ is the light cone consisting of all light lines through $E_1$. All of these light lines are null lines. $\mathfrak{S}(\mathcal{S}(M,0))=\{M\}\cup A_1$.\\
(2c) In the galilean space of dimension $n{\,>\,}1$, the sphere $\mathcal{S}(M,0)$ is the point $E_1$ with multiplicity $2$. But if we accepted the line $E_1\sqcup E_{n+1}$ as a null line (for example by demanding  $\mathbf{\varepsilon}^2=0$), then this line (with multiplicity $2$) would be $\mathcal{S}(M,0)$. Each point on this line is a symmetry point.\\ (2d) In the flag-space of dimension $n{\,>\,}1$  the sphere $\mathcal{S}(M,0)$ is the point $E_1$ with multiplicity $2$ unless the hyperplane $E_1\sqcup E_3\sqcup\dots\sqcup E_{n+1}$ is regarded a null plane.  All points of this hyperplane are symmetry points.\vspace*{1.7 mm}

An $(n{-}1)$-dimensional quadric $\mathcal{Q}$ in a metric affine CK space $\textrm{P}\boldsymbol{V}$ is called a \textit{horosphere} if there exists a regular plane $K$ of positive dimension such that $K$ is an axis of symmetry of $\mathcal{Q}$ and $\mathcal{Q}$ touches the hyperplane $A_1$ in a way that the points of $K\cap A_1$ are the points of tangency. We regard $K\cap A_1$ as the center of the horosphere.\\
\textit{Examples}: In the galilean CK space of dimension $n$ the horospheres passing through the point $E_1$ are the quadrics $\{(p_1{:}\dots{:}p_{n+1})$$|\,p_1p_{n+1}=a\,(p_2^2+\dots+p_{n}^2)\}$, $a\in \mathbb{R}^\times$, while in the $n$-dimensional flag-space these are the quadrics $\{(p_1{:}\dots{:}p_{n+1})$$|\,p_1p_{n+1}=a\,p_2^2\}$, $a\in \mathbb{R}^\times$. In the first case $E_{n+1}$
is the center of the horosphere, in the second $E_3\sqcup\dots\sqcup E_{n+1}$.\vspace*{-0.5 mm}
\subsection{} \textbf{($n{+}1$)-simplices in a metric affine CK space of dimension $n$.}\;\\
Let $P_1,\dots,P_{n+1}\in \textrm{P}\boldsymbol{V}\,{\smallsetminus}{\hspace{-4.3pt}\smallsetminus}\,A_1$ be independent regular points in a  metric affine CK space $(\textrm{P}\boldsymbol{V},{\mathfrak{A}_0},\dots,\mathfrak{A}_\rho)$ of dimension $n>1$. Then $\Sigma{\,:=\,}\{[q_1{:}\dots{:}q_{n+1}]\,|\, q_1{,\dots,}q_{n+1}\in \mathbb{R}^{\geq 0}\}$ is an $(n{+}1)$-simplex with vertices $P_1{,\dots,}P_{n+1}$. The hyperplane at infinity, $A_1$, consists of all points $Q=[q_1{:\dots:}q_{n+1}]$ with $q_1{\,+\dots+\,}q_{n+1}=0$.\vspace*{0.5 mm}\\
While all points of $\Sigma$ are regular, this is not the case for the mates of $\Sigma$.\vspace*{0.5 mm}\\
We assume that for each sideplane $U$ of $\Sigma$ the intersection of $U$ with $A_1$ is anisotropic with respect to ${\scriptstyle{\,[\mathfrak{A}_{1}]\,}}$. (It is possible to find such a collection of points $P_1,\dots,P_{n+1}$.)
It follows that none of the lines $P_i\sqcup P_j$, $1\le i<j\le n{+}1$, is a light-line, thus $d_{ij}{:=\,}\xi_{\mathbf{\varepsilon}}(P_i,P_{\!j})\not\in \mathbb{R}\,\mathbf{\varepsilon}^{2}$ for $1{\,\le\,}i{\,<\,}j{\,\le\,}n{+}1$. \\
If $Q=(q_1{:}\dots{:}q_{n+1}) \in \textrm{P}\boldsymbol{V}$ is a regular point, then the vector $\displaystyle Q^\circ=(1,\frac{q_2}{q_1},\dots,\frac{q_{n+1}}{q_1})\in \boldsymbol{V}$ can be assigned to $Q$.
If $Q$ is regular and $Q^{\circ}= q_1P_1^{\;\circ}+\dots+q_{n+1}P_{\!n+1}^{\;\,\circ}$, we write $Q=[q_1,\dots,q_{n+1}]$.
\vspace*{0.0 mm}
Given two regular points $Q=[q_1,\dots\,q_{n+1}]$ and $R=[r_1,\dots\,r_{n+1}]$, their squared distance can be calculated by $\,\xi_{\varepsilon}(Q,R)=-\sum_{1\le i<j\le n{+}1}d_{ij}(q_i-r_i)(q_j-r_j)$. (For a proof see \cite{Cox}.)\vspace*{1.3 mm}\\
\underline{Altidudes} of $\Sigma$: Put ${\mathcal{S}}:= \{P_1,\dots,P_{n+1}\}$, and for $i\in\{1,\dots,n{+}1\}$ put $S_i{\,:=\,}\textrm{span}({\mathcal{S}}{\smallsetminus}{\hspace{-4.3pt}\smallsetminus}\{P_i\})$. The absolute polar of $S_i$ is an anisotropic point $R_i$ in $A_{\rho}$.
The line $L_{i}:= P_i\sqcup R_i$ is called the \textit{altitude for} $S_i$. The squared distance of two regular points on an altitude is an infinitesimal number in $\mathbb{R}\,\mathbf{\varepsilon}^{\rho}$. Thus, the altitudes are light lines iff $\rho>1$.

Let $Q$ be the intersection of the hyperplane $S_{n+1}= P_1\sqcup\dots\sqcup P_n$ with its altidude $L_{n+1}$. Then $\psi({\mathcal{S}})= \psi(\{P_1{,\dots,}P_{n}\})\,\psi(P_{n+1},Q)\in \mathbb{R}\,\mathbf{\varepsilon}^{n+\rho-1}$. \\
If $\rho{\,=\,}1$,  $\psi({\mathcal{S}}) = \star_\mathbf{\varepsilon}(\det(1-\frac{1}{2} d_{ij})).$\vspace*{0.5 mm}\\
\textit{Remark}: If $\rho{\,=\,}1$, $\displaystyle\psi({\mathcal{S}})$ can also be calculated using the \textit{Cayley-Menger determinant}:\vspace*{1 mm}\\
\centerline{$\psi({\mathcal{S}}) = {(-1)^{n+1}2^{-n}\det(\left(\begin{array}{rr} 
0 & \boldsymbol{1\;\;\;\;\;\;\;\;}  \\  
\boldsymbol{1} & (d_{ij})_{1\le i,j\le n{+}1}\end{array}\right))}$\,.} 
\vspace*{1.2 mm}\\
\underline{The measure of dihedral angles:}\vspace*{2 mm}\\ \centerline{$\xi_{\mathbf{\varepsilon}}(S_i,S_j)\,\psi({\mathcal{S}}{\smallsetminus}{\hspace{-4.3pt}\smallsetminus}\{P_i\})\,\psi({\mathcal{S}}{\smallsetminus}{\hspace{-4.3pt}\smallsetminus}\{P_j\})=\psi({\mathcal{S}}{\smallsetminus}{\hspace{-4.3pt}\smallsetminus}\{P_i,P_j\})\,\psi({\mathcal{S}}),\; 1\le i<j\le n{+}1$.}\vspace*{3.8 mm}\\
Important centers of $\Sigma$.\vspace*{1.8 mm}\\
The proper \underline{centroid} $G$ of $\Sigma$ is the regular point $[1:\dots:1]$.  The line through a proper centroid $Q$ of a $k$-face of $\Sigma$ and 
the centroid $R$ of its opposite $(n{\,-\,}k{\,-\,}1)$-face passes through $G$, as was mentioned before;  but in addition we can state that\vspace*{0.8 mm}\\
\centerline{$(n-k)^2\xi_{\varepsilon}(G,Q)=(k+1)^2\xi_{\varepsilon}(G,R)$.}\vspace*{2 mm}\\
The \underline{circumcenter}.\;A regular point $O = (1:o_2:\dots:o_{n+1})$ is called a \textit{circumcenter} of $\Sigma$ if \vspace*{0.6 mm}\\
\centerline{$(O^\circ-P_1^{\;\circ}){\scriptstyle{\,[\mathfrak{A}_0{\,+\,}\mathfrak{A}_1]\,}}(O^\circ-P_1^{\;\circ})=(O^\circ-P_i^{\;\circ}){\scriptstyle{\,[\mathfrak{A}_0{\,+\,}\mathfrak{A}_1]\,}}(O^\circ-P_i^{\;\circ}),\, 1{\,\le\,}i{\,\le\,}n{+}1\;.$}\vspace*{0.4 mm}\\
$\Sigma$ has exactly one regular circumcenter iff $\rho{\,=\,}1$.\\ 
\textit{Proof.} We rearrange the last system of equations:\vspace*{0.6 mm}\\
\centerline{$(P_i^{\;\circ}-P_{\!1}^{\;\circ}){\scriptstyle{\,[\mathfrak{A}_0{\,+\,}\mathfrak{A}_1]\,}} O^\circ=(P_i^{\;\circ}-P_{\!1}^{\;\circ}){\scriptstyle{\,[\mathfrak{A}_0{\,+\,}\mathfrak{A}_1]}}(\frac{1}{2} (P_i^{\;\circ}+P_{\!1}^{\;\circ})),\, 1{\,\le\,}i{\,\le\,}n{+}1\;.$}\vspace*{0.4 mm}\\
This system consists of $n{+}1$ equations. The first of these equations is trivially fulfilled. The remaining $n$ equations form a system of maximal rank and there is a single solution for the n-tuple $(o_2,\dots,o_{n+1})$ precisely when $\rho=1$. 
We now assume that $\rho=1$. Then the circumcenter $O$ of $\Sigma$ is a regular point. We want to determine an $(n{+}1)$-tuple $(o_1{:}\dots{:}o_{n+1})$ of real numbers with $O=[o_1{:}\dots{:}o_{n+1}]$ 
and define an $(n{+}1)\times(n{+}1)$-matrix $\tilde{\mathfrak{D}}$ by $\displaystyle\tilde{\mathfrak{d}}_{ij}=1{\,-\,}\frac{1}{2}d_{ij}$.
We can find a tuple $(o_1,\dots,o_{n+1})$ such that:\vspace*{-0.8 mm}\\
\centerline{$(o_1,\dots,o_{n+1})\tilde{\mathfrak{D}}=(1,1,\dots,1)\hspace*{2.8 mm}(\star).$}\vspace*{1.4 mm}
Using Cramer's rule, we get $O=[\det\tilde{\mathfrak{D}}^{[1]}{\,:\,}\det\tilde{\mathfrak{D}}^{[2]}{\,:\,}\dots{\,:\,}\det\tilde{\mathfrak{D}}^{[n+1]}]$, where $\tilde{\mathfrak{D}}^{[i]}$ is the matrix formed by replacing the $i$-th row of  $\tilde{\mathfrak{D}}$ by the row $(1,\dots,1)$.\vspace*{1.4 mm}\\
The circumsphere consists of all points $X=[x_1,\dots,x_{n+1}]$ satisfying the equation\vspace*{0.4 mm}\\
\centerline{ $\sum\limits_{1{\le}i{<}j{\le}{n{+}1}} d_{ij}x_ix_j=0\hspace*{2.8 mm}(\star\star)\,.$}
\vspace*{-1 mm}\\
The squared radius of the circumsphere is $\displaystyle r=-\frac{{\det({d}}_{ij})}{2\det(\left(\begin{array}{rr} 
0 & \boldsymbol{1\;\;\;\;\;\;\;\;}  \\  
\boldsymbol{1} & (d_{ij})_{1\le i,j\le n{+}1}\end{array}\right))}.$ \vspace*{2.7 mm}\\
\textit{Examples}: (1) In a 2-dimensional CK space with $\rho=1$ we get\vspace*{0.4 mm}\\
\centerline{ $O=[d_{23} (d_{12} + d_{13} - d_{23}),d_{13} (d_{12} - d_{13} + d_{23}),d_{12} (-d_{12} + d_{13} + d_{23})]$.}\vspace*{0.4 mm}\\
(2) The two above formulae $(\star),(\star\star)$ lead to the correct result also for the 2-dimensional CK space with $\rho=2$ (the galilean plane):
A point $[q_1{:}q_2{:}q_3]$ in a metric affine CK space of dimension $2$ is a point on the line $A_1$ iff $q_1{+}q_2{+}q_3$ $=0$. It follows that the circumcenter  $O=$ $[d_{23} (d_{12}{+}d_{13}{-}d_{23}):d_{13} (d_{12}{-}d_{13}{+}d_{23}):d_{12} ({-}d_{12}{+}d_{13}{+}d_{23})]$ is a point on $A_1$ iff ${-\,}d_{12}^2{-}d_{13}^2{-}d_{23}^2{+}2(d_{12}d_{13}{+}d_{12}d_{23}{+}d_{13}d_{23})=0$. But this equation is a necessary and sufficient condition for the plane to be a galilean plane, cf \cite{Ev1}.\\ 
We will except $\{ [x_1:x_2:x_{3}]\,|\,\sum\limits_{1{\le}i{<}j{\le}{3}} d_{ij}x_ix_j=0\,\}$ as a circumcircle of the triangle $(P_1,P_2,P_3)$
even though we know that it is a horocircle.\vspace*{1 mm}\\

We give up the restriction $\rho = 1$. A plane $M\le A_1$ is called the \textit{circumcenter} of $\Sigma$ if there exists a horosphere with center $M$ passing through all vertices of $\Sigma$.\vspace*{0.8 mm}\\
For $\Sigma$ there is a uniquely determined plane in $\textrm{P}\boldsymbol{V}$ which is its circumcenter. As before, we use the letter $O$ to denote this plane. We call $G\sqcup O$ the \textit{Euler-plane} of $\Sigma$. \vspace*{-0.8 mm}\\


The \underline{incenter}.\;A point of $\Sigma$ is an \textit{incenter} of $\Sigma$ if it is the center of an hypersphere (called an insphere) which touches each of $(n{\,-\,}1)$-sideplanes of $\Sigma$. The squared radius $r$ of this insphere is a nonzero real infinitesimal $\in \mathbb{R}\boldsymbol{\varepsilon}^\rho$. Put $\mathcal{S}:= \{P_1,\dots,P_{n+1}\}$ and $\mathcal{S}_i:=\{P_1,\dots,P_{n+1}\}{\smallsetminus}{\hspace{-4.3pt}\smallsetminus}\{P_i\}$.\\
Necessary and sufficient for the existence of an incenter is that  $\psi(\mathcal{S}_i)\,\psi(\mathcal{S}_j){\,\geq\,}0$ for all $1{\le}i{<}j{\le}{n{+}1}$. If all these inequalities are true, the incenter is the point\\
\centerline{ $I=[\sqrt{(|\psi(\mathcal{S}_1)|}{:}\dots$ ${:}\sqrt{|\psi(\mathcal{S}_{n+1})|}\;]$.}\vspace*{-0.3 mm}\\ The squared radius of the insphere can be calculated by $\displaystyle r=\frac{\,\psi(\mathcal{S})}{\sum_{1\le i\le {n+1}} \psi(\mathcal{S}_i)}$.\vspace*{0.6 mm}\\

The \underline{Monge point}.\;We recall that $G$ is the proper centroid of $\Sigma$. 
For $1{\,\le\,}i{\,<\,}j{\,\le\,}n{+}1$ let $\check{H}_{ij}\le A_1$ be the total polar of the edge $P_i\sqcup P_j$ of $\Sigma$ and let $\tilde{H}_{ij}$ be the hyperplane $\displaystyle\check{H}_{ij}\sqcup Z_{ij}$ where $Z_{ij}:= \mathbb{R}(((n+1)G^\circ)-P_i^{\,\circ}-P_j^{\,\circ})$ is the centroid of the $(n{-}2)$-face opposite $P_i\sqcup P_j$. If $\rho=1$, the $\frac{1}{2}n(n+1)$ hyperplanes  $\tilde{H}_{ij}$ meet at one point $H$, which is called the \textit{Monge point} of $\Sigma$. For $G, O$ and $H$  the equation $\displaystyle H^\circ = \frac{n+1}{n-1}\,G^\circ - \frac{2}{n-1}\,O^\circ$ applies.\vspace*{1 mm}\\ 
\textit{Proof}: \vspace*{0.5 mm}\\ We have to show that there is exactly one point $H$ that satisfies the following system of linear equations:
{$(P_i^{\;\circ}-P_{\!j}^{\;\circ}){\scriptstyle{\,[\mathfrak{A}_0{\,+\,}\mathfrak{A}_1]\,}} H^\circ=(P_i^{\;\circ}-P_{\!j}^{\;\circ}){\scriptstyle{\,[\mathfrak{A}_0{\,+\,}\mathfrak{A}_1]}}\,Z_{ij}^{\;\circ},\, 1{\,\le\,}i{\,<\,}j{\,\le\,}n{+}1\;,$} and also satisfies the equation $H^\circ = \frac{n+1}{n-1}\,G^\circ - \frac{2}{n-1}\,O^\circ$. First, the linear system has at most one solution for $H^\circ$.  
The proof of this statement is quite analogous to the proof for the circumcenter and is therefore omitted. 
On the other hand, $H^\circ = \frac{n+1}{n-1}\,G^\circ - \frac{2}{n-1}\,O^\circ$ is a solution of the system:\\

\noindent\hspace*{0.6 mm}$(P_i^{\;\circ}{\,-\,}P_{\!j}^{\;\circ}){\scriptstyle{\,[\mathfrak{A}_0{\,+\,}\mathfrak{A}_1]\,}}((n{+}1)G^\circ-2O^\circ) = (P_i^{\;\circ}-P_{\!j}^{\;\circ}){\scriptstyle{\,[\mathfrak{A}_0{\,+\,}\mathfrak{A}_1]\,}}((n{+}1)G^\circ-(P_i^{\circ}+P_j^{\circ}))$\\
 $\hspace*{54.2 mm}= (P_i^{\;\circ}{\,-\,}P_{\!j}^{\;\circ}){\scriptstyle{\,[\mathfrak{A}_0{\,+\,}\mathfrak{A}_1]\,}}(n-1)Z_{ij}^{\;\circ},\,1{\,\le\,}i{\,<\,}j{\,\le\,}n{+}1.\;\Box$\vspace*{0.6 mm}\\
\vspace*{0.4 mm}\\
\textit{Remarks}: (1) If $n=2$, the hyperplanes $\tilde{H}_{12},\tilde{H}_{13},\tilde{H}_{23}$ are lines, which are called \textit{altitudes} of the triangle. 
These three altitudes meet at the point $H$, which, in this case, is called the \textit{orthocenter} of $\Sigma$. (2) For a tetrahedron in a $3$-dimensional euclidean space the point $H$ was discovered by G. Monge \cite{Mo}. (3) If $Q$ is a point on the circumsphere of $\Sigma$, then the point $R$ with $\displaystyle R^\circ=\frac{n{-}1}{n}H^\circ+\frac{1}{n}\,Q^\circ$ is a point on the \textit{Feuerbach sphere} of $\Sigma$, which is the $(n{-}1)$-sphere through the $n{+}1$ centers of the facets of $\Sigma$; see \cite{Fr} for a proof.\vspace*{-2mm}\\

There is a series of centers of $\Sigma$ whose construction is similar to that of O and H. We describe four of them. Again $\rho=1$ is assumed.\\
(1) Let $Q_{ij}$ be the mirror image of $P_i+P_j$ in the point $Z_{ij}$. The hyperplanes $\check{H}_{ij}\sqcup Q_{ij}$ meet at a point.
(2) Let $T_1,\dots,T_{n+1}$ be the reflections of the vertices $P_1,\dots,P_{n+1}$ in their opposite siteplanes. Put $T_{ij} = T_i+T_j$. The hyperplanes $\check{H}_{ij}\sqcup T_{ij}, 1\le i\le j\le n{+}1,$ meet at a point.
(3) Let $T_1,\dots T_{n+1}$ be the points introduced in (2). Let $\hat{H}_{ij}$  be the total polar of the line $T_i\sqcup T_j$. The hyperplanes $\hat{H}_{ij}\sqcup (P_i+P_j), 1\le i\le j\le n{+}1,$ meet at a point. (4) Let $T_1,\dots T_{n+1}$, $Z_{ij}, 1\le i<j\le n{+}1$, be the points and $\hat{H}_{ij}$  be the plane as in (1) and (3). The hyperplanes $\hat{H}_{ij}\sqcup Z_{ij}, 1\le i\le j\le n{+}1,$ meet at a point. \vspace*{-2 mm}\\

We calculate the barycentric coordinates of the centers $O, I, H$ and the point $(S_4\sqcap A_1)^{\mathfrak{A}_1}$ for an $4$-simplex in a euclidean space of dimension $3$:\\
Let $d_{ij}\boldsymbol{\varepsilon}$ be the squared distance of the vertices $P_i$ and $P_j$, then\vspace*{1 mm}\\
$O=[d_{12}(d_{23}+d_{24}-d_{34})d_{34}+d_{13}(d_{23}-d_{24}+d_{34})d_{24}+d_{14}(-d_{23}+d_{24}+d_{34})d_{23}+2d_{23}d_{24}d_{34}:\dots \textrm{cyclic}\dots]$,\vspace*{0.6 mm}\\
$I=[\sqrt{2(d_{23}d_{24}+d_{23}d_{34}+d_{24}d_{34})-(d_{23}d_{24}d_{34})}:\dots \textrm{cyclic}\dots]$,\vspace*{0.6 mm}\\
$H=[d_{23}(d_{14}-d_{12})(d_{14}-d_{13})+d_{24}(d_{13}-d_{12})(d_{13}-d_{14})+d_{34}(d_{12}-d_{13})(d_{12}-d_{14})-d_{23}d_{24}d_{34}:\dots \textrm{cyclic}\dots]$\vspace*{0.6 mm},\vspace*{0.6 mm}\\
$(S_4\sqcap A_1)^{\mathfrak{A}_1}=[-d_{14}d_{23}{+}d_{24}(-d_{12}{+}d_{13}{+}d_{23})+d_{34}(d_{12}{-}d_{13}{+}d_{23})+d_{23}(d_{12}{+}d_{13}{-}d_{23}):\\\hspace*{2.3 mm}-\,d_{24}d_{13}+d_{14}(-d_{12}+d_{13}+d_{23})+d_{13}(d_{12}-d_{13}+d_{23})+d_{34}(d_{12}+d_{13}-d_{23}):\\\hspace*{2.3 mm}-\,d_{34}d_{12}+d_{12}(-d_{12}+d_{13}+d_{23})+d_{14}(d_{12}-d_{13}+d_{23})+d_{24}(d_{12}+d_{13}-d_{23}):\\\hspace*{3 mm}{-\,}d_{12}^2-d_{13}^2-d_{23}^2+2(d_{12}d_{13}+d_{12}d_{23}+d_{13}d_{23})]$.\vspace*{-2 mm}\\


The \underline{medial simplex} and the \underline{anticomplementary simplex} of $\Sigma$. The \textit{medial simplex} $\Sigma^{[1]}$ of $\Sigma$ is the simplex whose vertices are the centroids of the $(n{\,-\,}1)$-faces of $\Sigma$. We get $\Sigma^{[1]}$ by applying a homothety with center $G$ (the centroid of $\Sigma$) and factor $\displaystyle -\frac{1}{n}$ to the simplex $\Sigma$. We can iterate this process getting a sequence of simplices where $\Sigma^{[n+1]}$ is the medial simplex of $\Sigma^{[n]}$.
The inverse of a homothety with center $G$ and factor $\displaystyle -\frac{1}{n}$ is a homothety with center $G$ and factor $\displaystyle {-\,}{n}$. The simplex $\Sigma^{[-1]}$ is called the \textit{anticomplementary} of $\Sigma$.\vspace*{0.5 mm}\\
If $T\le A_0$ is some \textit{center} \textit{(center plane)} of $\Sigma$ such as $G$ or $O$ or $G\sqcup O$ or $G\sqcup O\sqcup I$, then we denote the correspondent centers of $\Sigma^{[n]}$ by $T^{[n]}$, $n\in \mathbb{Z}$. There are centers (center planes) $T$ of $\Sigma$ such that $T = T^{[n]}$ for all $n\in \mathbb{Z}$; these include $G$, $G\sqcup O$, $G\sqcup I$ and $G\sqcup O \sqcup I$. 
\section{Addition: Clifford algebra on a semi CK space}
The foundations of Geometric Algebra (GA) were laid by H. Grassmann and W. Clifford. A description and appreciation of their works is given by D. Hestenes in \cite{He}. Due to research by Hestenes, GA experienced a revival and  is now an important part of mathematics and a number of its applications. Short introductions to GA offer \cite{Ch,KH,LS}.\vspace*{0 mm} 
\subsection{Grassmann algebras on $\mathbb{R}^{n+1}$ and on \textrm{P}$\mathbb{R}^{n+1}$} 
We make use of the Grassmann algebra $\Lambda \boldsymbol{V}$ on $\boldsymbol{V}=\mathbb{R}^{n+1}$. For $0\leq k\leq n{+}1$, the elements of $\Lambda^k \boldsymbol{V}$ are called \textit{multivectors} of grade $k$. Conventionally the wedge product $\wedge$ is used as the exterior product on $\Lambda \boldsymbol{V}$; but here we use the operator $\vee$ which fits better with the geometric interpretation. Starting from the canonical basis $(\boldsymbol{e}_1,\dots,\boldsymbol{e}_{n+1})$ of $\boldsymbol{V}$, we get ordered bases for $\Lambda^k \boldsymbol{V}$, $1\leq k\leq n{+}1$, as follows: The set of canonical basis elements of $\Lambda^k V$ is $\{\boldsymbol{e}_{i_1}\vee \dots \vee \boldsymbol{e}_{i_k}\,|\, i_1 < \dots < i_k\}$, and the elements within this set are placed according to the lexicographic order of their multi-indices. Instead of $\boldsymbol{e}_{i_1}\!\vee \dots \vee \boldsymbol{e}_{i_k}$, we also write $\boldsymbol{e}_{i_1,\dots,\,i_k}$. Let $\mathcal{S}_k$, $k\in \mathbb{N}$, denote the symmetric group on the set $\{1,\dots,k\}$.  A multivector $\boldsymbol{v}$ of rank $k \leq n{+}1$ can always be written \vspace*{-1.2mm}
\[\vspace*{-1mm}\boldsymbol{v}\, =\sum_{\substack{\sigma \in \mathcal{S}_{n+1}\\ \sigma(1)<\dots<\sigma(k)}}v_{\sigma(1),\dots,\sigma(k)}\,\boldsymbol{e}_{\sigma(1),\dots,\sigma(k)}\;.\vspace*{0mm}\]
We call $v_{\sigma(1),\dots,\sigma(k)}$ the \textit{component} of $\boldsymbol{v}$ \textit{with multi-index} $\sigma(1),\dots,\sigma(k)$.\\
Let $\boldsymbol{u}$ be a multivector of grade $r$ and $\boldsymbol{v}$ a multivector  of grade $s$, and we assume that $r{+}s\leq n{+}1$. Then we get a multivector  $\boldsymbol{w} = \boldsymbol{u}\vee \boldsymbol{v}$ of grade $r{+}s$ with\vspace*{-1mm}
\[{w}_{\sigma(1),\dots,\sigma(r+s)}\;=\sum_{\substack{\sigma \in \mathcal{S}_{r+s}\\ \sigma(1)<\dots<\sigma(r)\\ \sigma(r+1)<\dots<\sigma(r+s)}}\textrm{sgn}(\sigma)\,u_{\sigma(1),\dots,\sigma(r)}\,v_{\sigma(r+1),\dots,\sigma(r+s)}\;. \]\vspace*{-2mm}\\
The exterior product $\vee:\Lambda \boldsymbol{V}\times \Lambda \boldsymbol{V}\to \Lambda \boldsymbol{V}$ is associative.

The vector space $\Lambda V$ the can be projectivized and the exterior product $\vee$ can be transferred to an exterior product $\vee$ on $\textrm{P}\Lambda \boldsymbol{V}$ by $\mathbb{R}\boldsymbol{v} \vee \mathbb{R}\boldsymbol{w}:=\mathbb{R}(\boldsymbol{v}\vee\boldsymbol{w})$. 
The elements of $\textrm{P}\Lambda^k V$ are called $k$\textit{-flats}. A $k$-flat $F$ is called \textsl{decomposable} if there are $k$ independent points $P_1=\mathbb{R}\boldsymbol{v}_1,\dots, P_{k}=\mathbb{R}\boldsymbol{v}_k$ such that $F = \mathbb{R}(\boldsymbol{v}_1\vee{\dots}\vee \boldsymbol{v}_k)$. In this case we write $F = P_1\vee{\dots}\vee P_{k}$. If $\sigma\in \mathcal{S}_k$ is a permutation, $P_{\sigma(1)}\vee\dots\vee P_{\sigma(k)}=P_1\vee{\dots}\vee P_{k}$.\vspace*{0.5mm}

The points $P_1{,\dots,}P_{k}$ are independent precisely when $P_1\vee\dots\vee P_k$ $\ne 0$. If $P_1\vee\dots\vee P_k$ ${\ne 0}$, $P_1\sqcup\dots\sqcup P_k$ is  a $(k{-}1)$-plane in $\textrm{P}\boldsymbol{V}$ which consists of all points $Q \in\textrm{P}\boldsymbol{V}$ satisfying $P_1\vee\dots\vee P_k\vee Q = 0$. More generally, if $P_1\vee\dots\vee P_j\ne 0$ and $Q_1\vee\dots\vee Q_k \ne 0$, then $P_1\vee\dots\vee P_j\vee Q_1\vee\dots\vee Q_k \ne 0$ precisely when the planes $P_1\sqcup\dots\sqcup\,P_j$ and $Q_1\sqcup\dots\sqcup\,Q_k$ are disjoint.
\subsection{{The Pl\"ucker embedding of a plane}}
There is an embedding $\iota{:\,}\textrm{sub}(\textrm{P}\boldsymbol{V})\to \textrm{P}\Lambda \boldsymbol{V}$, the \textit{Pl\"ucker embedding}, defined by $\iota(U){\,:=\,}$ ${P_1\vee\dots\vee P_k}$, where $P_1,\dots,P_k$ is a minimal generating set of $U$. (It can be easily verified that this mapping is well-defined.) The image $\iota(U)$ is a cut out in $\textrm{P}\Lambda^k \boldsymbol{V}$ by a set of quadratic equations, so called \textit{Pl\"ucker relations}. \\
\textit{Example}: $\sum_{1\le i<j\le 4}\,l_{ij}\,{E_i{\,\vee\,}E_j} \in \textrm{P}\Lambda^2\mathbb{R}^{3}$ is the exterior product of two elements in $\textrm{P}\Lambda\mathbb{R}^{3}$
iff $l_{12}l_{34}+l_{13}l_{24}+l_{14}l_{23}=0$.
\subsection{Clifford's geometric product}
A (generalized) inner product ${\scriptstyle{[\mathfrak{A}]}}$ on $\boldsymbol{V}=\mathbb{R}^{n+1}$ can be expanded to a (generalized)  inner product on $\Lambda \boldsymbol{V}$: Given vectors $\boldsymbol{u},\boldsymbol{u}_1,\dots,\boldsymbol{u}_r,\boldsymbol{v}_1,\dots,\boldsymbol{v}_r\in \boldsymbol{V}$ and a multivector $\boldsymbol{w}\in \Lambda^k \boldsymbol{V}$, then \vspace*{-1 mm}  
\[
\begin{split}
\boldsymbol{u}{\scriptstyle{[\mathfrak{A}]}} (\boldsymbol{v}_1\vee \dots \vee \boldsymbol{v}_r) \, &=\sum_{1\leq j\leq r} (-1)^j\,(\boldsymbol{u}{\scriptstyle{[\mathfrak{A}]}} \boldsymbol{v}_j) (\boldsymbol{v}_1\vee \dots \vee \widehat{\boldsymbol{v}}_j\vee \cdots \vee\boldsymbol{v}_r),\\ 
(\boldsymbol{v}_1\vee \dots \vee \boldsymbol{v}_r){\scriptstyle{[\mathfrak{A}]}} \boldsymbol{w} &= \boldsymbol{v}_1{\scriptstyle{[\mathfrak{A}]}}(\boldsymbol{v}_2{\scriptstyle{[\mathfrak{A}]}}(\cdots {\scriptstyle{[\mathfrak{A}]}}(\boldsymbol{v}_r{\scriptstyle{[\mathfrak{A}]}} \boldsymbol{w})\cdots)),\; \textrm{if}\,r\leq k,\\ 
\textrm{and}\;\;\;(\boldsymbol{u}_1\vee \dots \vee \boldsymbol{u}_r){\scriptstyle{[\mathfrak{A}]}} (\boldsymbol{v}_r\vee \dots \vee \boldsymbol{v}_1)\; &= \textrm{det}(\boldsymbol{u}_i{\scriptstyle{[\mathfrak{A}]}} \boldsymbol{v}_j)_{1\leq i,j \leq r}\;.
\end{split}
\]
$\noindent\hspace*{3 mm}$ Given $\boldsymbol{v}, \boldsymbol{w} \in \boldsymbol{V}$, the \textit{geometric product} of these two vectors is defined by\\
\centerline{$\boldsymbol{v}\,\boldsymbol{w}:= \boldsymbol{v}{\scriptstyle{[\mathfrak{A}]}}\boldsymbol{w}+ \boldsymbol{v}\vee\boldsymbol{w}$.}\vspace*{0.8 mm}\\
As a direct consequence of this equation we get:\\
$\noindent\hspace*{15 mm}\displaystyle (1)\;\;\; \boldsymbol{v}{\scriptstyle{[\mathfrak{A}]}}\boldsymbol{w}= \frac{1}{2}(\boldsymbol{w}\,\boldsymbol{v}+\boldsymbol{v}\,\boldsymbol{w}),\hspace*{8 mm}(2)\;\;\;\boldsymbol{v}{\,\vee\,}\boldsymbol{w}=\frac{1}{2}(\boldsymbol{w}\,\boldsymbol{v}-\boldsymbol{v}\,\boldsymbol{w})\;$.\vspace*{-1.0 mm}\\
(3) If $\boldsymbol{v}$ is a vector with $\boldsymbol{v}^2{:=}\,\boldsymbol{v}\boldsymbol{v}{\,\ne}\,0$, then $\boldsymbol{v}$ has a multiplicative inverse\;$\displaystyle\boldsymbol{v}^{-1}{\,=\,}\frac{1}{\boldsymbol{v}^2}\boldsymbol{v}\,.$\\
(4) The geometric product, like the external product, is associative.
\subsection{}\textbf{Description of point reflections in a CK space using the geometric product.}\\
We start with a semi CK space $(\textrm{P}\boldsymbol{V}, \mathfrak{A})$ as introduced in Section 1. Let $Q=\mathbb{R}\boldsymbol{q}$ be any point and $R=\mathbb{R}\boldsymbol{r}$ and $S=\mathbb{R}\boldsymbol{s}$ be two anisotropic points in this space. The image of $Q$ under a reflection in $R$ is the point\vspace*{-0.5 mm}\\
$\displaystyle Q'= \mathbb{R}\big(\boldsymbol{q}-2\frac{\boldsymbol{q}{\scriptstyle{[\mathfrak{A}]}}\boldsymbol{r}}{\boldsymbol{r}{\scriptstyle{[\mathfrak{A}]}}\boldsymbol{r}}\boldsymbol{r}\big) = \mathbb{R}\big(\boldsymbol{q}-(\boldsymbol{q}\boldsymbol{r}+\boldsymbol{r}\boldsymbol{q})\frac{1}{\boldsymbol{r}{\scriptstyle{[\mathfrak{A}]}}\boldsymbol{r}}\boldsymbol{r}\big)
= \mathbb{R}\big(\boldsymbol{q}-\frac{\boldsymbol{q}\boldsymbol{r}^2}{\boldsymbol{r}{\scriptstyle{[\mathfrak{A}]}}\boldsymbol{r}}-\boldsymbol{r}\boldsymbol{q}\boldsymbol{r}^{-1}\big)=\mathbb{R}\big(\boldsymbol{r}\boldsymbol{q}\boldsymbol{r}^{-1}\big)$. \vspace*{0.5 mm}\\
The reflection of $Q'$ in $S$ delivers the point $Q''= \mathbb{R}\big(\boldsymbol{sr}\,\boldsymbol{q}\,\boldsymbol{(sr)}^{-1}\big)$.\vspace*{-3.0 mm}\\

We study the double reflection more closely and assume that $R\ne S$ and $\boldsymbol{r}^2, \boldsymbol{s}^2\in \{-1,1\}$. The antipode of $R$ on $R\,\sqcup\,S$ is a point $T=\mathbb{R}\boldsymbol{t}\in R\,\sqcup\,S$  
with $\boldsymbol{t}{\scriptstyle{\,[\mathfrak{A}]\,}}\boldsymbol{r}=0$.\\
If $R{\,\sqcup\,}S$ is isotropic, then $\boldsymbol{t}^2=0$ and we may assume that $\boldsymbol{r}^2=\boldsymbol{s}^2$. There exists one nonzero real number $d$ such that $\boldsymbol{s}=\boldsymbol{r}+d\boldsymbol{t}$, and we get
$\boldsymbol{s}\boldsymbol{r}=\boldsymbol{r}^2(1+d\boldsymbol{t}\boldsymbol{r}^{-1})=\boldsymbol{r}^2\exp(d\boldsymbol{t}\boldsymbol{r}^{-1})$\vspace*{-0.5 mm}\\ with  $\displaystyle\exp(\boldsymbol{v}):=\sum_{k=0}^{\infty}\frac{\boldsymbol{v}^k}{k!}$ \vspace*{-1.8 mm} {for} \vspace*{1 mm}{$\boldsymbol{v}\in \Lambda\boldsymbol{V}.$} \\
If $R\sqcup S$ is anisotropic, we may assume that $ \boldsymbol{t}^2 \in \{-1,1\}$. Let $d\in \mathbb{D}$ denote the distance of the points $R$ and $S$.\\
If $R\sqcup S$ is an elliptic line, $d$ is purely imaginary, thus $\tilde{d}:=-d\mathbf{i}\in\; ]\,0,\frac{1}{2}\pi\,[$. We have $\boldsymbol{r}^2=\boldsymbol{s}^2=\boldsymbol{t}^2=1$ and $\boldsymbol{s}\,=\,\cos(\tilde{d})\,\boldsymbol{r} + \sin(\tilde{d})\,\boldsymbol{t}$. Since $(\boldsymbol{t}\boldsymbol{r}^{-1})^2 = (\boldsymbol{t}\boldsymbol{r})^2 = -1$, we get \vspace*{0.5 mm}\\
\centerline{$\boldsymbol{s}\boldsymbol{r}=(\cos(\tilde{d})\,\boldsymbol{r} + \sin(\tilde{d})\,\boldsymbol{t})\,\boldsymbol{r}=\cos(\tilde{d})+\sin(\tilde{d})\,\boldsymbol{t}\boldsymbol{r}^{-1} = \exp(\,\tilde{d}\boldsymbol{t}\boldsymbol{r}^{-1}).$}\\
If $R\sqcup S$ is hyperbolic, the equations $\boldsymbol{t}^2= -\boldsymbol{r}^2$ and $\sqrt{\boldsymbol{s}^2}\,\boldsymbol{s}=\cosh(d)\sqrt{\boldsymbol{r}^2}\,\boldsymbol{r}-\sinh(d)\sqrt{-\boldsymbol{r}^2}\boldsymbol{t}\mathbf{i}$ apply. Since $(\boldsymbol{t}\boldsymbol{r}^{-1})^2 = (\boldsymbol{t}\boldsymbol{r})^2 = +1$, we get 
$\sqrt{\boldsymbol{s}^2}\,\boldsymbol{s}\boldsymbol{r}= \exp(\,{d}\boldsymbol{t}\boldsymbol{r}^{-1})$ in case of $\boldsymbol{r}^2 = 1$ and
$\sqrt{\boldsymbol{s}^2}\,\boldsymbol{s}\boldsymbol{r}= -\exp(\,- {d}\boldsymbol{t}\boldsymbol{r}^{-1})\,\mathbf{i}$ \,in case of $\boldsymbol{r}^2 = -1$.


\end{document}